\documentclass[10pt, a4paper, twoside]{article}

\usepackage{exscale, amsmath, amssymb, amsthm, dsfont, amsfonts, latexsym, url}
\usepackage[curve,matrix,arrow,ps]{xy}
\usepackage{hyperref}

\newtheorem{definition}{Definition}[section]
\newtheorem{theorem}[definition]{Theorem}
\newtheorem{proposition}[definition]{Proposition}
\newtheorem{lemma}[definition]{Lemma}
\newtheorem{corollary}[definition]{Corollary}
\newtheorem{conjecture}[definition]{Conjecture}

\newtheorem{deftheo}[definition]{Definition-Theorem}

\newcommand{\nd}{\noindent}

\newcommand{\dE}{{\mathds E}}
\newcommand{\dC}{{\mathds C}}

\newcommand{\dN}{{\mathds N}}

\newcommand{\dR}{{\mathds R}}

\newcommand{\dZ}{{\mathds Z}}
\newcommand{\dP}{{\mathds P}}

\newcommand{\dH}{{\mathds H}}

\newcommand{\bD}{{\mathbb D}}

\newcommand{\cA}{\mathcal{A}}
\newcommand{\cB}{\mathcal{B}}
\newcommand{\cC}{\mathcal{C}}
\newcommand{\cD}{\mathcal{D}}
\newcommand{\cE}{\mathcal{E}}
\newcommand{\cF}{\mathcal{F}}
\newcommand{\cG}{\mathcal{G}}
\newcommand{\cH}{\mathcal{H}}

\newcommand{\cK}{\mathcal{K}}
\newcommand{\cL}{\mathcal{L}}

\newcommand{\cN}{\mathcal{N}}
\newcommand{\cO}{\mathcal{O}}
\newcommand{\cP}{\mathcal{P}}
\newcommand{\cQ}{\mathcal{Q}}
\newcommand{\cR}{\mathcal{R}}
\newcommand{\cS}{\mathcal{S}}

\newcommand{\cU}{\mathcal{U}}
\newcommand{\cV}{\mathcal{V}}
\newcommand{\cW}{\mathcal{W}}

\newcommand{\SC}{\scriptstyle}

\newcommand{\HH}{{\mathcal H}}

\newcommand{\RR}{{\mathcal R}}

\newcommand{\UU}{{\mathcal U}}

\newcommand{\www}{\widetilde}
\newcommand{\oooo}{\overline}
\newcommand{\uuuu}{\underline}

\newcommand{\zdz}{{z\partial_z}}

\newcommand{\nnn}{\nabla}

\newcommand{\Cz}{\dC\{ z\}}

\newcommand{\Czm}{{\dC\{ z\}[z^{-1}]}}

\newcommand{\hiii}{H^\infty}
\newcommand{\hiir}{H^\infty_\dR}

\newcommand{\hiia}{H^\infty_{e^{-2\pi i\alpha}}}
\newcommand{\hiib}{H^\infty_{e^{-2\pi i\beta}}}

\newcommand{\Gr}{Gr}

\newcommand{\id}{\mathit{Id}}

\newcommand{\rank}{rank}

\begin{document}

\title{Nilpotent orbits of a generalization of Hodge structures}

\author{Claus Hertling\and Christian Sevenheck}

\date{March 23, 2006}

\maketitle

\begin{abstract}
\noindent
We study a generalization of Hodge structures which first appeared in the
work of Cecotti and Vafa. It consists of twistors,
that is, holomorphic vector bundles on $\dP^1$, with additional structure,
a flat connection on $\dC^*$, a real subbundle and a pairing.
We call these objects TERP-structures. We generalize to TERP-structures a correspondence of
Cattani, Kaplan and Schmid between nilpotent orbits of Hodge structures
and polarized mixed Hodge structures. The proofs use work of
Simpson and Mochizuki on variations of twistor structures
and a control of the Stokes structures of the poles at zero and infinity.
The results are applied to TERP-structures which arise
via oscillating integrals from holomorphic functions with
isolated singularities.

\end{abstract}
\renewcommand{\thefootnote}{}
\footnote{2000 \emph{Mathematics Subject Classification.}
14D07, 53C07, 32S40, 34M40\\
Keywords: mixed Hodge structures, nilpotent orbits, twistor structures,
$tt^*$ geometry, isomonodromic deformations, TERP-structures}

\tableofcontents

\setcounter{section}{0}

\section{Introduction}\label{c1}
\setcounter{equation}{0}

This paper studies generalizations of Hodge structures and
variations of them. They appeared in the work of
Cecotti and Vafa (\cite{CV1}\cite{CV2}) on
supersymmetric field theories. The abstract notion is of
general nature. It is studied under the name TERP-structure in
\cite{He2}.

In this paper we continue to investigate such TERP-structures. More specifically,
we study the relation between nilpotent orbits of them
and a corresponding generalization of mixed Hodge structures. This
result is an extension of a similar correspondence in Hodge theory
due to Cattani, Kaplan and Schmid.

Closely related objects, called twistor structures
(which also generalize Hodge structures)
appeared in the work of Simpson (\cite{Si5}).
TERP-structures are enriched twistors.
A twistor is simply a holomorphic vector bundle $\widehat H$ on $\dP^1$.
A twistor $\widehat H$ is called pure of weight $w\in\dZ$ if
it is semi-stable of slope $w$.
A pure twistor of weight 0 is polarized if there is a sesquilinear pairing
$\widehat S:\widehat{H}_z\times \widehat{H}_{-1/\oooo z}\to\dC$ for all $z\in\dP^1$
such that the induced pairing on $\Gamma(\dP^1,\widehat\cH)$
is hermitian and positive definite.
This can be generalized to any weight.
Note that a single twistor is a rather elementary object. However, if
one takes into account parameters, then the resulting structure is quite
involved, a variation of pure polarized twistors is actually equivalent
to a harmonic bundle on the parameter space (see \cite{Si1}\cite{Si2}\cite{Si4}).

A TERP-structure (the name stems from ``T'' for twistor, ``E'' for extension,
``R'' for real structure and ``P'' for pairing) is a twistor with additional data
which generalize all ingredients of a polarized Hodge structure, that is,
a real structure, a pairing and the Hodge filtration.
More precisely, a TERP-structure of weight 0 is a holomorphic vector bundle $H$ on $\dC$
with a flat connection $\nnn$ on $H_{|\dC^*}$ which has a pole of
order at most two at 0.
The real structure consists of a flat
real subbundle $H'_{\dR}$ on $\dC^*$. The pairing in this case is a pairing
$P:H_z\times H_{-z}\to \dC$ for $z\in\dC$,
which is symmetric and nondegenerate on $\dC$ and flat on $\dC^*$ and
which takes values in $\dR$ on the real subbundle.
In many applications it makes sense to consider the data $(H_{|\dC^*},\nnn,H'_{|\dR},P)$
as ``topological objects'' and the extension of the bundle to 0
as transcendent. This extension is the generalization of the Hodge
filtration. A key point now is that real and flat structures allow
to construct canonically an extension to infinity making up a twistor
(see chapter \ref{c3} for the precise construction).
The extension at $\infty$ then generalizes naturally the complex
conjugate of the Hodge filtration. By construction, the resulting twistor comes
equipped with a meromorphic connection with poles of order at most two
at zero and infinity. Moreover, the given pairing $P$
and the construction of the extension yield a pairing $\widehat{S}$
like the one described above.

A TERP-structure $H$ will be called pure (resp. polarized pure) if
the corresponding object $(\widehat{H},\widehat{S})$ is a pure
(resp. polarized pure) twistor. Pureness generalizes
the notion of opposite filtrations, i.e. of a Hodge structure.
Polarized pure TERP-structures generalize
polarized pure Hodge structures. The precise relation between
TERP-structures and twistor structures is reviewed in chapter \ref{c3}.
We also make some comments on the parameter case, i.e., variations
of TERP- resp. twistor structures, which were treated in great detail
in \cite{He2}.

\medskip

Important tools in the study of variations of Hodge structures are
the notions of nilpotent orbits of Hodge structures and mixed
Hodge structures, the first one being geometric, the second linear. There is a beautiful correspondence
relating them (theorem \ref{theoSchmidCorres}). One direction is due to
Schmid \cite{Sch}, the other one has been treated later in a series of papers
(\cite{CaK1}\cite{CaKS}\cite{CaK2}) by Cattani, Kaplan and Schmid. We give
a short reminder of these notions and results in chapter \ref{c2}.

The main purpose of the present paper is the generalization of this correspondence
to TERP-structures. The notion of a nilpotent orbit has a rather simple generalization
(definition \ref{t4.1}): We say that a
TERP-structure $(H,\nnn,H'_\dR,P)$ induces a nilpotent orbit
if $\pi^*_r(H,\nnn,H'_\dR,P)$ is a polarized pure TERP-structure
for any $r\in\dC^*$ with $|r|\ll 1$; here $\pi_r:\dC\to\dC,z\mapsto z\cdot r$, for $r\in\dC^*$.
It comes as a surprise that a simple rescaling of the
coordinate on $\dC$ changes the twistor which results
from the gluing procedure in an essential way.

The generalization of a polarized mixed Hodge structure is what we call a mixed TERP-structure.
Its definition (definition \ref{main-theorem}) is rather involved,
as we have to deal with a possibly irregular singularity of the connection $\nabla$ at zero.
Let us first explain the regular singular case. A TERP-structure $(H,\nnn,H'_\dR,P)$ is called regular singular
if the pole at 0 is so.
In that case there is a well known procedure to obtain a filtration
on the space $\hiii$ of multivalued flat global sections of the
bundle $H_{|\dC^*}$. Up to a twist, it was first considered
by Varchenko \cite{Va1} in the context of hypersurface singularities,
and refined later in \cite{SchSt}\cite{SM}.
For a regular singular TERP-structure, the condition to be mixed simply
means that this filtration is part of a polarized mixed Hodge structure.
If, however, the pole at zero becomes irregular, then we have to modify this condition.
First, we require that the formal decomposition of $(H,\nabla)$ can be done
without ramification. Moreover, we need
a compatibility condition between $H'_\dR$ and the Stokes structure defined by the irregular
pole. Under these two hypotheses, the regular singular factors which appear in the
formal decomposition are themselves TERP-structures and the main
condition imposed is that they induces polarized mixed Hodge structures as before.
The details are explained in the chapters \ref{c8} and \ref{c9}.
With this notion in mind, we can state the generalized correspondence as follows.
\begin{conjecture} (conjecture \ref{main-conjecture})
A TERP-structure which does not require a ramification
is a mixed TERP-structure iff it induces a nilpotent orbit.
\end{conjecture}

The main result of this paper is a proof of a good part of this conjecture, namely:
\begin{theorem}(theorem \ref{main-theorem})
\begin{enumerate}
\item The conjecture is true if the TERP-structure is regular singular.
\item The implication $\Rightarrow$ is true for any TERP-structure.
\end{enumerate}
\end{theorem}
The implication $\Rightarrow$ in the regular singular case
was shown in \cite[theorem 7.20]{He2}, using the analogous
implication in the correspondence between polarized mixed Hodge structures and nilpotent orbits
of Hodge structures. The opposite implication $\Leftarrow$ for regular singular TERP-structures
is proved in chapter \ref{c6} and
uses quite different techniques, namely, it relies on a recent result of
Mochizuki \cite{Mo2}: In that paper, he constructs, among other
things, for any tame harmonic bundle on the complement of a normal crossing divisor
a limit polarized mixed twistor structure.
This generalizes Schmid's limit polarized mixed Hodge structure
defined by a variation of Hodge structures.

For nilpotent orbits of TERP-structures, we only need the one-variable version
of the limit mixed twistor structure (\cite[theorem 12.1]{Mo2}).
To apply it, we establish in chapter \ref{c5}
a correspondence (lemma \ref{lemCorPMHS-PMTS})
between certain integrable polarized mixed twistor structures and polarized mixed Hodge structures
equipped with a semi-simple automorphism.
This correspondence extends similar correspondences in \cite[ch. 3]{Mo2}.

The last part of the paper (chapters \ref{c8} to \ref{c10}) deals with the general (i.e., irregular) case.
The implication $\Rightarrow$ in general is proved in chapter
\ref{c9}. It combines the regular singular case with a discussion of the Stokes structure.
In the end it comes down to a Riemann boundary value problem which we are able to solve
by an argument involving the Birkhoff decomposition
of the loop group $\Lambda GL_n(\dC)$ \cite[(8.1.2)]{PS}.

A particular case arises if the pole part of the connection
is a semi-simple endomorphism with pairwise different
eigenvalues. Such TERP-structures are called semi-simple and
the implication $\Rightarrow$
was already established by Dubrovin \cite[proposition 2.2]{Du} in that case.

\medskip

For semi-simple TERP-structures of rank two, families
$\bigcup_{r>0}\pi^*_{r^{-1}}(H,\nnn,H'_\dR,P)$ of TERP-structures
are closely related to solutions of the sinh-Gordon equation
$(\partial_r^2+\frac{1}{r}\partial_r)u(r)=\sinh u(r)$
\cite{CV1}\cite{CV2}\cite{Du} (and implicitly also \cite{IN}).
The implication $\Leftarrow$ in the semi-simple rank two case
is equivalent to the claim that the only solutions $u(r)$
which are smooth and real for large $r$ are the one parameter
family of solutions which were studied in \cite{MTW}.
This is very probably true, but the statements in \cite{IN}\cite{MTW}
do not imply it immediately. See the discussion at the end of chapter \ref{c10}
for more details.
The conjecture in general predicts that the singularity freeness
for large $r$ of certain systems of differential equations is equivalent to
``linear'' conditions which are neatly formulated
in terms of the Stokes data.

The semi-simple case is also interesting because of the apparent
simplicity in which the structure can be encoded.

\begin{lemma} (lemma \ref{t10.1}) \label{lemSemiSimpleUniv-Int}
Any data $(w, u_1,...,u_n,\xi,T)$ where $w\in\dZ$,
$u_i\in\dC$ with $u_i\neq u_j$ for $i\neq j$,
$\xi\in S^1$ with $\Re(\frac{u_i-u_j}{\xi})<0$ for $i<j$
and $T\in M(n\times n,\dR)$ upper triangular with $T_{ii}=1$ give
rise to a unique semi-simple mixed TERP-structure of weight $w$.
The numbers $u_1,...,u_n$ are the eigenvalues of the pole part of the
connection and the Stokes structure is equivalent to the matrix $T$ together with the point $\xi$.
Any semi-simple mixed TERP-structure arises in such a way (the choice of $\xi$ is not unique
and influences $T$).
\end{lemma}

An interesting question is to know which of these TERP-structures are pure and
polarized. The following conjecture proposes a partial answer.

\begin{conjecture}(see conjecture \ref{t10.2})
If the matrix $T+T^{tr}$ is positive definite,
then the TERP-structure is pure and polarized.
\end{conjecture}

The conjecture is proved in chapter \ref{c10} for Stokes data
which are related to the ADE-singularities.
One ingredient are beautiful results in \cite{De2} and \cite{Lo}
which roughly say that the semi-universal unfolding
of an ADE-singularity is a one to one atlas for all distinguished bases
up to signs and a finite to one atlas for possible Stokes data.
Another ingredient is a recent fundamental result of Sabbah
\cite[theorem 4.9]{Sa8} on TERP-structures for tame functions,
which applies to the ADE-singularities and their unfoldings.

Using oscillating integrals resp. by the Fourier-Laplace transform
of its Gauss-Manin system, any function germ $f:(\dC^w,0)\to (\dC,0)$
with an isolated singularity at zero and similarly any tame function $f:Y\to\dC$ on an affine manifold $Y$
gives rise to a mixed TERP-structure (theorem \ref{t11.1}).
This result relies on the work of many different people.
In the tame case, it has been established (although it is not
expressed in these terms) in \cite{Sa2} and \cite{DS}.
For function germs, the construction is described in \cite[8.1]{He2}.

Sabbah's new result \cite[theorem 4.9]{Sa8}
(see theorem \ref{t11.2}) states that for tame functions,
this TERP-structure is always pure and polarized. It seems
to be of fundamental importance for future study of tame functions.
In a sense, it should be seen as the analogue
to the fact that the primitive part of the
cohomology of a compact K\"ahler manifold
is a sum of polarized pure Hodge structures.
Even though this analogy is only a sort of
meta-theorem, it can hopefully be turned into more concrete results
using the point of view of mirror symmetry:
certain tame functions on affine manifolds
correspond to certain Fano manifolds.
We give some speculations in that direction in chapter \ref{c11}.

TERP-structures for tame functions, and even Sabbah's theorem,
are treated from the physicists' perspective in \cite{CV1}\cite{CV2}.
The functions are part of Landau-Ginzburg models, the positive
definite hermitian metric of the polarized pure TERP-structure
is a ground state metric. Also lemma \ref{lemSemiSimpleUniv-Int} is
already used implicitly in \cite{CV1}\cite{CV2}. Semi-simple
TERP-structures arises from massive field theories,
and it is an important feature in cit.loc. that these can be encoded
by simple data as in lemma \ref{lemSemiSimpleUniv-Int}.

It is an elementary computation that the TERP-structure $\mathit{TERP}(f)$ of a function $f$
(a germ or tame) satisfies $\mathit{TERP}(r\cdot f)=\pi_{r^{-1}}^*(\mathit{TERP}(f))$ for any parameter $r\in\dC^*$.
The fact that $\mathit{TERP}(f)$ is mixed and theorem \ref{main-theorem} 2.
imply that $\mathit{TERP}(r\cdot f)$ is pure and polarized for $|r|\gg 0$.
This proves the main part of the conjecture 8.3 in \cite{He2}.
It endows a part of the semi-universal unfolding space with a positive
definite hermitian metric. This will hopefully have applications
for moduli space questions or Torelli problems.

In \cite{CV1}\cite{CV2} the limit $r\to \infty$ is called infrared limit,
the limit $r\to 0$ is called ultraviolet limit. We also consider
$r\to 0$ and define a counterpart of a nilpotent orbit which is called Sabbah orbit
(definition \ref{t4.1}). The reason for this is that Sabbah
has defined in \cite{Sa2} for any (not necessarily regular singular) TERP-structure
a filtration $F^\bullet_{Sab}$ on $\hiii$.
In particular, it was proved that this filtration gives
a mixed Hodge structure for TERP-structures coming from tame functions.
However, polarizations were not considered in that paper.
The ``Sabbah orbit''-version of our correspondence reads as follows.
\begin{theorem} (see theorem \ref{theoMainResultSabbah})\label{theoMainResultSabbah-Int}
A TERP-structure induces a Sabbah orbit iff a twisted version of the filtration
$F^\bullet_{Sab}$ gives rise to a polarized mixed Hodge structure.
\end{theorem}
The implication $\Leftarrow$ is analogous to \cite[theorem 7.20]{He2},
the implication $\Rightarrow$ uses again \cite[theorem 12.1]{Mo2}.
This theorem \ref{theoMainResultSabbah-Int} and \cite[theorem 4.9]{Sa8}
show that $F^\bullet_{Sab}$, twisted appropriately, even makes up
a polarized mixed Hodge structure.

Besides \cite{CV1}\cite{CV2}, this paper owes a lot to \cite{Sa6} and \cite{Mo2}.
In \cite{Sa6}, Simpson's notion of a variation of twistor structures
is generalized to polarizable twistor $\cD$-modules.
Chapter 7 of cit.loc. is devoted to twistor structures with connection on $\dC^*$
and a flat hermitian pairing $\widehat S$ as above
(one might call them ``TEH''-structures). The extraordinary long paper
\cite{Mo2} provided us with a central degeneration result used
in the proofs of theorem \ref{theoMainResultRegSing} and theorem
\ref{theoMainResultSabbah}. In the third chapter of cit.loc., nilpotent orbits
of polarized twistor structures are considered, however,
these twistors are not equipped with connections on $\dC^*$,
so they are less rich than TERP-structures.

\bigskip

\textbf{Acknowledgements:}
We would like to thank Claude Sabbah for many helpful discussions.

\bigskip

\textbf{Notations:} For any complex manifold $M$, we denote
by $\overline{M}$ the same real manifold, equipped with the
conjugate complex structure, i.e., $\cO_{\overline{M}}:=\overline{\cO}_M$.
We will need at different places a total ordering on $\dC$ extending
the usual real ordering of $\dR$. This will be the lexicographic one, that
is, we will write $\alpha < \beta$ iff
either $\Re(\alpha)<\Re(\beta)$ or
$\Re(\alpha)=\Re(\beta)$ and $\Im(\alpha)<\Im(\beta)$.
However, the interval notation $[\alpha,\alpha']$ for
$\alpha, \alpha'\in\dR$ will continue to denote
all $\beta\in\dR$ with $\alpha\leq\beta\leq\alpha'$.
If $H$ is a holomorphic vector bundle on a complex manifold $M$, we
write $H\in\mathit{VB}_M$ for short and we will use the symbol
$\cH$ to denote its sheaf of holomorphic sections.
In \cite{He2} pure TERP-structures and their variations were
called (tr.TERP)-structures, polarized pure TERP-structures
and their variations were called (pos.def.tr.TERP)-structures.
The following maps will be used frequently in the paper.
$$
\begin{array}{c}
i:\dC^*\hookrightarrow\dC
\;\; ; \;\;\;\;
\widetilde{i}:\dC^*\hookrightarrow\dP^1\backslash\{0\}
\;\; ; \;\;\;\;
\widehat{i}:\dC^*\hookrightarrow\dP^1
\;\; ; \;\;\;\;
\pi_r:\dC\to\dC\;,\;\; z\mapsto r\cdot z \textup{ for }r\in\dC^*;
\\ \\
j:\dP^1\rightarrow \dP^1\;,\;\; j(z)=-z
\;\; ; \;\;\;\;
\gamma:\dP^1\rightarrow \dP^1\;,\;\; \gamma(z)=\overline{z^{-1}}
\;\; ; \;\;\;\;
\sigma:\dP^1\rightarrow \dP^1\;,\;\; \sigma(z)=-\overline{z^{-1}}.
\end{array}
$$

\section[Polarized mixed Hodge structures]
{Polarized mixed Hodge structures and
nilpotent orbits of Hodge structures}\label{c2}

This chapter recalls for the reader's
convenience some classical notions
from Hodge theory. We give the definition
of nilpotent orbits and state the correspondence
between nilpotent orbits of polarized
Hodge structures and (limit) polarized mixed Hodge structures.
This correspondence is one of the main motivation of
our work.

Throughout the whole chapter, $w$ will be an integer, $H$ a complex vector
space of finite dimension, $H_\dR$ a real subspace with
$H=H_\dR\oplus iH_\dR$, and $S$ a nondegenerate $(-1)^w$-symmetric pairing
on $H$ with real values on $H_\dR$.

\begin{definition}\label{defPHS}
A polarized Hodge structure of weight $w\in\dZ$
(abbreviation PHS) consists
of data $H,H_\dR$ and $S$ as above and an exhaustive decreasing
{\it Hodge filtration} $F^\bullet$ on $H$ with the following
properties.
\begin{eqnarray}
 F^p \oplus \oooo{F^{w+1-p}}&=&H, \label{2.1}\\
 S(F^p,F^{w+1-p})&=&0, \\
 i^{p-(w-p)}\cdot S(a,\oooo a)&>&0 \quad
\textup{ for }a\in \left(F^p\cap F^{w+1-p}\right) \backslash \{0\} .
\label{2.3}
\end{eqnarray}
\end{definition}

\textbf{Remarks:}
\begin{enumerate}
\item
A tuple $(H,H_\dR,F^\bullet,w)$ with \eqref{2.1} is a pure Hodge structure
of weight $w$.
\item
The filtrations $F^\bullet$ and $\oooo{F^{w-\bullet}}$ are called
opposite if they satisfy \eqref{2.1}. This condition is equivalent
to the Hodge decomposition
$H=\bigoplus_p H^{p,w-p}$, where
$H^{p,w-p} := F^p\cap \oooo{F^{w-p}}$.
Of course, $\oooo{H^{p,q}}=H^{q,p}$.
\item
Given a PHS, the Hodge subspaces $H^{p,w-p}$
are orthogonal with respect to the pairing $S(\;\;,\overline{\;\cdot\;})$
The pairing
\begin{eqnarray*}
h:H\times H&\to& \dC, \\
(a,b)&\mapsto& i^{p-(w-p)}\cdot S(a,\oooo b) \quad
\textup{for }a\in H^{p,w-p},b\in H\nonumber
\end{eqnarray*}
is hermitian and positive definite. This pairing distinguishes
a PHS from a Hodge structure.
\end{enumerate}

Mixed Hodge structures have been introduced by Deligne \cite{De} in order
to study the cohomology of singular or non-compact Kähler manifolds.
A MHS contains a second filtration called weight filtration such that
the Hodge filtration induces pure Hodge structures of appropriate weights on
the graded spaces with respect to the weight filtration. In Schmid's work,
a more specific variant of mixed Hodge structures is considered: The weight filtration
is always induced by a given nilpotent endomorphism on $H_\dR$ satisfying
\begin{eqnarray}\label{eqNisInfIso}
S(N\, a,b)+S(a,N\, b)=0,
\end{eqnarray}
i.e., which is an infinitesimal isometry of $S$. These data yield
a weight filtration $W_\bullet$ in the following way.
\begin{lemma} \label{defWeigthFilt} \cite[Lemma 6.4]{Sch}
Let $(H,H_\dR,S,N,w)$ be as above.
\begin{enumerate}
\item
There exists a unique exhaustive increasing filtration $W_\bullet$ on
$H_\dR$ such that $N(W_l)\subset W_{l-2}$ and such that
$N^l:\Gr_{w+l}^W\to \Gr_{w-l}^W$ is an isomorphism.
\item
The filtration satisfies
$S(W_l,W_{l'})=0$ for $l+l'<w$.
\item
A nondegenerate $(-1)^{w+l}$-symmetric bilinear form $S_l$ is well
defined on $\Gr_{w+l}^W$ for $l\geq 0$ by
$S_l(a,b):=S(\www a,N^l \www b)$ for $a,b\in \Gr_{w+l}^W$
with representatives $\www a,\www b\in W_{w+l}$.
\item
The primitive subspace $P_{w+l}\subset \Gr_{w+l}^W$ is defined by
$$
P_{w+l}:= \ker (N^{l+1}:\Gr_{w+l}^W\to \Gr_{w-l-2}^W)
$$
for $l\geq 0$ and by $P_{w+l}:=0$ for $l<0$. Then
\begin{eqnarray}\label{eqDecomPrimMHS}
\Gr_{w+l}^W=\bigoplus_{i\geq 0}N^iP_{w+l+2i},
\end{eqnarray}
and this decomposition is orthogonal with respect to $S_l$ if $l\geq 0$.
\end{enumerate}
\end{lemma}

\begin{definition}\label{defPMHS} \cite{CaK1}\cite{He1}
A {\it polarized mixed Hodge structure of weight $w$} (abbreviation PMHS)
consists of data $H,H_\dR,S,N$ and $W_\bullet$ as above and an exhaustive
decreasing {\it Hodge filtration} $F^\bullet$ on $H$ with the following
properties.
\begin{enumerate}
\item
The filtration $F^\bullet\Gr_k^W$ on $\Gr^W_k$
gives a pure Hodge structure of weight $k$,
\item
$N$ is a $(-1,-1)$-morphism of mixed Hodge structures,
\begin{eqnarray}\label{eqNisStrict}
N(F^p)\subset F^{p-1},
\end{eqnarray}
\item
\begin{eqnarray}\label{eqFisNOrthPMHS}
S(F^p,F^{w+1-p})=0.
\end{eqnarray}
\item
For $a\in \left(F^pP_{w+l}\cap \oooo{F^{w+l-p}P_{w+l}}\right) \backslash \{0\}$
\begin{eqnarray}
i^{p-(w+l-p)}S_l(a,\oooo a)>0.\label{eqPolPHMS}
\end{eqnarray}
\end{enumerate}
\end{definition}
\textbf{Remark:}
The conditions \eqref{eqNisInfIso}, \eqref{eqNisStrict}
and \eqref{eqFisNOrthPMHS} imply
that $
S_l(F^pP_{w+l},F^{w+l+1-p}P_{w+l})=0
$. This condition and condition \eqref{eqPolPHMS} say
the pure Hodge structure $F^\bullet P_{w+l}$ of weight $w+l$ on
$P_{w+l}$ is polarized by $S_l$.

\bigskip

Let us fix one reference polarized Hodge structure
$(H,H_\dR,S,F^\bullet_0)$ of weight $w$. The space
$$
\check{D} := \left\{\textup{filtrations }F^\bullet \mbox { on } H \; |\;
\dim F^p=\dim F^p_0,\;S(F^p,F^{w+1-p})=0\right\}
$$
is a closed submanifold of a product of Grassmannians, in particular
projective. It is also a complex homogeneous space. Consider the subspace
$$
D:= \left\{ F^\bullet\in \check{D} \; | \; F^\bullet
\textup{ gives rise to a PHS, i.e.,
satisfies \eqref{2.1} and \eqref{2.3}}\right\}
$$
which is an open complex submanifold and a real homogeneous space \cite{Sch}.
It classifies polarized Hodge structures with fixed Hodge numbers.

\begin{definition}\label{defNilOrbHS}
A tuple $(H,H_\dR,S,F^\bullet,N)$
is said {\it to give rise to a nilpotent orbit} if the following holds:
\begin{enumerate}
\item
$F^\bullet \in \check{D}$,
\item
the endomorphism $N$ of $H_\dR$ is nilpotent
and an infinitesimal isometry with $N(F^p)\subset F^{p-1}$,
\item
there exists a bound $b\in \dR$ such that
$$
e^{\rho N}F^\bullet \in D \textup{ for \ } \Im(\rho) > b .
$$
\end{enumerate}
Then the set $\{e^{\rho N}F^\bullet\ |\ \rho\in \dC\}$ is called {\it a nilpotent
orbit of Hodge structures}.
\end{definition}

Nilpotent orbits of Hodge structures
play a fundamental role in Schmid's work \cite{Sch}.
The following theorem gives a beautiful correspondence between
PMHS and nilpotent orbits of Hodge structures.
The main purpose of the whole paper is to generalize this
correspondence to TERP-structures.

\begin{theorem}\label{theoSchmidCorres}
\cite{Sch}\cite{CaK1}\cite{CaKS}\cite{CaK2}
Let $(H,H_\dR,S)$ be as above.
\begin{enumerate}
\item
The tuple $(H,H_\dR,S,F^\bullet,N)$ is a PMHS of weight $w$
if and only if $(H,H_\dR,F^\bullet,N)$
gives rise to a nilpotent orbit of Hodge structures.
\item
If $(H,H_\dR,S,F^\bullet,N)$ is a PMHS with $I^{q,p}=\oooo{I^{p,q}}$
then $e^{\rho N}F^\bullet \in D$ for $\Im (\rho)>0$.
\end{enumerate}
\end{theorem}
The direction `$\Longleftarrow$' in \textit{1.} is shown in \cite[Theorem 16.6.]{Sch}.
It is a consequence of the $SL_2$-orbit theorem.
`$\Longrightarrow$' in \textit{1.} is \cite[Corollary 3.13]{CaKS}.
Short proofs of both directions are given in \cite[Theorem 3.13]{CaK2}.
The special case \textit{2.} is proved in \cite[Proposition 2.18]{CaK1}
and in \cite[Lemma 3.12]{CaKS}.

The nilpotent orbit theorem \cite{Sch} says that any variation of PHS
on a punctured disk is approximated by a nilpotent orbit.
Therefore the correspondence above associates a (limit) PMHS
to any variation of PHS on a punctured disk.

\section{TERP-structures and twistor structures}\label{c3}
\setcounter{equation}{0}

This chapter introduces the central objects of this paper:
TERP-structures. This notion encapsulates a situation encountered
when studying Hodge theory for singularities. More precisely,
a TERP-structure arises when performing
a Fourier-Laplace transformation of the Gauss-Manin-system and the Brieskorn lattice
of a holomorphic function germ or a tame polynomial.
In a sense which will become clear later (chapter \ref{c6}),
(variations of) TERP-structures are natural generalizations of
(variations of) Hodge structures.
We will give in this chapter the definitions and some properties
of TERP-structures. The main point is the construction of
a bundle on $\dP^1$ starting from a given TERP-structure. After recalling
the notion of (polarized) twistor structure, we will see that this $\dP^1$-bundle
is a (polarized) twistor with some additional structure, called integrable twistor.
We only make some comments on how to extend these constructions to the case with
parameters, i.e., for variations of TERP/twistor structures. In a sense, this
chapter is a short version of the second chapter of \cite{He2}
with some additional notations and comparison results.

\begin{definition}\label{defTERP}
A TERP-structure (``twistor, extension, real structure, pairing'')
of weight $w\in\dZ$
is a tuple  $(H,H'_\dR, \nabla, P, w)$ where $H$ is a
holomorphic vector bundle on $\dC$,
equipped with a flat meromorphic connection $\nabla$ with a pole of order at most two at zero,
a flat real subbundle $H'_\dR \subset H':=H_{|\dC^*}$ of the restriction to $\dC^*$
satisfying $H'=H'_\dR\otimes\dC$ and a flat, bilinear, $(-1)^w$-symmetric,
nondegenerate pairing
$$
P:H_z \times H_{-z} \longrightarrow \dC \quad \textup{ for }z\in\dC^*
$$
with the following two properties.
\begin{enumerate}
\item
For any $z\in \dC^*$, we have
\begin{equation}\label{eqPReal}
P:(H'_\dR)_z\times (H'_\dR)_{-z}\rightarrow i^w\dR.
\end{equation}
\item
The pairing induced on sections satisfies
\begin{equation}\label{eqPPoleZero}
P:\cH \otimes j^*\cH \longrightarrow z^w\cO_\dC,
\end{equation}
and the pairing $z^{-w}P$ is nondegenerate at $0$.
\end{enumerate}
\end{definition}
\begin{definition}[extension to infinity]
\label{defHhat}
Consider a TERP-structure $(H,H'_\dR,\nabla,P,w)$. Let $\gamma:\dP^1\rightarrow\dP^1; z\mapsto \overline{z}^{-1}$
and define for any $z\in\dC^*$ the following two anti-linear involutions.
$$
\begin{array}{c}
\begin{array}{rcl}
\tau_{real}: H_z & \longrightarrow & H_{\gamma(z)}\\
s & \longmapsto & \nabla\textup{-parallel transport of }\overline{s}\\ \\
\tau: H_z & \longrightarrow & H_{\gamma(z)}\\
s & \longmapsto & \nabla\textup{-parallel transport  of }\overline{z^{-w}s}
\end{array}
\end{array}
$$
The induced maps on sections by putting
$s\mapsto\left(z\mapsto\tau s(\overline{z}^{-1})\right)$ resp.
$s\mapsto\left(z\mapsto\tau_{real}s(\overline{z}^{-1})\right)$
will be denoted by the same letter.
They can either be seen as morphisms
$\tau,\tau_{real} :\cH' \rightarrow \overline{\gamma^*\cH'}$
which fix the base,
or as morphisms $\tau,\tau_{real}:\cH'\rightarrow\cH'$
which map sections in $U\subset \dC^*$ to
sections in $\gamma(U)\subset \dC^*$.
Note that due to the two-fold conjugation (in the base and in fibres),
$\tau$ and $\tau_{real}$ are morphisms of holomorphic bundles over $\dC^*$.
Denote by $\widehat{H}\in\mathit{VB}_{\dP^1}$ the bundle
obtained by patching $\cH$ and $\overline{\gamma^*\cH}$ via the identification
$\tau$.
\end{definition}
Notice that the pairing $P$ does not enter
in the construction of the bundle $\widehat{H}$. However, the sole fact
that the bundle $H$ is equipped with a pairing with the above properties
puts restrictions on $\widehat{H}$ as the following lemma shows.
Let us denote the sheaf $\cO(\widehat{H}_{|\dP^1\backslash\{0\}})$ by
$\widetilde{\cH}$ for short.

\begin{lemma}\label{lemDegHzero}
\begin{enumerate}
\item
The connection naturally extends
with a pole of order two at infinity.
\item
The pairing $P$ satisfies
$P:\widetilde{\cH}\otimes j^*\widetilde{\cH}
\rightarrow z^w\cO_{\dP^1\backslash\{0\}}$,
and $z^{-w}P$ is nondegenerate at $\infty$.
\item The bundle $\widehat{H}$ has degree zero.
\end{enumerate}
\end{lemma}
\begin{proof}
We will need the following two equalities of
endomorphisms of $\cH'$, which express the flatness property of
$\tau_{real}$.
They are immediate consequences of
$\overline{\gamma}^*\left(\frac{dz}{z}\right)=-\frac{dz}{z}$.
\begin{equation}\label{eqTauFlat}
\nabla_{z\partial_z} \circ \tau_{real}=
\tau_{real}\circ \nabla_{-z\partial_z}
\;\;\; ; \;\;\;
\nabla_{z\partial_z} \circ \tau=
\tau\circ\left(\nabla_{-z\partial_z}+w\mathit{Id}\right) .
\end{equation}

Consider $\cH$ as a subsheaf of $i_*\cH'$ where
$i:\dC^*\hookrightarrow \dC$. By definition, we have
$\widetilde{\cH}=\tau\cH\subset \widetilde{i}_*\cH$ with
$\widetilde{i}:\dC^*\hookrightarrow \dP^1\backslash\{0\}$.
This gives immediately, using the above formula, that
$$z^{-1}\nabla_{\partial_{z^{-1}}}\widetilde{\cH}
=-\nabla_{\zdz}\tau\cH
=\tau(\nnn_\zdz -w\mathit{Id})\cH = \tau(\frac{1}{z}\cH)=z\www\cH.$$
To detect the order of $P$ at infinity, consider the following calculation:
\begin{equation}\label{eqPPoleInfty}
\begin{array}{rccl}
\overline{z^{-w}\cdot P(a(z),b(-z))} & = &
(-\overline{z})^w\cdot \overline{P(z^{-w}a(z),(-z)^{-w}b(-z))} \\ \\
& = & \overline{z}^w\cdot P(\overline{z^{-w}a(z)},\overline{(-z)^{-w}b(-z)})
& [(-1)^w\mbox{ because of condition \eqref{eqPReal}}]
\\ \\
&=& \overline{z}^w\cdot P\left(\tau(a)(\overline{z}^{-1}),\tau(b)(-\overline{z}^{-1})\right)
& [P\textup{ is flat}].
\end{array}
\end{equation}
The order of $P$ at zero and $\widetilde{\cH}=\tau\cH$ yield $P:\widetilde{\cH}\otimes
j^*\widetilde{\cH}\rightarrow z^w\cO_{\dP^1\backslash\{0\}}$ as required.
To prove the third point, we first consider the case where $\rank(H)=1$. Then $\widehat{\cH}\cong\cO_{\dP^1}(k)$
for some $k\in\dZ$. Choose non-vanishing holomorphic sections
$\sigma\in\Gamma(\dC,\widehat{\cH})$ and $\widetilde{\sigma}\in\Gamma(\dP^1\backslash\{0\},\widehat{\cH})$
satisfying $\widetilde{\sigma}=z^k\sigma$ on $\dC^*$. This implies that
$z\mapsto P(\sigma(z),\sigma(-z))$ defines a non-vanishing holomorphic function on $\dC^*$
with a zero of order $w$ at zero (by equation \eqref{eqPPoleZero} in the definition of TERP-structures) and a zero of order
$2k-w$ at infinity by the above computation. Consequently, $k=0$. Now for the general case, we remark that
given any TERP-structure $(H,H'_\dR,\nabla,P,w)$, then the determinant (line) bundle $\det(H)$
is naturally a TERP-structure of weight $w\cdot\rank(H)$. Moreover,
it is clear that $\widehat{\det(H)}=\det(\widehat{H})$. This implies that $\deg(\widehat{H})=0$.
\end{proof}

The next step is to investigate more closely the case where $\widehat{H}$ is a
trivial $\dP^1$-bundle. This implies that we have a canonical identification
of any fibre with the space $H^0(\dP^1,\widehat{\cH})$ and thus also
a canonical identification of all fibres.
\begin{lemma}\label{lemhIsHerm}
Let $\widehat{H}$ be trivial and consider the identification
$H_0 \stackrel{\cong}{\longrightarrow} H^0(\dP^1, \widehat{\cH})$.
The morphism $\tau$ acts on this space as an anti-linear involution
and the pairing $z^{-w}P$ is symmetric and has constant values on it.
Define
$$
h:H_0\times H_0  \longrightarrow  \dC\;\;\;;\;\;\;(a,b)
\longmapsto  z^{-w}P(a,\tau b) .
$$
Then $h$ is a hermitian pairing on $H_0$.
\end{lemma}
\begin{proof}
In order to see that $\tau$ defines an anti-holomorphic involution
on $H^0(\dP^1,\widehat{\cH})$, consider the extension $\widehat{i}_*\cH'$,
where $\widehat{i}:\dC^*\hookrightarrow \dP^1$. We have
$\widetilde{i}_*\cH\subset\widehat{i}_*\cH'$, $i_*\widetilde{\cH}\subset\widehat{i}_*\cH'$,
and $H^0(\dP^1, \widetilde{i}_*\cH\cap i_*\widetilde{\cH})$ is precisely
the finite-dimensional space $H^0(\dP^1,\widehat{\cH})$ which is of dimension
equal to the rank of $H$ if $\widehat{H}$ is trivial. The morphism $\tau$
acts on $\cH'$ and therefore on $\widehat{i}_*\cH'$. It maps
$\widetilde{i}_*\cH$ isomorphically to $i_*\widetilde{\cH}$ and vice versa.
This shows that it acts on $H^0(\dP^1,\widehat{\cH})$.

For any two global sections $a,b\in H^0(\dP^1,\widehat{\cH})$, putting
$z\mapsto P\left(a(z),b(-z)\right)$ defines a holomorphic function with zero of order
$w$ at the origin and pole of order $w$ at infinity. Therefore, the function
$z\mapsto z^{-w}P\left(a(z),b(-z)\right)$ is holomorphic on $\dP^1$ and thus constant.
The symmetry property follows from
$$
z^{-w}P\left(a,b\right)(z)=z^{-w}P\left(a(z),b(-z)\right)
=z^{-w}(-1)^wP\left(b(-z),a(z)\right)=z^{-w}P\left(b,a\right)(z).
$$
In order to show that $h$ is hermitian, we apply computation
\eqref{eqPPoleInfty} to global sections,
which gives $\overline{z^{-w}P(a,b)} = z^{-w}P(\tau a,\tau b)$
for $a,b\in H^0(\dP^1,\widehat{\cH})$. This implies that
$$
h(a,b)=z^{-w}P(a,\tau b)=\overline{z^{-w}P(\tau a, b)}
=\overline{z^{-w}P(b, \tau a)}=\overline{h(b,a)}
$$
which is what we need.
\end{proof}
The last lemma motivates the following definition.
\begin{definition}
A TERP-structure is called pure iff the bundle $\widehat{H}$ is trivial.
A pure TERP-structure is called polarized iff the hermitian
form $h:H_0\times H_0\rightarrow \dC$ is positive definite.
\end{definition}

The next result introduces one of the most interesting objects attached
to a pure TERP-structure, namely, an endomorphism of $H_0$ that
was considered in \cite{CFIV} under the name ``new supersymmetric index''.
Its eigenvalues are related to and can be considered (in the regular singular
case, see chapter \ref{c6})
as a generalization of the spectral numbers of $H,\nabla$.
Let $\cU$ be the pole part of the connection $\nabla$ on $H$, i.e.
$\cU$ is an endomorphism of the fibre $H_0$ defined
by $\cU=[z\nabla_{z\partial_z}]$.
\begin{lemma}
Suppose that $(H,H'_\dR,\nabla,P,w)$ is a pure
TERP-structure. Then there exists an endomorphism $\cQ$ of $H_0$ such
that for any $\omega\in H^0(\dP^1, \widehat{\cH})\cong H_0$, we have
$$
\nabla_{z\partial_z}\omega = \left(\frac{1}{z}\cU+
\left(\frac{w}{2}\mathit{Id}-\cQ\right)-z\overline{\cU}\right)\omega
$$
where $\overline{\cU}$ denotes the adjoint of $\cU$ with respect to $h$
and satisfies $\overline{U}=\tau\circ U\circ\tau$.
$\cQ$ is $h$-selfadjoint and anti-commutes with $\tau$.
If the TERP-structure is polarized, $\cQ$ is semi-simple and
its eigenvalues are real and symmetric with respect to zero.
\end{lemma}
\begin{proof}
It is obvious that there are endomorphisms $A,B$ such that
$\nabla_{z\partial_z}\omega=\left(\frac{\cU}{z}+A+zB\right)\omega$. We need
to show that $B=-\overline{\cU}$, that $\cQ:=\frac{w}{2}\mathit{Id}-A$ is
$h$-selfadjoint
and anti-commutes with $\tau$ and that $\overline{U}=\tau\circ U\circ\tau$.
All of these properties follow using the equations \eqref{eqTauFlat},
$\tau(z\omega)=z^{-1}\tau(\omega)$
and the fact that for $a,b\in H^0(\dP^1,\widehat{\cH})$ we have
$z^{-w}P(\cU a,b)=z^{-w}P(a,\cU b)$.
If $h$ is positive definite, $\cQ$ is semi-simple with real eigenvalues.
They are symmetric because of $\tau\circ\cQ=-\cQ\circ\tau$.
\end{proof}
We will give in the following definition/theorem a brief reminder on how to extend
the notion of a TERP-structure to the relative case, where parameters have to be taken into account.
The main reference is \cite[chapter 2]{He2}, in particular sections 2.4 to 2.7 of cit.loc.
\begin{deftheo}
Let $M$ be a complex manifold. A variation of TERP-structures over $M$ is a
tuple $(H,H'_\dR,\nabla,P,w)$ where $H\in\mathit{VB}_{\dC\times M}$
and $H'_\dR$ a maximal real subbundle
of the restriction $H':=H_{|\dC^*\times M}$.
The connection
$\nabla:\cH'\rightarrow\cH'\otimes\Omega^1_{\dC\times M}$
is flat and meromorphic
with a pole of Poincaré rank one along $\{0\}\times M$. The pairing
$P:\cH\otimes j^*\cH\rightarrow z^w\cO_{\dC\times M}$ is non-degenerate,
$(-1)^w$-symmetric, flat and sends $H'_\dR$ to $i^w\dR$.

The following facts hold:
\begin{enumerate}
\item
The construction of the extension to infinity generalizes and yields
a complex vector bundle $\widehat{H}$ over $\dP^1\times M$ with holomorphic
structure in $\dP^1$-direction, in other words, a locally free
$\cC^\infty_M\cO_{\dP^1}$-module.
The connection extends to $\widehat{H}$ with a pole of Poincaré rank one
along $\{0,\infty\}\times M$.
\item
The notion of pure resp. polarized pure TERP-structures
is defined as in the absolute case.
Given a variation of $(H,H'_\dR,\nabla,P,w)$ of pure TERP-structures,
the objects $h,\tau,\cU,Q$ are defined on
$p_*\cC^{\infty h}(\widehat{H})\cong \cO(H_{|z=0})\otimes_{\cO_M}\cC^\infty_M$.
Moreover precisely, the connection $\nabla$ takes the following form:
Let $\omega\in p_*\cC^{\infty h}(\widehat{H})$,
then
$$
\nabla\omega=
\left(D+\frac1z C+z\overline{C}+
\left(\frac1z \cU+\left(\frac{w}{2}\mathit{Id}-\cQ\right)
-z\overline{\cU}\right)\frac{\textup{d}z}{z}\right)\omega
$$
where $D$ is the Chern connection for $h$, $C$ a Higgs field
on $H_{|z=0}$ and $\overline{C}$ its $h$-adjoint. The operators
$D,C,\overline{C}, h,\tau,\cU,Q$ satisfy a couple of compatibility conditions
(\cite[equations 2.50-2.61]{He2}) making up what was called
$CV$-structure in cit.loc.
\end{enumerate}

\end{deftheo}

In the remaining part of this chapter, we will discuss the
relation of the notion of (variation of) TERP-structures
with (variation) of polarized integrable twistor structures.
Polarized twistor structures were defined in \cite{Si5}, and
the term integrable was first used in \cite{Sa6}.
We briefly recall the definitions.

As in \cite{Si5}, we denote for any
$\cO_{\dP^1}$-module $\cE$ by $\sigma^*\cE$ the sheaf defined by
$\Gamma(U,\sigma^*\cE):=\overline{\Gamma(\sigma(U),\cE)}$.
As before, the conjugate
complex structure is needed to ensure that $\sigma^*\cE$ is again
a sheaf of $\cO_{\dP^1}$-modules. Note that the convention here
differs from the one used for the map $\gamma$, but we prefer to
be compatible with the notations both in \cite{He2} and \cite{Si5}.
\begin{definition}\label{t3.13}
\begin{itemize}
\item
A twistor is a holomorphic bundle on $\dP^1$.
\item
A twistor $\widehat{H}\in\mathit{VB}_{\dP^1}$ is integrable
if it comes equipped
with a meromorphic connection with poles of order at most two at zero
and infinity.
\item
$\widehat{H}$ is called pure of weight $w$ iff it is semi-stable of slope
$w$, i.e., isomorphic to a sum
$\oplus_{i=1}^{\rank(\widehat{H})} \cO_{\dP^1}(w)$.
\item
A pairing on a twistor $\widehat{H}$ is a non-degenerate $(-1)^k$-symmetric morphism
$\widehat{S}:\cO(\widehat{H})\otimes_{\cO_{\dP^1}}\sigma^*\cO(\widehat{H})
\rightarrow \cO_{\dP^1}(2k)$. In case $\widehat{H}$
is integrable $\widehat{S}$ is required to be flat.
\item
A pure twistor $\widehat{H}$ of weight $w$ is called polarized by a pairing
$\widehat{S}$
iff the induced pairing
$$
\widehat{S}_w:\left(\cO(\widehat{H})\otimes\cO_{\dP^1}(-w)\right)
\otimes_{\cO_{\dP^1}}\sigma^*\left(\cO(\widehat{H})
\otimes\cO_{\dP^1}(-w)\right) \longrightarrow \cO_{\dP^1}
$$
induces a positive definite hermitian pairing on the space of global sections.
The pair $(\widehat{H},\widehat{S})$ is called polarized twistor structure
(abbreviation PTS).
\item
A twistor $\widehat{H}$ is mixed iff it is equipped with an
increasing filtration $\widehat{W}_\bullet$ by subbundles
such that each graded piece $Gr_k^{\widehat{W}}(\widehat{H})$
is a pure twistor of weight $k$.
\end{itemize}
\end{definition}
With these definitions in mind, we can state the following comparison lemma.
\begin{lemma} \label{t3.14}
Let $(H,\nabla, H_\dR, P,w)$ be a TERP-structure. Then
$(\widehat{H},\nabla)$ is an integrable twistor.
There is a naturally defined pairing $\widehat{S}:
\cO(\widehat{H})\otimes_{\cO_{\dP^1}}\sigma^*\cO(\widehat{H})
\rightarrow \cO_{\dP^1}$. $\widehat{H}$ is pure (and then
automatically of weight zero) iff $H$ is pure TERP and in that case
$\widehat{S}$ gives a polarization iff $H$ is polarized pure TERP.
\end{lemma}
\begin{proof}
The statements about $(\widehat{H},\nnn)$
are obvious from what has been said before,
it was shown in lemma \ref{lemDegHzero}
that the connection extends to $\widehat{H}$
as required and that $\deg(\widehat{H})=0$.
The pairing $\widehat{S}$ is defined as
$$
\widehat{S}(a,b):=(-1)^w P(a,\tau_{real} b)=z^{-w}P(a,\tau b) .
$$
$P$ has by definition a zero of order $w$ at the origin
and (as was shown) a pole of order $w$ at infinity, which implies
that $\widehat{S}$ maps to $\cO_{\dP^1}$. The flatness of $\widehat{S}$
follows from the flatness of $P$ and $\tau_{real}$.
The only thing that remains
to discuss is that if $\widehat{H}$ is pure then
it is polarized by $\widehat{S}$ precisely iff the TERP-structure we started with
is polarized pure TERP. But this is a tautology: $\widehat{S}$ polarizes $\widehat{H}$
iff the induced form on $H^0(\dP^1, \cO(\widehat{H}))$ is positive definite hermitian,
but this form is exactly $z^{-w}P(-,\tau(-))$. The positive definiteness of this
form on the space of global sections was the defining property for a pure TERP-structure
to be polarized pure TERP.
\end{proof}
\textbf{Remark:} A TERP-structure $(H,\nabla, H'_\dR, P)$ comes equipped with an
integer $w$, its weight. However, lemma \ref{lemDegHzero} shows that the twistor
$\widehat{H}$ constructed from such a TERP-structure is a $\dP^1$-bundle which always
has degree zero. Consequently, for a pure TERP-structure, the twistor $\widehat{H}$
is pure of weight zero, regardless of the value of $w$.
The reason for this is that the Tate twist, which is used to transform a pure twistor
of some weight in a twistor of weight zero is already implicitly contained in our gluing constructing
from definition \ref{defHhat}, namely, $\widehat{H}$ is defined by patching
$H$ and $\overline{\gamma^*H}$ via the map $\tau$, and not via $\tau_{real}$.

\bigskip

In \cite{Sa6}, the related notion of an $\cR$-triple was introduced.
The following lemma gives the comparison. We omit the proof, which is more or less straightforward.
\begin{lemma}
Given a TERP-structure $(H,\nnn,H'_\dR,P,w)$,
the tuple $(\cO(H),\cO(H),\widehat{S})$
where $\widehat{S}$ is the above pairing restricted to
$\cO(H')$
is a smooth $\cR$-triple. It is equal
to its hermitian adjoint $(\cO(H),\cO(H),\widehat{S}^*)$, so it is
polarized by $(\mathit{Id},\mathit{Id})$. It is even an object in $\cR\mathit{int}(\mathit{pt})$, the vertical connection
being $\nabla_z$ (and using the flatness of $\widehat{S}$).
\end{lemma}
As before, we need a relative version of the above notions
taking parameters into account. It can be formulated in two
equivalent ways.
\begin{definition}[Variation of twistor structures, harmonic bundles]\label{t3.18}
Let $M$ be a complex manifold.
\begin{enumerate}
\item
Consider a $\cC^\infty$ vector bundle on $\dP^1\times M$ together
with an integrable operator
$\overline{\partial}_{\dP^1}:\cC^\infty(E)\rightarrow\cC^\infty(E)\otimes_{\cC^\infty_{M\times\dP^1}}
\cC^\infty_M\cA^{0,1}_{\dP^1}$ defining a locally free sheaf $\cE$ of $\cC_M^\infty\cO_{\dP^1}$-modules.
Then $\cE$ is called a variation of twistor structures (abbreviation VTS)
if it comes equipped with an $\cO_{\dP^1}$-linear operator
$\bD:\cE\rightarrow \cE\otimes_{\cC_M^\infty\cO_{\dP^1}} \xi\cA^1_M$ satisfying
the following Leibniz rule:
$$
\bD(fe)=f\bD(e)+\mathbf{d}(f)e
$$
where $\xi\cA^1$ is the twistor $\cA^1_M \otimes_\dC \cO_{\dP^1}(1)$
(see \cite{Si5} or \cite{Mo2}) and $\mathbf{d}:\xi \cC^\infty_M\cong\cO_{\dP^1}\otimes_\dC \cC_M^\infty
\rightarrow \xi\cA_M^1$ is the natural ``twistor derivative''.
A variation of polarized twistor structures (abbreviation VPTS) is a variation equipped with
a polarization as in definition \ref{t3.13} which is flat with respect to $\bD$.
Similarly, one defines variations of pure resp. mixed twistor structures.
\item
A $\cC^\infty$-bundle $E$ on $M$ is called
harmonic iff there are operators $\partial,\overline{\partial}, \theta, \overline{\theta}$
and a pairing $h$ with
$$
\begin{array}{rrccl}
\partial, \theta&:&\cC^\infty(E) & \longrightarrow & \cC^\infty(E) \otimes_{\cC^\infty_M } \cA^{1,0}_M \\
\overline{\partial}, \overline{\theta}&:&\cC^\infty(E) & \longrightarrow & \cC^\infty(E) \otimes_{\cC^\infty_M } \cA^{0,1}_M \\
h& :& \cC^\infty(E)\otimes_{\cC^\infty_M}\cC^\infty(\overline{E})& \longrightarrow & \cC^\infty_M
\end{array}
$$
where $\theta, \overline{\theta}$ are $\cC^\infty_M$-linear and
$\partial, \overline{\partial}$  are satisfying the Leibniz-rule, such that
$(\partial+\overline{\partial}+\theta+\overline{\theta})^2=0$. The pairing $h$ is
positive definite, $\overline{\theta}+\theta$ is $h$-self-adjoint and
$(\partial+\overline{\partial})$ is $h$-metric, i.e.
$dh(a,b)=h(\partial a, b) + h(a, \overline{\partial}b)$ and $\overline{d}h(a,b)
=h(\overline{\partial}a,b)+h(a,\partial b)$.
\end{enumerate}
\end{definition}

A basic result due to Simpson (\cite{Si5}) is that the category of variations of pure polarized twistor structures
of weight zero is equivalent to the category of harmonic bundles. This correspondence
will be used implicitly several times in the sequel. As a matter of notation,
for a variation of twistors $\widehat{E}$ on $M$,
we denote by $\cC^{\infty h}(\widehat{E})$ the sheaf of $\cC^\infty$-sections
of $\widehat{E}$ which are holomorphic in the $\dP^1$-direction, i.e., annihilated
by the operator $\overline{\partial}_{\dP^1}$.

The following extension of lemma \ref{t3.14} to the relative case
is a condensed version of \cite[theorem 2.19]{He2} (see also
\cite[corollary 7.2.6]{Sa6}).
\begin{lemma}\label{t3.19}
For any variation of polarized pure TERP-structures $(G,G'_{|\dC^*\times M},\nabla,P,w)$,
$\widehat{G}$ has the structure of a variation of pure polarized integrable twistors structures of weight
zero yielding a harmonic bundle $p_*\cC^{\infty h}(\widehat{G})$ on $M$ equipped
with operators $\cU, \cQ$ satisfying \cite[equations 7.2.5]{Sa6}.
\end{lemma}
\begin{proof}
We only remark how to define the connection in the parameter direction: By definition,
a variation of TERP-structures furnishes a ``horizontal'' connection
$\nabla_M:\cO(G)\rightarrow\cO(G)\otimes z^{-1}\cO_\dC\Omega^1_M$ which we
can of course extend to
$$
D'_M:\cC^{\infty h}(G) \longrightarrow \cC^{\infty h }(G)\otimes
\left(z^{-1}\cO_\dC\cC^\infty_M\cA^{1,0}_M\oplus\cO_\dC\cC^\infty_M\cA^{0,1}_M\right)
$$
Similarly, the morphism $\gamma^*\nabla_M$ induces the operator
$$
D'':\cC^{\infty h}(\overline{\gamma^*G})\longrightarrow
\cC^{\infty h}(\overline{\gamma^*G}) \otimes
\left(z \cO_{\dP^1\backslash\{0\}}\cC^\infty_M \cA^{0,1}_M
\oplus\cO_{\dP^1\backslash\{0\}}\cC^\infty_M\cA^{1,0}\right)
$$
Then we put $\bD:\cC^{\infty h}(\widehat{G})\rightarrow\cC^{\infty h}(\widehat{G})
\otimes\xi\cA^1_M$ by $\bD=zD'+z^{-1}D''$ where $z$ and $z^{-1}$ are seen
as global sections of $\cO_{\dP^1}(1)$.
\end{proof}
\textbf{Remark:} It should be more or less clear from
what has been said that the notion of integrable twistor
also extends to the relative case, namely, given
a variation  of twistor structures $(\cE,\bD)$ like before
then we call it integrable if on both charts $\dC\times M$
and $(\dP^1\backslash\{0\})\times \overline{M}$ the relative meromorphic connections
defined by $\bD$ can be completed to an absolute meromorphic
connection having poles of Poincaré rank one at $\{0\}\times M$ and
$\{\infty\}\times \overline{M}$. This also explains the term integrable.
Each restriction $\cE_{\dP^1\times\{x\}}$ for $x\in M$ is then naturally
an integrable twistor as defined before. For more details,
see \cite[chapter 7]{Sa6}.

\section{Nilpotent orbits of TERP-structures}\label{c4}
\setcounter{equation}{0}

After the general discussion of TERP-structures and variations
of them in the previous chapter, we will introduce now a particular
class of such variations over one-dimensional bases, these are called
nilpotent orbits. The name is derived from corresponding objects in
Hodge theory (\cite{Sch}), and it will become clear later (see lemma \ref{lemOrbTERPOrbHS})
that nilpotent orbits of TERP-structures give rise
to nilpotent orbits of Hodge-structures in the classical sense.

We will simultaneously consider two points of view: Starting with
a single TERP-structure, there is a canonical way to construct
variations over punctured discs, on the other hand, given such a variation,
we will give a simple criterion to decide whether it is a nilpotent orbit.

\begin{definition}\label{t4.1}
Consider the following holomorphic maps: for any $r\in \dC^*$, let $\pi_r:\dC\rightarrow\dC$
be the multiplication by $r$ and define
$\pi, \pi':\dC\times\dC^* \rightarrow \dC $ by $\pi(z,r)=zr$ and $\pi'(z,r)=zr^{-1}$.
\begin{enumerate}
\item
Let $(H,H'_\dR,\nabla,P,w)$ be a TERP-structure. Then  $\pi_r^*(H,H'_\dR,\nabla,P,w)$ is
also TERP, and we say that $(H,H'_\dR,\nabla,P,w)$ induces a nilpotent orbit iff
$\pi_r^*(H,H'_\dR,\nabla,P,w)$ is a polarized pure TERP-structure for any $r\in\dC^*$
with $|r| \ll 1 $. Similarly, $(H,H'_\dR,\nabla,P,w)$ is said to induce a Sabbah orbit
iff $\pi_r^*(H,H'_\dR,\nabla,P,w)$ is polarized pure TERP for $|r|\gg 0$, i.e.
iff $\pi_{r^{-1}}^*(H,H'_\dR,\nabla,P,w)$ is polarized pure TERP for $|r|\ll 1$
\item
Let $M=\Delta^*$ be a punctured unit disk with coordinate $r$,
and $(G, G'_\dR,\nabla, P, w)$ be a variation of TERP-structures on $M$ (it might even be defined on
a larger disc or on the whole of $\dC^*$). Then we call it nilpotent orbit (resp. Sabbah orbit) iff
\begin{enumerate}
\item
The sheaf $\cO(G)$ is stable under $\nabla_{z\partial_z-r\partial_r}$ (resp. under
$\nabla_{z\partial_z+r\partial_r}$).
\item
For any $r\in M$ with $|r|\ll 1$, the restriction $(G,\nabla, H'_\dR, P,w)_{|\dC\times\{r\}}$
is a polarized pure TERP-structure.
\end{enumerate}
\end{enumerate}
\end{definition}
We will show that if a TERP-structure induces a nilpotent orbit
resp. a Sabbah orbit, then the family $\pi^*(H,\nnn,H'_\dR,P,w)$
(resp. $(\pi')^*(H,\nnn,H'_\dR,P,w)$) is a nilpotent orbit resp. Sabbah orbit according to
the second definition above and that vice versa, any nilpotent resp. Sabbah orbit
is of this type.  For this purpose and
also for later use, we will discuss the different types of monodromy involved
in this situation. For any flat bundle on $\dC^*\times \Delta^*$, we call the monodromy
of the loop $(z_0e^{i\varphi},r_0)$ vertical, and that of $(z_0,r_0e^{i\varphi})$ horizontal
(here $z_0,r_0\neq 0$ and $\varphi\in[0,2\pi)$).
We will mainly treat nilpotent orbits, and only comment on the case of Sabbah orbits which
is quite analogous.
\begin{lemma}\label{t4.2}
Suppose that $(H,\nabla, H'_\dR, P,w)$ induces a nilpotent orbit. Let $G:=\pi^*H$.
Consider the restriction of $G_{|\pi^{-1}(z_0)}$ to a fibre of $\pi$ (which is
isomorphic to $\dC^*$). This bundle
is canonically trivialized and we denote by $\rho_{r,z_0}:G_{(z_0 r^{-1},r)}\rightarrow G_{(z_0,1)}$
the identification of the fibres of $G$.
For $z_0\in\dC^*$ this trivialization is given by the flat structure of $G_{|\pi^{-1}(z_0)}$.
In particular, the monodromy of $G_{|\pi^{-1}(z_0)}$ is
trivial. This implies that vertical and horizontal monodromy of $G_{|\dC^*\times M}$ coincide.
\end{lemma}
\begin{proof}
The first two statements are clear, they follow from the definition of the pull-back
of a bundle with flat connection (see, e.g., \cite{Sa4}). The monodromy of a flat bundle
on $\pi^{-1}(z_0)\cong\dC^*$ which has a basis of flat sections is obviously trivial.
Moreover, a counter-clockwise oriented loop inside a fibre $\pi^{-1}(z_0)$ is homotopic to the composition
of $(z_0e^{i\varphi},r_0)$ and $(z_0,r_0e^{-i\varphi})$, so that
vertical and horizontal monodromy must be equal.
\end{proof}
\begin{lemma}\label{t4.3}
A TERP-structure $(H,\nabla, H'_\dR, P,w)$ induces a nilpotent orbit (resp. a Sabbah orbit) iff the variation
$\pi^*(H,\nabla, H'_\dR, P,w)$ (resp. $\pi'^*(H,\nabla, H'_\dR, P,w)$) is a nilpotent orbit
(resp. a Sabbah orbit). Any nilpotent orbit (resp. Sabbah orbit) on a
punctured disk containing $1$  is induced
from its restriction to $\dC\times\{1\}$ by the map $\pi$ (resp. $\pi'$).
\end{lemma}
\begin{proof}

Suppose that $(H,\nabla,H'_\dR,P,w)$ induces a nilpotent orbit. The restriction
of $\pi^*(H,\nabla,H'_\dR,P,w)$ to $\dC\times\{r\}$ is equal to $\pi_r^*(H,\nabla,H'_\dR,P,w)$ by definition and
therefore a polarized pure TERP-structure. We need to show that the family $G:=\pi^*H$ is a variation
of TERP-structures, i.e, that the connection has a pole of Poincaré rank at most one along $\{0\}\times\Delta^*$
and that $\cO(\pi^*H)$ is stable under $\nabla_{z\partial_z-r\partial_r}$.
It is readily checked that the fibres of $\pi$ are precisely the integral curves
of the vector field $z\partial_z-r\partial_r$. Lemma \ref{t4.2} implies
that for any $\sigma\in\cO(H)$, the pullback $\pi^*\sigma\in\pi^{-1}\cO(H)\subset
\cO_{\dC\times M}\otimes\pi^{-1}\cO(H)=\cO(\pi^*H)$
satisfies $\nabla_{z\partial_z-r\partial_r}\pi^*\sigma=0$. By definition,
$\cO(\pi^*H)$ is generated by such sections $\pi^*\sigma$, therefore,
$\cO(\pi^*H)$ must be stable under $\nabla_{z\partial_z-r\partial_r}$.
Moreover, it is clear that $\nabla_{z^2\partial_z}\cO(\pi^*H)\subset\cO(\pi^*H)$,
as this is true for the restriction to any $\dC\times\{r\}$. Finally, it follows from the identity
$z\partial_r=\frac{1}{r}z^2\partial_z- \frac{z}{r}\left(z\partial_z-r\partial_r\right)$
that $\cO(\pi^*H)$ is also stable under $\nabla_{z\partial_r}$. In conclusion,
we have shown that $\pi^*(H,\nabla,H'_\dR,P,w)$ is a nilpotent orbit.

Conversely, given any bundle $G\in\mathit{VB}_{\dC\times \Delta^*}$ which underlies
a variation of polarized pure TERP-structures and satisfies
$\nabla_{z\partial_z-r\partial_r}\cO(G)\subset \cO(G)$, we need to see that it
is of the form $G=\pi^* G_{|\dC\times\{1\}}$. This follows from the
calculation done in lemma 7.19 of \cite{He2}, i.e., the fact that
if $\nabla_{z\partial_z-r\partial_r}$ sends $\cO(G)$ to itself then
any base $\underline{e}$ of this sheaf is related to the base
$\pi^*(\underline{e}_{|\dC\times\{1\}})$ by a base change in $\mathit{GL}(\cO_{\dC\times\Delta^*})$.
\end{proof}

In the following lemma, we consider the bundle $\widehat{H}$ over $\dP^1$ obtained
from a TERP-structure by the construction of definition \ref{defHhat}. If we are given
a variation of TERP-structures, then by \cite{He2}, lemma 2.14 (e), one ends up with
a real analytic family of holomorphic $\dP^1$-bundles. In the case of a variation
$\pi^*H$ of the above type, the following statement shows that the canonical identification
$\rho_{r,z_0}$ of fibres of $\pi^*H$
over points $(z,r)\in\pi^{-1}(z_0)$ extends to an identification of the $\dP^1$-bundles
$\widehat{\pi^*H}_{|\dP^1\times\{r\}}$ for all $r\in S^1$ and more generally
we have $\widehat{\pi^*H}_{|\dP^1\times\{r\}}\cong\widehat{\pi^*H}_{|\dP^1\times\{r'\}}$
if $|r|=|r'|$.
This means that the ``interesting'' part of the variation $\pi^*H \rightarrow \dC\times M$
is the restriction to rays in $M$ with fixed argument.

\begin{lemma}\label{t4.4}
For any $r\in S^1$, there is a
canonical bundle isomorphism $\widehat{\pi_r^*H}\rightarrow \widehat{H}$ which restricts
over any $z_0\in\dC\subset\dP^1$ to the isomorphism of fibres
$\rho_{r,z_0}:(\pi^*H)_{|(z_0\cdot r^{-1}, r)} \rightarrow H_{|z_0}$ considered in lemma \ref{t4.2}.
\end{lemma}
\begin{proof}
The pointwise isomorphisms $\rho_{r,z_0}$ glue to a bundle isomorphisms
$\phi_r:\pi_r^*H \rightarrow H$.
Moreover, for $r\in S^1$ we have $1/r=\overline{r}$, and we can also glue the maps
$\rho_{r,\overline{z_0}^{-1}}:(\pi^*H)_{|(\overline{r\cdot z_0^{-1}}, r)}
\rightarrow H_{|\overline{z_0}^{-1}}$
for all $z_0\in\dP^1\backslash\{0\}$ to an isomorphism
$\widetilde{\phi}_r:\overline{\gamma^*(\pi_r^*H)}
\rightarrow \overline{\gamma^*H}$. If we twist $\widetilde{\phi}_r$ by
multiplication with $r^w$, we obtain the following commuting diagram, which shows
that $\widehat{\pi_r^*H}\cong\widehat{H}$.
$$
\xymatrix@!0{
{\SC \textup{flat shift of }\overline{(z\frac1r)^{-w}b}}\hspace*{1cm}&&
\overline{\gamma^*\pi_r^*H}  \ar[rrr]_{\cong}^{r^w\cdot \widetilde{\phi}_r}
& &  &
\overline{\gamma^*H}
& &
{\SC \textup{flat shift of }\overline{z^{-w}a}}
\\ \\ \\
{\SC b\ar@{|->}[uuu]}&&\pi_r^*H  \ar[uuu]^\tau \ar[rrr]^{\cong}_{\phi_r}  & &  & H \ar[uuu]_\tau
&& {\SC a}\ar@{|->}[uuu]
}
$$
\end{proof}
The lemma shows in particular that for any $r\in S^1$, $\widehat{\pi_r^*H}$ is a trivial bundle
iff $\widehat{H}$ is so. We obtain as a consequence:
\begin{corollary}
For any bundle $G\in\mathit{VB}_{\dC\times \Delta^*}$ underlying a variation of TERP-structures
with the property that $\nabla_{z\partial_z-r\partial_r}\cO(G)\subset\cO(G)$, the subset
$\left\{r\in \Delta^*\,|\, \widehat{G}_{|\dP^1\times\{r\}}\textup{ is not trivial}\right\}$
is either the whole of $\Delta^*$ or a discrete union $\bigcup S_c$ with $S_c=\left\{r\in\dC^*\,|\,|r|=c\right\}$
or empty. On the complement, one obtains a real analytic hermitian form $h$ by the procedure of
lemma \ref{lemhIsHerm}. Its signature is constant on any connected component of this complement.
\end{corollary}
We finish this chapter by describing some properties of the harmonic bundle $p_*\cC^{\infty h}(\widehat{H})$
associated to a nilpotent orbit or Sabbah orbit (here $p:\dP^1\times\Delta^*\rightarrow \Delta^*$ is the projection).
Any harmonic bundle $(E,D,C,h)$ on a punctured disk with coordinate $r$ is called \emph{tame} iff
the eigenvalues of $C_{r\partial_r}\in\cE\!nd(E)$ are bounded at zero (see \cite{Si2}).
In our situation, the endomorphism $C_{r\partial_r}$ of the harmonic bundle associated
to a nilpotent orbit or Sabbah orbit $(G,\nabla,G'_\dR,P,w)$ is
given by $C_{r\partial_r}=[r\nabla_{z\partial_r}] \in
{\cE\!nd}_{\cO_M}(G_{|z=0})$. Recall also that
the pole part $\cU$ is defined as $\cU=[z\nabla_{z\partial_z}]\in {\cE\!nd}_{\cO_M}(G_{|z=0})$.
\begin{lemma}\label{lemTameness}
\begin{enumerate}
\item
The harmonic bundle associated to a nilpotent orbit of TERP-structures is tame iff
the pole part $\cU$ of the TERP-structure is nilpotent.
\item
The harmonic bundle associated to a Sabbah orbit of TERP-structures is always tame.
\end{enumerate}
\end{lemma}
\begin{proof}
Let $(H,\nabla,H'_\dR,P,w)$ be a TERP-structure inducing a
nilpotent orbit $G:=\pi^*H$, then the fact that $\nabla_{z\partial_z-r\partial_r}\pi^{-1}\cO(H)=0$
implies that $C_{r\partial_r}=\cU\in {\cE\!nd}_{\cO_M}(G_{|z=0})$. If
$(H,\nabla,H'_\dR,P,w)$ induces a Sabbah orbit, then the same argument shows
$C_{r\partial_r}=-\cU\in {\cE\!nd}_{\cO_M}(G_{|z=0})$. Therefore, it suffices
in both cases to study the behavior of the eigenvalues of $\cU$.

If $G$ is a nilpotent orbit, then for any
$r\in M$ we have that $G_{|\dC\times\{r\}}=\pi_r^*H$, which implies that
$\cU_{|r}=r^{-1} \cdot \rho_{0,r}^{-1} \circ \cU_{|1} \circ \rho_{0,r}$
whereas for a Sabbah orbit,
$G_{|\dC\times\{r\}}=\pi_{r^{-1}}^*H$ so that
$\cU_{|r}=r \cdot \rho_{0,r}^{-1} \circ \cU_{|1} \circ \rho_{0,r}$.
This shows that in the first case, the eigenvalues are bounded iff $\cU_{|1}$ is nilpotent
whereas in the second case they tend to zero as $r$ approaches the origin.
\end{proof}

\section{PMHS and integrable PMTS}
\label{c5}
This chapter is devoted to establish correspondences
between (polarized) mixed Hodge structures and
some particular TERP-structures resp. integrable twistors.
This is merely an extension to integrable twistors of the correspondence
due to Simpson (\cite{Si5} and \cite{Mo2}). This will allow us to use
Mochizuki's main result (\cite[theorem 12.1]{Mo2},
the theorem of the limit mixed twistor structure) in chapter \ref{c6}.

We start by giving a very brief reminder on how to associate
linear algebra data to the ``topological part'' of a TERP-structure.
The term ``topological'' refers to the following fact: If a TERP-structure
$(H,\nabla,H'_\dR,P, w)$ arises from singularity theory (see chapter \ref{c11}),
then the restriction $H':=H_{|\dC^*}$ is a topological invariant of
the singularity and the extension $H$ to a bundle on $\dC$ with meromorphic
connection is of transcendental nature.
We thus stick to $(H',\nabla,H'_\dR,P,w)$ for a moment. The following objects are either well
known or have been discussed extensively in \cite[7.2]{He2}.
\begin{itemize}
\item
$H^\infty\supset H_\dR^\infty$: the spaces of complex resp. real flat
multivalued global sections of $(H',\nabla)$.
\item
$M\in\mathit{Aut}(H^\infty_\dR)$: the monodromy of $\nabla$, which decomposes
into $M=M_s\cdot M_u$ with $M_s$ semi-simple and $M_u$ unipotent, $N:=\log(M_u)$
the nilpotent part.
\item
For any $\lambda\in \dC$, the generalized eigenspace $H^\infty_\lambda$ of
$M$ for the eigenvalue $\lambda$ and the corresponding flat subbundle
$H'_\lambda\subset H'$. Moreover, let $H^\infty_{\arg=0}=\oplus_{\arg(\lambda)=0} H^\infty_\lambda$,
$H^\infty_{\arg\neq 0}=\oplus_{\arg(\lambda)\neq0}
H^\infty_\lambda$ and $H^\infty_{\neq 1}=\oplus_{\lambda\neq 1} H^\infty_\lambda$.
\item
Elementary holomorphic sections of $H'$: fix $A\in H^\infty_\lambda$ and $\alpha\in\dC$
with $\lambda=e^{-2\pi i\alpha}$, then define
$es(A,\alpha):=z^{\left(\alpha\mathit{Id}-\frac{N}{2\pi i}\right)}A(z):=
e^{\left(2\pi i \alpha \mathit{Id}-N\right)\zeta}A(\zeta)$ with $\zeta$
a coordinate on the universal covering of $\dC^*$ such that $e^{2\pi i\zeta}=z$; the space of elementary
sections of order $\alpha$ is denoted by $C^\alpha$, and sending
$A\mapsto es(a,\alpha)$ defines an isomorphism $\psi_\alpha: H^\infty_{e^{-2\pi i \alpha}}
\rightarrow C^\alpha$. The connection acts as follows:
$\nabla_{z\partial_z} es(A,\alpha) = \alpha \cdot es(A,\alpha)-es\left(\frac{N}{2\pi i}A,\alpha\right)$.
\item
Polarizations: The pairing $P$ induces two bilinear forms on $H^\infty$, which are
equivalent to each other. The first one, called $L$ here, corresponds to the Seifert form
in singularity theory: first fix any $z\in \dC^*$ and define
$L:H_z\times H_z\rightarrow \dC; (a,b)\mapsto P(a,\gamma_\pi(b))$,
where $\gamma_\pi$ is the counter-clockwise flat shift from the fibre at $z$
to the fibre at $-z$. It is readily checked that $L$ is monodromy invariant, so that
we get a pairing on $H^\infty$. Then the following formulas, which might seem
artificial at first sight, define a pairing $S$ on $H^\infty_\dR$:
\begin{eqnarray}\label{5.9}
S(a,b) &:=& (-1)(2\pi i)^wL\left(a,\frac{1}{M-\mathit{Id}}b\right)
\;\;\;\;\;\;
\forall a,\forall b\in H^\infty_{\arg\neq 0},\nonumber \\ \\
S(a,b) & := & (2\pi i)^w
L\left(a,\frac{2\pi i (-\beta)\mathit{Id} + N}{M-\mathit{Id}} b\right)
\;\;\;\;\;\; \forall a\in H^\infty,\forall b\in H^\infty_{e^{-2\pi i\beta}}
\nonumber
\end{eqnarray}
where $\beta\in i\dR$. If $\beta=0$ (i.e., $b\in H_1$) we put
$$
\frac{2\pi i (-\beta) + N}{M-\mathit{Id}} :=
\left(\sum\limits_{k\geq1}\frac{1}{k!}N^{k-1}\right)^{-1}.
$$
The pairing $S$ is monodromy invariant, nondegenerate , $(-1)^{w-1}$-symmetric on $H_{\arg\neq 0}$
and $(-1)^w$-symmetric on $H_{\arg=0}$. It takes real values on $H_\dR$.
\item
A topological version of the Fourier-Laplace transformation:
Let $G^{(\alpha)}$ for $\alpha$ with $\Re(\alpha)>0$
be the automorphism of $H^\infty_{e^{-2\pi i \alpha}}$ defined as:
$$
G^{(\alpha)} := \sum_{k\geq 0}\frac{1}{k!}\Gamma^{(k)}(\alpha)
\left( \frac{-N}{2\pi i}\right)^k
=: \Gamma \left(\alpha\cdot \id - \frac{N}{2\pi i}\right) .
$$
Here $\Gamma^{(k)}$ is the $k$-th derivative of the gamma function.
For notational convenience, we let
\begin{equation}\label{eqGamma}
G:=\sum_{\Re(\alpha)\in(0,1]} G^{(\alpha)} \in
\mathit{Aut}\left(H^\infty=\oplus_{\alpha}
H^\infty_{e^{-2\pi i\alpha}}\right) .
\end{equation}
The following identities hold true.
\begin{itemize}
\item
Let $\tau$ be another coordinate on $\dC$. Then
\begin{eqnarray}
\int_0^{\infty\cdot z}
e^{-\tau/z}\cdot es(A,\alpha-1)(\tau)d\tau=
es(G^{(\alpha)}A,\alpha)(z).
\end{eqnarray}
Here the left hand side means that the values of the section for different $\tau$'s
are shifted using $\nabla$ to the fibre over $z$ and then summed up (as an integral).
\item
Let $\alpha,\beta \in (0,1)+i\dR$ and $A\in \hiia, B\in \hiib$. Then
\begin{eqnarray} \label{5.18}
P\left( es(G(A),\alpha),\, es(G(B),\beta)\right)(z)
=z\frac{1}{(2\pi i)^{w-1}}\cdot S(A,B).
\end{eqnarray}
\item Let
$\alpha,\beta\in 1+i\dR$ and
$A\in \hiia, B\in \hiib$. Then
\begin{eqnarray} \label{5.19}
P\left( es(G(A),\alpha),\, es(G(B),\beta)\right)(z)
=z^2\frac{-1}{(2\pi i)^{w}}\cdot S(A,B).
\end{eqnarray}
\end{itemize}
In fact, the automorphism $G$ and the pairing $S$ were defined in \cite[7.2]{He2} in a slightly
more restricted situation, namely, it was assumed in that paper that
the eigenvalues of the monodromy are in $S^1$. However, the generalization
we consider here can be shown by the same type of calculations.
\end{itemize}

It is also possible to go the other way round, i.e., to construct a
tuple $(H',\nabla, H'_\dR,P,w)$ starting from the vector spaces
$H^\infty\supset H^\infty_\dR$, an automorphism $M$ and a pairing $S$ as above.
The first part is obvious as a flat bundle $(H'_\dR,\nabla)$ is
equivalent to the data $(H^\infty_\dR,M)$ (similarly for the complexifications).
Thus the only thing to show is how to define $P$ starting from
$(H^\infty,H^\infty_\dR,M\in\mathit{Aut}(H_\dR), S)$. First remark
that the formulas \eqref{5.9} can be reversed to define $L$ starting from $S$.
Note that $S$ takes real values on $H^\infty_\dR$ but $L$ sends $H^\infty_\dR$
to $i^w\dR$. Then put $P(a,b)(z):=L\left(a(z),\gamma_{-\pi} b(-z)\right)$.
This gives a flat pairing as $L$ is defined on the space of flat sections.
The symmetry property for $P$ follows from that of $L$ (i.e,
one needs to check that $L(a,b)=(-1)^wL(Mb,a)$ which is straightforward).
We therefore arrive at the following basic result.
\begin{lemma}[Correspondence of topological data]\label{lemCorrTop}
There is a one to one correspondence between tuples consisting of
$(H^\infty,H_\dR^\infty, M, S,w)$ with the above properties
and flat bundles $H'\supset H'_\dR$ on $\dC^*$ equipped with a flat
$(-1)^w$-symmetric pairing $P:\cO(H')\otimes j^*\cO(H') \rightarrow \cO_{\dC^*}$
sending (pointwise) $H'_\dR$ to $i^w\dR$.
\end{lemma}

We will  gradually enrich this result to get eventually a correspondence between
sums of two polarized mixed Hodge structures with an automorphism having eigenvalues
in $S^1$ on the one hand and integrable polarized mixed twistor structures ``generated by
elementary sections'' on the other hand. The first step is to construct a twistor
from one additional piece of data, namely, a filtration. This is the result of the next
lemma which is quite close to \cite[lemma 7.12]{He2}.
\begin{lemma}
\label{lemCorrElTwist}
Let $(H^\infty,H^\infty_\dR,M)$ be as above and $F^\bullet$
an exhaustive decreasing filtration on
$H^\infty$. Suppose that it is invariant under $M_s$ and satisfies
$NF^\bullet\subset F^{\bullet-1}$.
Consider the flat (complex) bundle $H'$ on $\dC^*$ corresponding to
$(H^\infty,M)$.
Define the following sheaves
\begin{eqnarray}
\begin{array}{rcl}
\cH&:=&\sum\limits_{p\in\dZ, \Re(\alpha)\in(0,1]}
\cO_\dC z^{\left(\alpha+w-1-p\right)\mathit{Id}-\frac{N}{2\pi i}}
F^p H^\infty_{e^{-2\pi i \alpha}} \\ \\
&=& \sum\limits_{
\begin{array}{c}
\SC p\in\dZ, \Re(\alpha)\in(0,1]\\
\SC A\in F^pH^\infty_{e^{-2\pi i \alpha}}
\end{array}} \cO_\dC \mathit{es}(A,\alpha+w-1-p),\\ \\
\widetilde{\cH}&:=&\sum\limits_{q\in\dZ, \Re(\alpha)\in(0,1]}
\cO_{\dP^1\backslash\{0\}} z^{\left(q+\alpha-[\Re(\alpha)]\right)
\mathit{Id}-\frac{N}{2\pi i}}\overline{F}^q H^\infty_{e^{-2\pi i \alpha}} \\ \\
&=& \sum\limits_{
\begin{array}{c}
\SC q\in\dZ, \Re(\alpha)\in(0,1]\\
\SC A\in \overline{F}^qH^\infty_{e^{-2\pi i \alpha}}
\end{array}} \cO_{\dP^1\backslash \{0\}} \mathit{es}(A,q+\alpha-[\Re(\alpha)]).
\end{array}
\end{eqnarray}
Theses are locally free extensions of $\cO(H')$ to zero and infinity defining
a bundle $\widehat{H}\in\mathit{VB}_{\dP^1}$. The connection extends to
$\widehat{H}$ with poles of order at most two at $0$ and $\infty$, therefore $\widehat{H}$
is integrable. Twistors of this type are called generated by elementary sections.

The connection has logarithmic poles at zero and infinity if $M$ is semi-simple.
$\widehat{H}$ is pure of weight k iff
$(H_{\arg\neq 0}, H_{\dR,\arg\neq 0}, F^\bullet)$ is a Hodge structure of
weight $w+k-1$ and $(H_{\arg=0}, H_{\dR,\arg=0}, F^\bullet)$ is a Hodge structure of
weight $w+k$.
\end{lemma}
\begin{proof}
The inclusions $\cH \subset i_*\cO(H')$ and $\widetilde{\cH}\subset \widetilde{i}_*\cO(H')$, where
$i:\dC^*\rightarrow \dC$ and $\widetilde{i}:\dC^*\rightarrow \dP^1\backslash\{0\}$, are obvious.
Consequently, we obtain a bundle $\widehat{H} \in \mathit{VB}_{\dP^1}$.
If $M$ is the identity, then this is precisely the Rees construction (\cite{Si5}, \cite{Mo2}).
The connection on $H'$ extends to $\widehat{H}$ (with poles at zero and infinity) using the Leibniz rule.
Then for any $A\in F^pH_{e^{-2\pi i \alpha}}$ we have
$$
\begin{array}{c}
(z^2\nabla_z)(e^{\log(z)\left(\left(\alpha+w-1-p\right)\mathit{Id}-\frac{N}{2\pi i}\right)}A)=
z\left(\left(\alpha+w-1-p\right)\mathit{Id}-\frac{N}{2\pi i}\right)
z^{\left(\alpha+w-1-p\right)\mathit{Id}-\frac{N}{2\pi i}}A\\ \\
=
z^{\left(\alpha+w-p\right)\mathit{Id}-\frac{N}{2\pi i}}
\left(\left(\alpha+w-1-p\right)\mathit{Id}-\frac{N}{2\pi i}\right)A
\in z^{\left(\alpha+w-p\right)\mathit{Id}-\frac{N}{2\pi i}} F^{p-1}H_{e^{-2\pi i \alpha}} \subset \cH.
\end{array}
$$
A similar calculation shows that for $A\in\overline{F}^q H_{e^{-2\pi i \alpha}}$ we have
$$
z^{-2}\nabla_{\partial_{z^{-1}}}(z^{\left(q+\alpha-[\Re(\alpha)]\right)\mathit{Id}-\frac{N}{2\pi i}}A)
\in \widetilde{\cH}.
$$
In both cases the fact that $NF^\bullet\subset F^{\bullet-1}$ is essential.
We see that the connection has a pole of order at most two at zero and infinity.
For $N=0$ it follows immediately that $\cH$ is stable under $z\nabla_{\partial_z}$
and $\widetilde{\cH}$ is stable under $z^{-1}\nabla_{\partial_{z^{-1}}}$.

In general, for $A\in (F^p\cap \overline{F}^{(w+k)-1-p})H^\infty_{arg\neq 0}$
we have that
$$
z^{\left(\alpha+w-1-p\right)\mathit{Id}-\frac{N}{2\pi i}}
A \in z^{-k}z^{\left((w+k)-1-p+\alpha-[\Re(\alpha)]\right)\mathit{Id}-\frac{N}{2\pi i}}\overline{F}^{(w+k)-1-p}H^\infty_{e^{-2\pi i \alpha}}
\subset z^{-k}\widetilde{\cH}
$$
and similarly for $A\in (F^p\cap \overline{F}^{(w+k)-p})H^\infty_{arg=0}$ (i.e., $\alpha\in1+i\dR$)
$$
z^{\left(\alpha+w-1-p\right)\mathit{Id}-\frac{N}{2\pi i}}
A \in z^{-k}z^{\left((w+k)-p+\alpha-[\Re(\alpha)]\right)\mathit{Id}-\frac{N}{2\pi i}}\overline{F}^{(w+k)-p}H^\infty_{e^{-2\pi i \alpha}}
\subset z^{-k}\widetilde{\cH}.
$$
Suppose that $(H^\infty,H^\infty_\dR,F^\bullet)$ is a sum of two Hodge structures of weights $w+k-1$ and $w+k$, i.e.,
$$
H^\infty=\bigoplus_{p} \left(F^p\cap \overline{F}^{(w+k)-1-p}\right) H^\infty_{\arg\neq 0} \oplus
\bigoplus_{p} \left(F^p\cap \overline{F}^{(w+k)-p}\right) H^\infty_{\arg = 0}.
$$
Then by choosing an appropriate basis, we see from the last two formulas that the constructed bundle $\widehat{H}$
is semi-stable of slope $k$, thus we obtain an integrable pure twistor
of weight $k$.
Conversely, if we know that $\widehat{H}$ is isomorphic to
$\cO^{\rank(H)}_{\dP^1}(k)$ then reading the last calculation
backwards gives us a basis of $H^\infty$ showing that
the filtration $\overline{F}^\bullet$, shifted appropriately, splits $F^\bullet$ .
\end{proof}

Consider now the general situation of a tuple $(H^\infty,H^\infty_\dR,M\in\mathit{Aut}(H^\infty_\dR),F^\bullet H^\infty,w)$,
where $F^\bullet$ and $\overline{F}^\bullet$ are not necessarily $w+k$ resp. $w+k-1$-opposed. Then
the following construction is a refinement of the above correspondence.

Let $N$ be the nilpotent part of $M$, and let $W_\bullet$ be the weight filtration defined by $N$ as described
in lemma \ref{defWeigthFilt}, but centered at zero. By definition of $W_\bullet$, for any $k$ the endomorphism $M$ restricts
to an element in $\mathit{Aut}(W_k)$.
\begin{lemma}\label{t5.3}
The above construction applied to the tuple $(W_k, W_k \cap H^\infty_\dR, M_{|W_k}, F^\bullet\cap W_k, w)$
defines a subbundle $\widehat{W}_k\subset \widehat{H}$. We obtain a filtration of $\widehat{H}$
by subbundles. There is a naturally defined nilpotent morphism
$\widehat{N}:\cO(\widehat{H})\rightarrow\cO(\widehat{H})\otimes
\cO_{\dP^1}(2)$ sending $\widehat{W}_k$ to $\widehat{W}_{k-2}\otimes
\cO_{\dP^1}(2)$.
The quotient bundle $Gr_k^{\widehat{W}}\widehat{H}$
is an integrable twistor and is canonically isomorphic to the
integrable twistor corresponding to $(Gr^W_k H^\infty, Gr^W_k H^\infty_\dR, M\in\mathit{Aut}
(Gr^W_kH^\infty), F^\bullet Gr^W_kH^\infty, w)$. The induced $M$ on $Gr_k^W H^\infty$ is semi-simple, so that
$(Gr_k^{\widehat{W}}(\widehat{H}),\nabla)$ is logarithmic.
\end{lemma}
\begin{proof}
The only point to understand is the definition of $\widehat{N}$: this is merely
the functoriality of the above construction (which we will not discuss in detail here).
The map $N:H^\infty\rightarrow H^\infty$ satisfies by assumption $N(F^\bullet)
\subset F^{\bullet-1}$ and $N(\overline{F}^\bullet)
\subset \overline{F}^{\bullet-1}$. Define $\widehat{N}es(A,\alpha):=es(NA,\alpha)$,
then $\widehat{N}\cH \subset z^{-1}\cH $ and $\widehat{N}\widetilde{\cH}\subset z\widetilde{\cH}$ and,
so that $\widehat{N}:\cO(\widehat{H})\rightarrow\cO(\widehat{H})\otimes\cO_{\dP^1}(2)$.
By definition of the weight filtration $W_\bullet$ on $H$, $NW_k\subset W_{k-2}$, this
implies the same property for $\widehat{N}$ with respect to $\widehat{W}_\bullet$.
Note that the filtration $\widehat{W}$ may also be defined
from $\widehat{N}$ as in lemma \ref{defWeigthFilt}.
\end{proof}
The following is now an immediate consequence of lemma \ref{lemCorrElTwist}.
\begin{lemma}\label{t5.4}
$(\widehat{H}, \widehat{W}_\bullet)$ is a mixed twistor, that is, each
$Gr_k^{\widehat{W}}(\widehat{H})$ is
pure of weight $k$, if and only if the tuples
$(H^\infty_{\arg \neq 0},H^\infty_{\dR,\arg\neq 0},W_{\bullet+w-1}, F^\bullet)$
and $(H^\infty_{\arg=0},H^\infty_{\dR,\arg=0},W_{\bullet+w}, F^\bullet)$ are mixed
Hodge structures. Moreover, given a twistor generated by elementary sections,
one recovers the filtration $F^\bullet$ on $H^\infty$ by
\begin{equation}\label{5.5}
F^pH^\infty_{e^{-2\pi i \alpha}}:=\psi_\alpha^{-1}\left(z^{p+1-w}(C^{\alpha+w-1-p}\cap\cH)\right)
\end{equation}
for $\Re(\alpha)\in(0,1]$. We obtain a one to one correspondence between
mixed twistors generated by elementary sections and sums of two mixed
Hodge structures (on $H^\infty_{\arg\neq 0}$ and $H^\infty_{\arg=0}$).
\end{lemma}

The next step is to take into account the pairing $S$. We
suppose a tuple $(H^\infty,H^\infty_\dR, M, F^\bullet H^\infty, S)$ be given.
Remember that the data $(H^\infty,H^\infty_\dR, M, S)$ were already shown
to be equivalent to a flat bundle $(H',\nabla)$ over $\dC^*$ and
a flat pairing $P$ on opposite fibres.
For a moment, we will only consider the extension $\cH\in\mathit{VB}_\dC$
(i.e., we will not make use of $\overline{F}^\bullet$) and show
that if $F^\bullet$ satisfies a particular orthogonality property with respect to $S$, then
the pairing $P$ will have a zero of order $w$ on the extension $\cH$ making
$(H,\nnn,H'_\dR,P,w)$ into a TERP-structure. Define the filtration $\widetilde{F}^\bullet$ by
$\widetilde{F}^\bullet := G^{-1} F^\bullet$ where $G$ is the
automorphism defined by formula \eqref{eqGamma}. The said orthogonality property will be formulated
using this twisted filtration $\widetilde{F}^\bullet$.
\begin{lemma}\label{t5.5}
Let $(H^\infty,H^\infty_\dR,M,F^\bullet,w,S)$ with all the properties
of lemma \ref{lemCorrTop} and \ref{lemCorrElTwist}
be given. Suppose that
\begin{equation}\label{eqOrthogonal}
\left(\widetilde{F}^p H^\infty_{\arg\neq 0}\right)^\bot = \widetilde{F}^{w-p}H^\infty_{\arg\neq 0}
\;\;\;\mbox{and}\;\;\;
\left(\widetilde{F}^p H^\infty_{\arg = 0}\right)^\bot = \widetilde{F}^{w+1-p}H^\infty_{\arg =0}.
\end{equation}
Here $^\bot$ is the orthogonal complement with respect to $S$.
Then $(\cH,\nnn,H'_\dR,P,w)$ is a TERP-structure.
\end{lemma}
\begin{proof}
In view of what has been said before, we need to show that $P$ sends
the germ $\cH_0$ to $z^w\cO_{\dC,0}$ and is nondegenerate when multiplied by $z^{-w}$. So take
two generators of this germ, i.e., let $A\in F^p H^\infty_{e^{-2\pi i\alpha}}, B\in F^q H^\infty_{e^{-2\pi i\beta}}$
(suppose first that $\alpha,\beta\in(0,1)+i\dR$), and compute
$$
\begin{array}{c}
P(es(A,\alpha+w-1-p),es(B,\beta+w-1-q))(z) =
P(z^{\left(\alpha+w-1-p\right)\mathit{Id}-\frac{N}{2\pi i}}A,
(-z)^{\left(\beta+w-1-q\right)\mathit{Id}-\frac{N}{2\pi i}}B)\\ \\
´=z^{2w-2-(p+q)}(-1)^{w-1-q}P(z^{\alpha\mathit{Id}-\frac{N}{2\pi i}}A,
(-z)^{\beta\mathit{Id}-\frac{N}{2\pi i}}B)
\\ \\
= z^{2w-2-(p+q)}(-1)^{w-1-q}\frac{z}{2\pi i} S\left(G^{-1}A,G^{-1}B\right) =
z^{2w-1-(p+q)}\frac{(-1)^{w-1-q}}{2\pi i}S\left(G^{-1}A,G^{-1}B\right).
\end{array}
$$
As $S\left(G^{-1}A,G^{-1}B\right)\neq 0$ only if $p+q\leq w-1$, we get
$P(es(A,\alpha+w-1-p),es(B,\beta+w-1-q))(z) \in z^w\cO_{\dC,0}$. On
the other hand, for any $A\in F^{p}H_{\arg\neq 0}$ there is
a $B\in F^{w-p-1}H_{\arg\neq 0}$ with $S(G^{-1}A,G^{-1}B) \neq 0$,
so that $P$ is nondegenerate. Sections coming from elements in $H_{\arg=0}$ are treated
similarly.
\end{proof}
We will call TERP-structures as above \emph{generated by elementary sections}.
The following lemma shows that this notion is consistent with the corresponding one
for twistor structures.
\begin{lemma}
Let $(H,H'_\dR, \nabla,P,w)$ be a TERP structure generated by
elementary sections. Then the bundle $\widehat{H}$ constructed
by the procedure from definition \ref{defHhat} is
a twistor generated by elementary sections as defined in lemma \ref{lemCorrElTwist}.
\end{lemma}
\begin{proof}
This is an immediate consequence of the following formula, which is
shown in \cite[(7.86)]{He2}.
\begin{eqnarray}
\tau \left(es(A, \alpha)\right) = es\left(\overline{A}, w-\overline{\alpha}\right).
\end{eqnarray}
\end{proof}
It is an easy exercise to check that the converse of lemma \ref{t5.5} holds,
so that we get the following extension of the correspondence of topological data.
Note that any TERP-structure generated by elementary sections canonically defines
a filtration $F^\bullet$ on $H^\infty$ by formula \eqref{5.5}.
\begin{lemma}\label{lemCorElTERP}
There is a one to one correspondence between TERP-structures generated by elementary sections
and tuples $(H^\infty, H_\dR^\infty, M, F^\bullet, S, w)$ as above such that the twisted filtration
$\widetilde{F}^\bullet$ satisfies condition \eqref{eqOrthogonal}.
\end{lemma}
The last step is now a (common) extension of the correspondence of lemma \ref{t5.4}
and of the last result putting together (mixed) twistors and polarizations. We cite
the following definition from \cite[definition 3.48]{Mo2}. We restrict here to weight
zero, which is what we need, but the definition extends to any weight.
\begin{definition}
Let $(\widehat{H}, \widehat{N}, \widehat{W}_\bullet)$ be a mixed twistor structure
(with $\widehat{W}_\bullet$ generated by $\widehat{N}$), suppose $\deg(\widehat{H})=0$,
and let $\widehat{S}:\cO(\widehat{H})\otimes_{\cO_{\dP^1}} \sigma^*\cO(\widehat{H})\rightarrow
\cO_{\dP^1}$ be a non-degenerate pairing. Then $\widehat{S}$ is called a polarization
if the following holds.
\begin{itemize}
\item
$\widehat{S}$ is a morphism of mixed twistors and satisfies $\widehat{S}(\widehat{N}-,-)
+S(-,\sigma^*\widehat{N}(-))=0$. This implies in particular that $\widehat{S}$
induces a morphism
$$
\begin{array}{rcl}
\widehat{S}_{k}:\left(\mathit{Gr}_k^{\widehat{W}}
\otimes\cO_{\dP^1}(-2k)\right)\otimes \left(\sigma^*\mathit{Gr}_k^{\widehat{W}}\right)
& \longrightarrow & \cO_{\dP^1}\\ \\
(a,b)&\longmapsto & \widehat{S}\left(\widehat{N}^ka,b\right).
\end{array}
$$
and then also a morphism
$$
\begin{array}{rcl}
i^{-k} \widehat{S}_{k}:\left(\mathit{Gr}_k^{\widehat{W}}\otimes\cO_{\dP^1}(-k)\right)\otimes
\sigma^*\left(\mathit{Gr}_k^{\widehat{W}}\otimes\cO_{\dP^1}(-k)\right) & \longrightarrow & \cO_{\dP^1}\\ \\
(a,b)&\longmapsto &i^{-k}\widehat{S}_k(a,b).
\end{array}
$$
The term $i^{-k}$ comes from the identification $\sigma^*\cO_{\dP^1}(-k)\cong\cO_{\dP^1}(-k)$.
\item
Let $\widehat{P}_k:={{\cK}\!er}\left(\widehat{N}^{k+1}:\mathit{Gr}^{\widehat{W}}_k\rightarrow
\mathit{Gr}^{\widehat{W}}_{-k-1}\right)$ be the primitive part of weight $k$. Then $i^{-k}\widehat{S}_k$ polarizes
$\widehat{P}_k$, i.e., for all global sections $a\in H^0(\dP^1, \widehat{P}_k\otimes \cO_{\dP^1}(-k))\backslash\{0\}$, $i^{-k}\widehat{S}_k(a,a)>0$.
\end{itemize}
$(\widehat{H}, \widehat{N}, \widehat{W}_\bullet, \widehat{S})$ is called a
polarized mixed twistor structure of weight zero (PMTS).
\end{definition}
The final result of this chapter is the following.
\begin{lemma}[Correspondence between PMHS and integrable PMTS]\label{lemCorPMHS-PMTS}
There is a one to one correspondence between tuples
$(H^\infty,H^\infty_\dR, M, F^\bullet, S, w)$ with the property that
$(H^\infty_{\arg\neq 0}, H^\infty_{\dR, \arg\neq 0}, -N, \widetilde{F}^\bullet, S)$ resp.
$(H^\infty_{\arg=0}, H^\infty_{\dR,\arg=0}, -N, \widetilde{F}^\bullet, S)$ are PMHS's of weight $w-1$ resp. $w$ and
integrable polarized mixed twistors $(\widehat{H}, H'_\dR, \widehat{W}_\bullet, \widehat{N}, \widehat{S}, \nabla)$
with flat real substructure on $\dC^*$ (and such that
$\widehat{W}_\bullet $ and $\widehat{N}$ are induced by the monodromy $\nabla$) which are generated by
elementary sections.
The eigenvalues of $M$ resp. of the monodromy have absolute value one.
\end{lemma}
\begin{proof}
Let us first show the last statement:
Suppose that
$(H^\infty, H_\dR^\infty, -N, \widetilde{F}^\bullet, S)$ is a sum of two
PMHS's. The semi-simple part of the monodromy
acts on the primitive part of each $Gr^W_kH^\infty$ and respects
the hermitian form induced by $S$ on these spaces.
This forces its eigenvalues to be in $S^1$.
Note that the space $Gr^W_kH^\infty$ is decomposed according to formula
\eqref{eqDecomPrimMHS} which implies that each eigenspace has a primitive part.
On the other hand, for an integrable polarized mixed twistor $(\widehat{H},H'_\dR,
\widehat{W}_\bullet, \widehat{N}, \widehat{S}, \nabla)$
as above it is easy to see that the semi-simple part of the monodromy gives
rise to a bundle map $\widehat{M}_s\in\mathit{Aut}(\widehat{H})$ compatible
with the weight filtration $\widehat{W}_\bullet$ and with the pairing
$\widehat{S}$. The induced map on
$H^0(\dP^1, \widehat{P}_k\otimes \cO_{\dP^1}(-k))$ therefore respects the
positive definite hermitian form $i^{-k}\widehat{S}_k$. As before,
we can conclude that its eigenvalues are of absolute value  one.
For the remaining part of this proof, we can therefore assume
that the various logarithms of the monodromy eigenvalues are real numbers.

Remark that in order to define MHS on $H^\infty$ one can work with $\widetilde{F}^\bullet$
as well as with $F^\bullet$, because these two filtrations coincide by definition
on the quotients $Gr^W_\bullet$. Therefore, we know by lemma \ref{t5.4} that $\widehat{H}$ is a mixed
twistor iff the tuples $(H^\infty_{\neq 1}, (H^\infty_\dR)_{\neq 1}, \widetilde{F}^\bullet, W_{\bullet+w-1})$ resp.
$(H^\infty_{1}, (H^\infty_\dR)_1, N, \widetilde{F}^\bullet, W_{\bullet+w})$ are mixed Hodge structures.
The orthogonality of $\widetilde{F}^\bullet$ with respect to $S$
is equivalent to the fact that $\widehat{S}$ takes values in $\cO_{\dP^1}$
(this follows from lemma \ref{lemCorElTERP} and lemma \ref{t3.14}),
so that the correspondence is true iff the positive definiteness
properties of the two polarizations are equivalent.
On the one hand, we have PMHS's of weights $w-1$ resp. $w$ iff
$$
i^{p-(w-1+[\alpha]+k-p)}S(A,(-N)^k\overline{A}) > 0
$$
for all non-vanishing classes
$$
[A]\in \left(\widetilde{F}^p\cap
\overline{\widetilde{F}^{w-1+[\alpha]+k-p}}\right)\left(Gr_k^W H^\infty_{e^{-2\pi i \alpha}}\right)_{\mathit{prim}}=
\left(F^p\cap
\overline{F^{w-1+[\alpha]+k-p}}\right)\left(Gr_k^W H^\infty_{e^{-2\pi i \alpha}}\right)_{\mathit{prim}}
$$
i.e., all $[A]\neq 0 $ in  that space with $N^{k+1}[A]=0$.

On the other hand, $(\widehat{H}, \widehat{W}_\bullet, \widehat{N},\widehat{S})$ is
a PMTS iff $i^{-k}\widehat{S}_k$ is positive on  classes $[es(A,\alpha+w-1-p)]\in \widehat{P}_k$, where
again $A\in W_kH^\infty_{e^{-2\pi i \alpha}}$ with $[A]\in \left(F^p\cap
\overline{F^{w-1+[\alpha]+k-p}}\right)\left(Gr_k^W H^\infty_{e^{-2\pi i \alpha}}\right)_{\mathit{prim}} \backslash\{0\}$.
Therefore we need to show that for all these elements $A$, $i^{p-(w-1+[\alpha]+k-p)}S(A,(-N)^k\overline{A}) > 0$
is equivalent to $i^{-k}\widehat{S}(\widehat{N}^k es(A,\alpha+w-1-p), es(A,\alpha+w-1-p))>0$. This is a consequence
of the following computation.
$$
\begin{array}{c}
i^{-k}\widehat{S}\left(\widehat{N}^k es(A,\alpha+w-1-p), es(A,\alpha+w-1-p)\right) \\ \\
=i^{-k}z^{-w}P\left(es(N^k A, \alpha+w-1-p),es(\overline{A},p+1-\alpha)\right)\\ \\
=i^{-k}(-1)^{p+1}(-z)^{-1-[\alpha]}P\left(es(N^kA,\alpha),es(\overline{A},1+[\alpha]-\alpha)\right)=
\frac{i^{-k}(-1)^p}{(2\pi i)^{w-1+[\alpha]}} S\left(G^{-1}N^k A, G^{-1}\overline{A}\right)\\ \\
=\frac{i^{p-\left(w+k-1+[\alpha]-p\right)}}{(2\pi)^{w-1+[\alpha]}} S\left(G^{-1} A, (-N)^kG^{-1}\overline{A}\right)
=\frac{1}{\Gamma\left(1+[\alpha]-\alpha\right)\Gamma(\alpha)}\cdot
\frac{i^{p-\left(w+k-1+[\alpha]-p\right)}}{(2\pi)^{w-1+[\alpha]}} \cdot
S\left(A, (-N)^k \overline{A}\right) .
\end{array}
$$

The last equality is a result of the following three facts, which are deduced directly from the definitions.
$$
G^{-1}A\equiv \frac1{\Gamma(\alpha)} A \;\;\;\mbox{mod}\;\;\;W_{k-1}
\;\;\;\; ; \;\;\;\;\;
G^{-1}\overline{A}\equiv \frac1{\Gamma(1+[\alpha]-\alpha)} A \;\;\;\mbox{mod}\;\;\;W_{k-1}
\;\;\;\; ; \;\;\;\;\;
S(W_k,W_{-k-1})=0
$$
This shows the equivalence of the positive definiteness properties of the two polarizations.
\end{proof}

\section{Regular singular TERP-structures}\label{c6}
\setcounter{equation}{0}

In this chapter we investigate more closely
the case of TERP-structures $(H, H'_\dR,\nabla, P, w)$ and nilpotent orbits
of them which are \emph{regular singular}, i.e.,
such that $\nabla$ has a regular singularity on $\cO(H)$ at zero (see the definition below).
In this case we will obtain a direct generalization of Schmid's correspondence
(theorem \ref{theoSchmidCorres}). This generalization will consist in a
correspondence between two types of TERP-structures: those inducing nilpotent
orbits and those which define polarized mixed Hodge structures on $H^\infty$ (theorem \ref{theoMainResultRegSing}).
To start with, we recall the definition of some well-known objects
associated to regular singular connections, namely, the Kashiwara-Malgrange filtration
and the Deligne lattices $V^\alpha$. More precisely,
let $(H',\nabla)$ be a flat bundle over $\dC^*$. Remember the definition of
the subspaces $C^\alpha\subset \cO(H')$ from chapter \ref{c5}.
Recall also that we have chosen a total order on $\dC$ which will be used below.
\begin{definition}
Put
$$
\begin{array}{c}
V^\alpha:=\sum\limits_{\beta \geq \alpha } \cO_\dC \cdot C^\beta
=\bigoplus\limits_{\alpha \leq \beta<\alpha +1} \cO_\dC \cdot C^\beta \\ \\
V^{>\alpha}:=\sum\limits_{\beta>\alpha } \cO_\dC \cdot C^\beta
=\bigoplus\limits_{\alpha<\beta \leq \alpha+1} \cO_\dC \cdot C^\beta \\ \\
V^{>-\infty} :=\sum\limits_\beta\cO_\dC \cdot C^\beta
\end{array}
$$
$V^\alpha$ and $V^{>\alpha}$ are locally free $\cO_\dC$-modules
whereas $V^{>-\infty}$ is locally free over $\cO_\dC[z^{-1}]$, i.e., a meromorphic bundle.
Obviously, all of them are subsheaves of $i_*\cO(H')$. The decreasing filtration
defined by $V^\alpha$ on $V^{>-\infty}$ is the Kashiwara-Malgrange or
V-filtration. There is a natural
identification $Gr_V^\alpha:=V^\alpha/V^{>\alpha}\cong C^\alpha$.
\end{definition}

The following definition introduces the notions of regular singular connections and explains some
important object and invariants attached to them.
The idea is essentially due to Varchenko (\cite{Va1}, \cite{SchSt}) who used it to define
a filtration making up a mixed Hodge structure on the cohomology of the Milnor fibre of an isolated
hypersurface singularity (see corollary \ref{corFiltRegSing}).

\begin{definition}\label{defHel}
\begin{itemize}
\item
Let $H$ be a bundle over $\dC$ and $\nabla$ a meromorphic connection with pole at zero.
Then $\nabla$ is said to have a regular singularity at zero if $\cH$ is a subsheaf of $V^{>-\infty}$, in other words, if the
sections in $\cH$ have moderate growth at the origin when expressed in a basis of flat sections.
A TERP-structure $(H,\nabla,H'_\dR,P,w)$ is called regular singular if
$(H,\nabla)$ is so. For a regular singular connection, the V-filtration
induces a filtration on $\cH$ which we denote by the same letter.
\item
Any section $\omega\in\cH_0$ can be written as a (possibly infinite) sum
$$
s=\sum_{\beta \geq \alpha} s(\omega, \beta)
$$
where $s(\omega, \beta)\in C^\beta$ and $s(\omega,\alpha)\neq 0$. $\alpha$
is called the order of $\omega$ and the section
$s(\omega, \alpha)$ the principal part of $\omega$.
\item
Let $\omega_1,\ldots,\omega_n$ be a set of generating sections
for $\cH_0$ with linearly independent principal parts $s(\omega_i,\alpha_i)$.
Put
$$
\cH^{el}:=\oplus_i\cO_\dC \cdot s(\omega_i,\alpha_i)
\cong \sum_\alpha\cO_\dC \cdot Gr^\alpha_V\cH.
$$
Then $\cH^{el}\subset V^{>-\infty}$ defines a vector bundle $(H^{el},\nnn)$ over
$\dC$ generated by elementary sections.
\item
Define the spectrum $\mathit{Sp}(H,\nabla)$
as a ``subset of $\dC$ with multiplicities'', i.e.:
$$
\mathit{Sp}(H,\nabla) = \sum_{\alpha \in\dC}\nu(\alpha)\alpha\in\dZ[\dC]
\;\;\;
\mbox{with}
\;\;\;
\nu(\alpha):=\dim_\dC\left(\frac{Gr^\alpha_V \cH}{Gr^\alpha_V z\cH}\right)
$$
We order the (possibly repeated) spectral numbers $\alpha_1, \ldots, \alpha_{\rank(H)}$
by the total order on $\dC$ defined in the beginning.
(In \cite[7.2]{He2} the spectral numbers are called exponents.)
\end{itemize}
\end{definition}

\begin{lemma}\label{lemHelTERP}
Let $(H,H'_\dR,\nabla,P,w)$ be a regular singular TERP-structure.
Then $(H^{el},H'_\dR,\nabla,P,w)$ is also a TERP-structure, generated by
elementary sections, therefore also regular singular.
We denote by $(\widehat{H^{el}}, \widehat{W}_\bullet, \widehat{N}, \widehat{S})$
the corresponding twistor.
The spectral numbers obey the symmetry $\alpha_i+\alpha_{\rank(H)+1-i}=w$.
\end{lemma}

\begin{proof}
The first thing to show is that the connection has still a pole of order
at most two on $\cH^{el}$: This follows directly from the definition
of the V-filtration. Namely, let $\omega\in V^\alpha\cap\cH$ be a section having
principal part $s(\omega,\alpha)$. If
$z\nabla_{z\partial_z}s(\omega,\alpha)\in C^{\alpha+1}$ is non-zero,
then it must be the principal part of $z\nabla_{z\partial_z}\omega\in V^{\alpha+1}\cap\cH$.
This proves that $z\nabla_{z\partial_z}s(\omega,\alpha)\in\cH^{el}$ as required.
The main point now is to see that the pairing $P$ sends $\cH^{el}$ to
$z^w\cO_\dC$ and that $z^{-w}P$ is nondegenerate on $\cH^{el}$. Then we know that $(H^{el},H'_\dR,\nabla,P,w)$
is a TERP-structure. By definition it is generated by elementary sections. This implies
its regularity as elementary sections have moderate growth.
The symmetry property of the spectrum will come out of the proof as a by-product.

First note that because $z^{-w}P$ is nondegenerate on $\cH$,
we have an isomorphism $(z^{-w}\cH,\nabla)\cong j^*(\cH^*,\nabla^*)$, where
$\cH^*:={\cH\!om}_{\cO_\dC}(\cH,\cO_\dC)$ is the dual module with the induced connection
and $j:z\mapsto -z$ as before. This implies that
$
Gr_V^\alpha (\cH/z\cH)\cong Gr_{V^*}^{\alpha-w}(\cH^*/z\cH^*)$.
On the other hand, for any bundle $(H,\nabla)$ with regular singular connection, we have that
$Gr_{V^*}^\alpha (\cH^*/z\cH^*) \cong \left(Gr_V^{-\alpha}(\cH/z\cH)\right)^*$.
This can be shown as in \cite[Remark 3.6]{Sa2} or
\cite[III.1.18]{Sa4} (Note that the indices must be shifted by one compared to loc.cit.).
Putting these two equations together, we conclude that $P$ induces a non-degenerate pairing
\begin{equation}\label{eqPnonDeg}
P:
Gr_V^\alpha\frac{\cH}{z\cH}
\otimes
Gr_V^{w-\alpha}\frac{\cH}{z\cH}\longrightarrow z^w\dC ,
\end{equation}
which gives, by the identification
$$
\sum_a \cO_\dC Gr_V^\alpha \frac{\cH}{z\cH}\cong
\frac{\sum_\alpha \cO_\dC Gr_V^\alpha\cH}{\sum_\alpha z\cO_\dC Gr_V^\alpha \cH} \cong\frac{\cH^{el}}{z\cH^{el}}
$$
the desired non degenerateness of $z^{-w}P$ on $\cO(H^{el})$. In particular,
\eqref{eqPnonDeg} implies that
$
\dim_\dC\left(Gr_V^\alpha\frac{\cH}{z\cH}\right)
=\dim_\dC\left(Gr_V^{w-\alpha}\frac{\cH}{z\cH}\right),
$
which is precisely the symmetry property of the spectrum.
\end{proof}
The next result is a rather trivial but important consequence of the last lemma.
\begin{corollary} \label{corFiltRegSing}
Let $(H, H'_\dR, P, w)$ be a TERP-structure. Then the formula
\begin{eqnarray}
F^pH^\infty_{e^{-2\pi i \alpha}}:=\psi_\alpha^{-1}\left(z^{p+1-w}Gr_V^{\alpha+w-1-p}\cH\right)
\;\;\;\;\;\;\;\;\;\mbox{with }\;\;\Re(\alpha)\in(0,1]
\end{eqnarray}
defines a filtration $F^\bullet$ on the space $H^\infty$ coinciding
with the one obtained from $(H^{el}, H'_\dR, \nabla, P, w)$ using
lemma \ref{lemCorElTERP},
i.e.,
$F^pH^\infty_{e^{-2\pi i \alpha}} = \psi_\alpha^{-1}\left(z^{p+1-w}(C^{\alpha+w-1-p}\cap\cH^{el})\right)$.
In particular, the orthogonality conditions \eqref{eqOrthogonal}
for the twisted filtration $\widetilde{F}^\bullet$ are satisfied,
so that $\widetilde{F}^\bullet$ gives an element in $\check{D}$, the classifying space
of Hodge-like filtrations on $H^\infty$ from chapter \ref{c1}.
\end{corollary}
\textbf{Remark:}
In case a TERP-structure arises as Fourier-Laplace
transform of the Brieskorn lattice of an isolated hypersurface singularity (see chapter \ref{c11}),
$\widetilde{F}^\bullet$ is the filtration defined by Steenbrink \cite{SchSt} as part of
a polarized mixed Hodge structure on the cohomology of the Milnor fibre (the space $H^\infty$).

\bigskip

The following lemma motivates the introduction
of nilpotent orbits of TERP-structures. Let
$(G,G'_\dR,\nabla,P,w)$ be a variation of TERP-structures on $\Delta^*$
such that $\cO(G)$ is stable under $\nabla_{z\partial_z-r\partial_r}$.
First it follows that if
the restriction $G_{\dC\times\{r\}}$ to any $r\in \Delta^*$ is regular singular,
then the whole bundle is regular singular along $z=0$. In this situation,
consider the space $G^\infty$ of flat multivalued sections of $G':=G_{|\dC^*\times\Delta^*}$.
Then $\widetilde{F}^\bullet$ induce a family of filtrations on $G^\infty$
parameterized by the universal cover $\dH$ of $\Delta^*$.

\begin{lemma}\label{lemOrbTERPOrbHS}
This family
satisfies $\widetilde{F}^\bullet(\rho) = \exp(-\rho N)\widetilde{F}^\bullet(0)$
(with $\rho$ a coordinate on $\dH$ such that $r=e^{2\pi i\rho}$). In particular, it
yields a holomorphic map $\phi:\dH\rightarrow \check{D}$ into the classifying space, and this
map is a nilpotent orbit of Hodge structures
with nilpotent endomorphism $-N$ iff $\widetilde{F}^\bullet(\rho)\in D$
for $\Im(\rho)\gg0$.
\end{lemma}
\begin{proof}
The first statement is shown in \cite[theorem 7.20]{He2}. It implies that
the dimensions of the various $\widetilde{F}^p$'s are constant in $\rho$ and
by the last corollary we know that $\widetilde{F}^\bullet(\rho)$ satisfies the orthogonality conditions
for any $\rho$. This proves the second statement.
\end{proof}
The following theorem is the main result of this chapter. It generalizes
the correspondence \ref{theoSchmidCorres}.
\begin{theorem}\label{theoMainResultRegSing}
Let $(H,\nnn,H'_\dR,P)$ be a regular singular TERP-structure.
The following two conditions are equivalent.
\begin{enumerate}
\item
$(H,\nnn,H'_\dR,P)$ induces a nilpotent orbit.
\item
$(H^\infty, H_\dR^\infty, -N, S, \widetilde{F}^\bullet)$ defines
a PMHS of weight $w-1$ resp. $w$ on $H^\infty_{\arg\neq 0}$ resp.
on $H^\infty_{\arg=0}$.
\end{enumerate}
\end{theorem}
Before entering into the proof, let us make some comments on
the content of this result. Suppose that 2.) holds, then the family $\widetilde{F}^\bullet(\rho)$
will eventually end up in the interior $D$ of the classifying space
by Schmid's correspondence, which by chapter \ref{c5} amounts to say that $(\widehat{G^{el}}, \widehat{S})$
is a family of pure polarized twistors. However, it is by no means obvious that the same holds for
$(\widehat{G}, \widehat{S})$. It was shown in \cite[theorem 7.20]{He2}
that this is indeed the case, mainly because the two objects
$\widehat{G^{el}}$ and $\widehat{G}$ tend ``to each other'' in a suitable
sense when $r$ approaches the origin. Therefore  2.) --> 1.) of the above
theorem is precisely what has been shown in \cite[theorem 7.20]{He2}. We need to prove the converse.
The same difficulty occur: From the fact that $G$ underlies
a variation of polarized pure TERP-structure we know (almost by definition) that
$(\widehat{G},\bD,\widehat{S})$ is a VPTS,
but it is not at all clear a priori that this is also true for
$(\widehat{G^{el}},\bD,\widehat{S})$. The proof will actually be quite different:
We will use one of the main results of \cite{Mo2} to obtain from the variation
$G$ a ``limit object'' which we can identify with the twistor
$(\widehat{H^{el}}, \widehat{W}_\bullet, \widehat{N}, \widehat{S})$
constructed from our original TERP-structure. Mochizuki's theorem says that
this limit is a PMTS so that by the correspondence of lemma \ref{lemCorPMHS-PMTS} we can conclude that
$(-N,\widetilde{F}^\bullet)$ gives rise to a (sum of) PMHS on $H^\infty$.

In order to state the next theorem, recall
from definition \ref{t3.18}
that any VPTS $(\widehat{G}, \bD, \widehat{S})$
of weight zero on a base space $M$
gives rise to the structure of a harmonic bundle on
$E=p_*\cC^{\infty h}(\widehat{G})$, where $p:\dP^1\times M\rightarrow M$ is the projection.
If $M$ a punctured disk, Mochizuki constructs
(following the idea in \cite[section 6]{Si5}) a ``limit'' twistor
(i.e., a vector bundle on $\dP^1\times\{0\}$).
It should be noticed that
in the one dimensional situation we are looking at, the limit object
is already considered in the work of Simpson. Mochizuki actually works in
a much more general situation, namely, he considers harmonic bundles in any dimensions
defined on the complement of a normal crossing divisor. As we will see, the nilpotent orbit assumption also
restricts the class of harmonic bundles which can occur.

To be more precise, Mochizuki defines for any $a\in(0,1]$ a tuple
$$
\left(S_{(a,0)}^{can}(E), \widehat{\cN}, \widehat{\cW}_\bullet,\widehat{\cS}\right)
$$
where
$S_{(a,0)}^{can}(E)
\in \mathit{VB}_{\dP^1\times\{0\}}$,
$\widehat{\cN}$ is a nilpotent bundle endomorphism generating a weight filtration $\widehat{\cW}_\bullet$
and $\widehat{\cS}$ is a pairing as the one discussed in chapter \ref{c5}.
Theorem 12.1 in \cite{Mo2} then states that
$
\left(S_{(a,0)}^{can}(E), \widehat{\cN}, \widehat{\cW}_\bullet,\widehat{\cS}\right)$
is a polarized mixed twistor structure.
Therefore, the proof of theorem \ref{theoMainResultRegSing}
reduces by lemma \ref{lemCorPMHS-PMTS} to the following comparison result.
\begin{theorem}\label{theoMainCompMochizuki}
Consider the harmonic bundle $E:=p_*\cC^{\infty h}(\widehat{G})$ on $\Delta^*$,
where $G:=\pi^*H$ is a nilpotent orbit of TERP-structures.
Then $\widehat{H^{el}}\in\mathit{VB}_{\dP^1}$ is canonically isomorphic
to the bundle $\bigoplus_{0<a \leq 1} S_{(a,0)}^{can}(E)$
and this isomorphism identifies $\widehat{W}_\bullet, \widehat{S}$
from lemma \ref{lemHelTERP} with the objects $\widehat{\cW}_\bullet,\widehat{\cS}$
and sends $\widehat{N}$ to $2\pi \cdot\widehat{\cN}$.
\end{theorem}

{\bf Proof:}

We will mainly have to adapt the proof of
\cite[theorem 12.1]{Mo2} to our situation. In fact, as the setup considered
by Mochizuki is much more general, everything will simplify to a large extent
but still the abundance of objects and notations in \cite{Mo2} makes
this translation somewhat painful. We will proceed in three steps.

\bigskip
\textbf{(I) Canonical extensions over $r=0$ and $(z,r)$-elementary sections}

Let us consider for a moment the restriction $G':=G_{|\dC^*\times \Delta^*}$ of
our family of TERP-structures. We will use a slight generalization
of the constructions of the beginning of chapter \ref{c5} to the
case of a flat bundle on the complement of a normal crossing divisor. The divisor
here is just $\{r=0\}\cup\{z=0\}\subset\dC\times\Delta$, our bundle is $G'$
and as has been shown, the two monodromies are equal.
For each fixed $z$, we can consider the restriction
$G_{|\{z\}\times \Delta^*}$ which is flat and we have a space of multivalued flat sections
of this bundle which comes equipped with a monodromy operator.
If $z$ varies, they patch together to a flat bundle $\cH^\infty$ on $\dC^*\times\{0\}$,
on which the horizontal monodromy acts. We denote by $\cH^\infty_\lambda$ the
flat generalized eigenbundle for this monodromy. It is obviously the same
as the flat subbundle for the "intrinsic" vertical monodromy of $\cH^\infty$.
For any fixed $r$, $(\cH^\infty, \nabla)$ can be identified
with the restriction $(G'_{|\dC^*\times\{r\}}, \nabla_z)$, in particular, $(\cH^\infty, \nabla)$
is isomorphic to $(H',\nabla)$, the restriction over $\dC^*$ of the original TERP-structure.

The bundle $\cH^\infty$ can also be defined more intrinsically by putting $\cH^\infty:=\psi_r (G'^\nabla)$, where
$G'^\nabla$ is the local system on $\dC^*\times \Delta^*$ corresponding to the flat bundle $G'$,
and $\psi$ is the functor of nearby cycles of Deligne (see \cite{Di}). In this simple situation, the complex
$\psi_r (G'^\nabla)$ has cohomology only in degree zero. This implies directly that $\cH^\infty$
is a local system on $\{r=0\}=\dC^*\times\{0\}$ with monodromy,
and by abuse of notation we denote by $(\cH^\infty,\nabla)$
the flat bundle corresponding to it and by $\cH^\infty_\lambda$ the eigenbundle of the monodromy
operator. Note that as we have a canonical identification of
local systems
$$
(\overline{\gamma^*G'})^{\gamma^*\nabla} \cong (G')^\nabla
$$
we can also describe $\cH^\infty$ by
$\psi_r((\overline{\gamma^*G'})^{\gamma^*\nabla})$. This fact will be used in a moment.
As before in the absolute situation, we define for any
$\alpha\in\dC$ by $\cV^\alpha$ the bundle
$$
\cV^\alpha = \bigoplus_{\alpha\leq\beta <\alpha+1} \cO_{\dC^*\times\Delta} r^{\beta\mathit{Id}-\frac{N}{2\pi i}}
{\cH^\infty_{e^{-2\pi i \beta}}}
$$
By definition, $\cV^\alpha$ is the unique $\cO_{\dC^*\times\Delta}$-free extension of $\cO(G')$ having
a logarithmic pole along $\dC^*\times\{0\}$ with residue eigenvalues $\beta$ such that
$\alpha\leq\beta<\alpha+1$.
There is a natural isomorphism
$$
\cH^\infty_{e^{-2\pi i \alpha }}\longrightarrow Gr^\alpha_\cV :=\cV^\alpha/\cV^{>\alpha}
$$
where, as before $\cV^{>\alpha}:=\bigcup_{\beta>\alpha}\cV^\beta$.
Similarly, one defines
$$
\overline{\cV}^\alpha:=
\bigoplus_{\alpha\leq\beta <\alpha+1} \cO_{\dC^*\times\overline{\Delta}} \overline{r}^{\overline{\beta}\mathit{Id}+\frac{N}{2\pi i}}
{\cH^\infty_{e^{2\pi i \overline{\beta}}}}
$$
yielding an isomorphism $\cH^\infty_{e^{2\pi i \overline{\alpha} }}\rightarrow Gr^{\alpha}_{\overline{\cV}}$.
In the sequel, we need an explicit description of generating sections of $G$, which
is provided by the following lemma.
\begin{lemma}
Consider a basis of the original TERP-structure $H$ given by sections which
we decompose into elementary sections, i.e.:
$$
\cO(H) = \bigoplus_{i=1}^{\rank(H)} \cO_\dC\left(\sum_{j\in \dN} es(A_{ij}, \alpha_{ij})\right)
$$
with $A_{ij}\in H^\infty$ and $\alpha_{i1}<\alpha_{i2}<\ldots$ for all $i$.
Then we have
\begin{eqnarray}
\cO(G)=\bigoplus_{i=1}^{\rank(H)}\cO_{\dC\times\Delta^*}
\left(\sum_{j\in\dN}r^{\alpha_{ij}\mathit{Id}-\frac{N}{2\pi i}}es(A_{ij}, \alpha_{ij})\right)
\end{eqnarray}
\end{lemma}
\begin{proof}
This follows directly from the definition of $G=\pi^*H$, namely,
$$
\nabla_{z\partial_z-r\partial_r}
\left(
\sum_{j\geq 1}r^{\alpha_{ij}\mathit{Id}-\frac{N}{2\pi i}}es(A_i, \alpha_{ij})
\right)=0
$$
holds.
\end{proof}

\vspace*{1cm}

\textbf{(II) KMSS-spectrum of nilpotent orbits}

Let us recall in brief the notion of the parabolic filtration
associated to a VPTS on $\Delta^*$.
It is mainly due to Simpson (\cite{Si2}), but we need a slightly more
precise version as in \cite{Mo2} (see also \cite[Corollary 5.3.1]{Sa6}).
\begin{deftheo}
\label{deftheoKMSS}
Let $(E, \overline{\partial},\theta ,h)$ be a tame harmonic bundle on $\Delta^*$, take the pullback
$E':=p^*E_{|\dC\times \Delta^*}$ and define
$\cE:={\cK\!er}\left(\overline{\partial}+z\overline{\theta}:\cC^{\infty h}(E')\rightarrow\cC^{\infty h}(E')
\otimes\cO_\dC \cA_{\Delta^*}^{0,1}\right)\in\mathit{VB}_{\dC\times\Delta^*}$. Consider the following extensions
to sheaves over $\dC\times\Delta$
(let $i:\dC\times\Delta^*\hookrightarrow\dC\times\Delta$):
$$
\begin{array}{c}
{_a\cE} := \left\{s\in i_*\cE\;|\;|s|_{p^*h} \in
O\left(|r|^{-\epsilon-a}\right)\;\;\forall \epsilon>0\right\} ,
 \\ \\
_*\cE := \bigcup_{a\in\dR} {_a\cE}=
\{s\in i_*\cE \;|\;|s|_{p^*h}\mbox{ has moderate growth along }\dC\times\{0\} \}.
\end{array}
$$
The increasing filtration
of $_*\cE$ by the subsheaves $_a\cE$ is called parabolic filtration.
We have that $r \cdot {_a \cE} \cong {_{a-1}\cE}$ endowing $_a\cE$ with
a $\cO_{\dC\times\Delta}$-module structure. The sheaf $_*\cE$ is a locally free
$\cO_{\dC\times\Delta}[r^{-1}]$-module and for any $z\in\dC$,
the restrictions ${_a\cE^z}:=j_z^*({_a\cE})$ are $\cO_{\Delta}$-locally free extension
of $j_z^* \cE$ over $(z,0)$ (here $j_z:\{z\}\times\Delta\hookrightarrow
\dC\times\Delta$). However, ${_a\cE}$ is not $\cO_{\dC\times\Delta}$-free in general.
The $z$-connection $(z\partial+\theta)$ sends $_a\cE^z$ to $_{a+1}\cE^z$, in particular
for any $z\in \dC^*$, $\partial+z^{-1}\theta$ has a logarithmic pole on $_a\cE^z$
Let $Gr^\cP_a (_*\cE^z):={_a\cE^z}/_{<a}\cE^z$ with $_{<a}\cE^z=\bigcup_{b<a} {_b\cE^z}$.
The residue of $z\partial+\theta$ at $r=0$ acts on the graded pieces. Denote by $\dE_\alpha
Gr^\cP_a (_*\cE^z)$ its generalized eigenspace decomposition and by $\cN^z_E$ its nilpotent
part. Define
$$
\mathit{KMSS}(\cE^z):=\left\{(a,\alpha)\in\dR\times \dC\,|\, \dE_\alpha Gr^\cP_a  (_*\cE^z) \neq \{0\}\right\}
$$
to be the \textbf{Kashiwara-Malgrange-Sabbah-Simpson spectrum} of $\cE^z$. The multiplicity
of $(a,\alpha)\in \mathit{KMSS}(\cE^z)$ is by definition $\dim_\dC \dE_\alpha Gr^\cP_a  (\cE^z)$.
In particular, $(a,\alpha)\in\mathit{KMSS}(\cE^0)$ only if $\alpha$ is an eigenvalue
of the residue endomorphism of the Higgs field $\theta$.
The KMSS-spectrum is determined by the reduced KMSS-spectrum
$$
\mathit{KMSS}_{red}(\cE^z):=\mathit{KMSS}(\cE^z)\cap (0,1]\times \dC,
$$
because it is invariant under the shift $(a,\alpha)\mapsto (a+1,\alpha-z)$.
\end{deftheo}
We quote the following properties of the KMSS-spectrum.
\begin{lemma}
\label{lemPropKMSS}
Let $(E,\overline{\partial}, \theta, h)$ be a tame harmonic bundle on $\Delta^*$.
\begin{enumerate}
\item
The KMSS-spectrum satisfies:
$$
(a,\alpha)\in\mathit{KMSS}(\cE^0)
\;\;\;\;\;\;
\mbox{iff}
\;\;\;\;\;\;
(a + 2\Re(z\cdot \overline{\alpha}),\alpha-a z-\overline{\alpha}z^2)\in \mathit{KMSS}(\cE^z)
\;\;\;\;\;\;\;\;\;\forall z\in \dC
$$
\item
The KMSS-spectrum
of the harmonic bundle $(E,\partial, \overline{\theta}, h)$ on the conjugate
complex manifold $\overline{\Delta^*}$ is given by
$$
\mathit{KMSS}(\overline{\cE}^z) =
\{
-a + 2\Re(z^{-1}\cdot \alpha),\overline{\alpha}+a z^{-1}-\alpha z^{-2})\;|\;(a,\alpha)\in\mathit{KMSS}(\cE^0)
\}
\;\;\;\;\;\;\;\;\;\forall z\in \dP^1\backslash\{0\},
$$
where $\overline{\cE}$ is constructed from $E$ as above, but using
the pullback $p^*E_{|(\dP^1\backslash\{0\})\times \overline{\Delta^*}}$.
\item
For any $(a,\alpha)\in\mathit{KMSS}(\cE^0)$, there is a unique holomorphic
bundle $\cG_{(a,\alpha)}$ on $\dC\times\{0\}$ whose restriction to $z\times \{0\}$ satisfies
$$
(\cG_{(a,\alpha)})_{|z} = \dE_{\alpha-a z-\overline{\alpha}z^2}\left(Gr^\cP_{a + 2\Re(z\cdot \overline{\alpha})}(_*\cE^z)\right)
$$
and a unique holomorphic bundle $\overline{\cG}_{(-a,\overline{\alpha})}$ on $(\dP^1\backslash\{0\})\times\{0\}$
whose restriction to $z \times \{0\}$ satisfies
$$
(\overline{\cG}_{(-a,\overline{\alpha})})_{|z} = \dE_{\overline{\alpha}+a z^{-1} -\alpha z^{-2}}\left(Gr^\cP_{-a + 2\Re(z^{-1} \cdot \alpha)}(\overline{_*\cE}^z)\right)
$$
for any $z\in \dP^1\backslash\{0\}$.
\end{enumerate}
\end{lemma}
\begin{proof}
The first statement is \cite[proposition 1.8 and corollary 7.71]{Mo2}.
The second point then simply reduces to the fact that $\mathit{KMSS}(\overline{\cE}^\infty)=
\left\{(-a,\overline{\alpha})\;|\;(a,\alpha)\in\mathit{KMSS}(\cE^0)\right\}$ which is obvious.
The third point is explained in detail in \cite[sections 1.3.6, 8.9.1 and 11.2.3]{Mo2}.
\end{proof}
We can apply this result to our more special situation. Remember that if
a harmonic bundle $E$ is given by $p_*\cC^{\infty h}(\widehat{G})$ where
$G$ underlies a variation of TERP-structures on $\Delta^*$, then the sheaf
$\cE$ from definition-theorem \ref{deftheoKMSS} is nothing else then $\cO(G)$ and
the connection $\partial+z^{-1}\theta$ on $\cE$ is canonically identified with the horizontal part
$\nabla_r:\cO(G)\rightarrow\cO(G)\otimes z^{-1}\cO_{\dC}\Omega^1_{\Delta^*}$ of $\nabla$.
In the same way, $(\overline{\cE}, \overline{\partial}
+z\overline{\theta})$ can be identified with
$(\cO(\overline{\gamma^*G}), \gamma^*\nabla_r)$.
\begin{lemma}\label{lemPropKMSSNilpOrb}
Consider a harmonic bundle $E$ constructed from a nilpotent orbit
of TERP-structures as above.
\begin{enumerate}
\item
The reduced KMSS-spectrum is $\mathit{KMSS}_{red}(\cE^z)=\left\{(a,-z\cdot a)\,|\, a \in L\subset(0,1]\right\}$,
where $L$ is finite.
\item
The eigenvalues of the monodromy of either $\nabla_{\partial_z}$ or $\nabla_{\partial_r}$
are equal to $e^{2\pi i a}$ with $a\in L$. In particular, they are elements of $S^1$.
\item
For any $a\in L$, the restriction $({_a}\cE)_{|\dC^*\times\Delta}$ coincides with
the sheaf $\cV^{-a}$.
\item
$_a\cE$ is $\cO_{\dC\times\Delta}$-free for all $a\in L$.
\end{enumerate}
\end{lemma}
\begin{proof}
\begin{enumerate}
\item
For a nilpotent orbit, it was
shown in lemma \ref{lemTameness} that the Higgs field is nilpotent iff the polar part
$\cU$ of the vertical connection $\nabla_z$ of our original TERP-structure is
nilpotent. By \cite[theorem II.4.1]{Sa4}, this is the case for a regular singular
connection, so that the eigenvalues of the Higgs field are zero. This means
that $\mathit{KMSS}_{red}(\cE^0)=\left\{(a,0)\,|\, a \in L\subset (0,1]\right\}$ so that
by the preceding lemma we get that $\mathit{KMSS}_{red}(\cE^z)=
\left\{(a,-z\cdot a)\,|\, a \in L\subset(0,1]\right\}$.
\item
It was already shown in lemma \ref{t4.2} that the nilpotent orbit condition
implies that vertical and horizontal monodromy coincide. We know that
$\nabla_{\partial_r}=\partial+z^{-1}\theta$
is logarithmic on any $_a\cE^z$, so that the eigenvalues of the monodromy are of the
form $e^{-2\pi i \gamma}$ where $\gamma$ is an eigenvalue of the residue endomorphism
of $\nabla_{\partial_r}$ acting on $_a\cE^z/(r\cdot{_a\cE^z})$. By the first point,
the eigenvalues of the residue of the $z$-connection $z\nabla_{\partial_r}=z\partial+\theta$
on $Gr^\cP_a(_*\cE^z)$
are of the form $-za$ with $a\in L$. Therefore the monodromy eigenvalues are precisely
the numbers $e^{2\pi i a}\in S^1$.
\item
It follows from 1. and the defining property of $\cV^{-a}$
that $_a\cE^z\cong \cV^{-a}_{|\{z\}\times \Delta}$ for all $z\in\dC^*$. This implies
$_a\cE_{|\dC^*\times\Delta}\cong\cV^{-a}$.
\item
On $\dC\times\Delta^*$, the sheaf $_a\cE$ coincides by definition with $\cE$
which is $\cO_{\dC\times\Delta^*}$-free.
As was just shown, on $\dC^*\times\Delta$ it is isomorphic to $\cV^{-a}$ which is obviously free over
$\cO_{\dC^*\times\Delta}$. Hence we only need to show that the germ
$_a\cE_{(0,0)}$ is a free $\cO_{\dC\times\Delta,(0,0)}$-module. Suppose first that $a\notin L+\dZ$,
then this follows from \cite[proposition 1.11]{Mo2}. If, however, $a\in L+\dZ$, then by
discreteness of $L+\dZ$ in $\dR$ we can choose an $\epsilon_0>0$ such that
$(a,a+\epsilon_0)\cap (L+\dZ)=\emptyset$. Then for any smaller $0< \epsilon <\epsilon_0$, we have
an equality
$$
{_{a+\epsilon}\cE}_{|\dC\times\Delta\backslash\{(0,0)\}} =
{_a\cE}_{|\dC\times\Delta\backslash\{(0,0)\}}.
$$
The sheaves $_{a+\epsilon}\cE$ are all free over $\cO_{\dC\times\Delta}$ by construction, and
a classical result (\cite{Se}) asserts that for any locally free sheaf defined on a complement
of a subvariety of codimension at least two, there is at most one locally free extension
to the entire space. This implies that for any two $\epsilon, \epsilon'\in(0,\epsilon_0)$ as above,
$_{a+\epsilon}\cE={_{a+\epsilon'}}\cE$ on the whole of $\dC\times\Delta$. But by definition,
if $s\in{_{a+\epsilon}\cE}$, then for any $\delta>0$, $|s|_{p^*h}\in O(|r|^{-a-\epsilon-\delta})$.
The last formula is true for any $\epsilon$ and $\delta$ because of $_{a+\epsilon}\cE={_{a+\epsilon'}}\cE$,
so that we already have $s\in{_a\cE}$. This gives $_{a+\epsilon}\cE={_a\cE}$ and thus the desired
freeness property.
\end{enumerate}
\end{proof}
\textbf{Remark:} It turns out that the last lemma holds
true in a more general context, namely, for any \emph{integrable} VPTS
over a curve with tame behavior at the singularities. In other
words, the nilpotent orbit condition is needed only to get tameness
of the corresponding harmonic bundle as in lemma
\ref{lemTameness}. This follows from the simple observation
that if we have an integrable VPTS, then the monodromy of
$\partial+z^{-1}\theta$ is constant in $z$ (isomonodromic situation), so that
the eigenvalues of the residue of this operator along $r=0$ on
$Gr^\cP_a(_*\cE)$ are equal
to $\frac{\alpha}{z}-a-\overline{\alpha}z$ which can be constant only if $\alpha=0$.
Therefore $\theta$ is nilpotent which implies all other statements of the lemma.

\vspace*{1cm}

\textbf{(III) Identification of twistors}

The main playing characters of this last part of the proof
are the bundles $\cG_{(a,\alpha)}\in\mathit{VB}_{\dC}$ and
$\cG_{(-a,\overline{\alpha})}\in\mathit{VB}_{\dP^1\backslash\{0\}}$
defined above.
According to lemma \ref{lemPropKMSSNilpOrb}, only pairs $(a,0)$ (resp. $(-a,0)$) occur
as indices in our situation, therefore we denote the corresponding
bundles simply by $\cG_a$ resp. $\overline{\cG}_{-a}$.
The following lemma gives an explicit description of theses objects.
\begin{lemma}
Let $E$ as above be constructed from a nilpotent orbit $G$ of TERP-structures.
Then we can consider the graded pieces $Gr^\cP_a(_*\cE)$ globally on $\dC$ and there
is an identification $\cG_a\cong Gr^\cP_a({_*}\cE)$
for any $a\in L$. Similarly, $\overline{\cG}_{-a}\cong Gr^\cP_{-a}(_*\overline{\cE})$.
It follows that the restriction $(\cG_a)_{|\dC^*}$
is canonically isomorphic to $Gr_\cV^{-a}\cong\cH^\infty_{e^{2\pi i a}}$
and that $(\overline{\cG}_{-a})_{|\dC^*}\cong Gr_{\overline{\cV}}^a\cong\cH^\infty_{e^{2\pi i a}}$.
\end{lemma}
\begin{proof}
The quotients $Gr^\cP_a({_*}\cE)$ are bundles since $_a\cE$ is locally free.
For any fixed $z\in \dC$, the third part of lemma \ref{lemPropKMSS} shows that
$(\cG_{(a,0)})_{|z}=\dE_{-za}\left(Gr^\cP_a(_*\cE^z)\right)$. But the fact that
only pairs $(a,-za)$ occur as elements of $\mathit{KMSS}(\cE^z)$ shows that
$\dE_{-za}\left(Gr^\cP_a(_*\cE^z)\right)=Gr^\cP_a(_*\cE^z)$. The same arguments applied
to the sheaves $\overline{\cE}^z$ yield $\overline{\cG}_{-a}\cong Gr^\cP_{-a}(_*\overline{\cE})$.
Concerning the second statement, the identification $(_a\cE)_{|\dC^*\times\Delta} \cong \cV^{-a}$
from lemma \ref{lemPropKMSSNilpOrb} shows $(\cG_a)_{|\dC^*}\cong Gr_\cV^{-a}$. Similarly one
checks that $(_{-a}\overline{\cE})_{|\dC^*\times\Delta} \cong \overline{\cV}^a$
yielding $(\overline{\cG}_{-a})_{|\dC^*}\cong Gr^a_{\overline{\cV}}$.
\end{proof}
\begin{corollary}
The morphism $\Phi^{can}_{(a,0)}$ from \cite[section 10.4.1]{Mo2} can be expressed as
the composition of the isomorphisms
$$
\cH^\infty_{e^{2\pi i a}} \longrightarrow
Gr_\cV^{-a} \longrightarrow \left(Gr_a^\cP(_*\cE)\cong \cG_a\right)_{|\dC^*}
$$
Similarly, the morphism $\Phi^{\dag\;can}_{(-a,0)}$ is given by
$$
\cH^\infty_{e^{2\pi i a}} \longrightarrow
Gr_{\overline{\cV}}^a \longrightarrow \left(Gr_{-a}^\cP(_*\overline{\cE})\cong \overline{\cG}_{-a}\right)_{|\dC^*}.
$$
\end{corollary}
\begin{proof}
This is a direct consequence of the last lemma as soon as we know that
the filtration $\cF_\bullet$ on $\cH^\infty$ considered by Mochizuki
is trivial in our situation which follows
from the behavior of the graded pieces of this filtration (the ``KMSS-spectrum
of flat multivalued sections''), see
\cite[proof of lemma 2.4, sections 7.4.2 and 9.1.5]{Mo2}.
\end{proof}

The next result is the central step in the comparison of Mochizuki's objects
with the ones we are interested in. Note that by the last lemma,
we can consider $\cG_a$ as a subsheaf of $i_*\cH^\infty_{e^{2\pi i a}}$
and $\overline{\cG}_{-a}$ as a subsheaf of $\widetilde{i}_*\cH^\infty_{e^{2\pi i a}}$,
where, as before, $i:\dC^*\hookrightarrow \dC$ and $\widetilde{i}:\dC^*\hookrightarrow \dP^1\backslash\{0\}$.
\begin{lemma}\label{lemIdentG-Hel}
Consider an identification $\cH^\infty_{e^{2\pi i a}}\longrightarrow \cO(H'_{e^{2\pi i a}})$.
Under this isomorphism,
$\cG_a$ is mapped to
$$
\sum_{k\in\dZ} \cO_\dC \left(Gr^{k-a}_V \cO(H)\right)
\subset i_*\cO(H'_{e^{2\pi i a}})
$$
Therefore, $\oplus_{a \in L} \cG_a \cong \cO(H^{el})$.
Similarly,  $\oplus_{a \in L} \overline{\cG}_{-a}
\stackrel{\cong}{\longrightarrow}
\cO(\widehat{H}^{el}_{|\dP^1\backslash\{0\}})=\cO(\overline{\gamma^*H^{el}})\subset\widetilde{i}_*\cO(H')$.
\end{lemma}
\begin{proof}
As was shown in part (I)
$$
\cO(G)=\bigoplus_{i=1}^{\rank(H)}\cO_{\dC\times\Delta^*}
\left(\sum_{j\geq 1}r^{\alpha_{ij}\mathit{Id}-\frac{N}{2\pi i}}es(A_{ij}, \alpha_{ij})\right)
$$
Lemma \ref{lemPropKMSSNilpOrb} gives
$$
_a\cE=\bigoplus_{i=1}^{\rank(H)}\cO_{\dC\times\Delta}
\left(\sum_{j\geq 1}r^{\left(\alpha_{ij}-[\alpha_{i1}+a]\right)\mathit{Id}-\frac{N}{2\pi i}}es(A_{ij}, \alpha_{ij})\right)
$$
so that
$$
\cG_a\cong Gr_a(_*\cE)=
\bigoplus_{i:\alpha_{i1}+a\in\dZ}\cO_\dC
r^{-a\mathit{Id}-\frac{N}{2\pi i}}es(A_{i1}, \alpha_{i1})
$$
If $\alpha_{i1}+a\in\dZ$ then
$es(A_{i1}, \alpha_{i1})\in \oplus_{k\in\dZ}Gr^{k-a}_V\cO(H)$. On the other hand,
the composed map
$$
(\cG_a)_{|\dC^*}\longrightarrow Gr_\cV^{-a}\longrightarrow\cH^\infty_{e^{2\pi i a}}\longrightarrow \cO(H'_{e^{2\pi i a}})
$$
sends $r^{-a\mathit{Id}-\frac{N}{2\pi i}}es(A_{i1}, \alpha_{i1})$ to $es(A_{i1}, \alpha_{i1})$
which shows the first part of the lemma. The statement on $\overline{\cG}_{-a}$
is proved in the same way.
\end{proof}
The twistor $S^{can}_{(a,0)}(E)$ is constructed in \cite{Mo2} by
patching $\cG_a$ and $\overline{\cG}_{-a}$ via their identifications
with $\cH^\infty_{e^{2\pi i a}}$ on $\dC^*$. Therefore we obtain
the following result which shows the first part of theorem \ref{theoMainCompMochizuki}.
\begin{corollary}
Let $(H,H'_\dR,\nabla,P,w)$ be a TERP-structure inducing a nilpotent orbit
and $E=p_*\cC^{\infty h}(\widehat{\pi^* H})$ the corresponding harmonic bundle.
Then the $\dP^1$-bundles $\widehat{H^{el}}$ and $\oplus_{a\in L}S^{can}_{(a,0)}(E)$ are
isomorphic.
\end{corollary}
It remains to identify $(\widehat{W}_\bullet, \widehat{S}, \widehat{N})$ on $\cO(\widehat{H^{el}})$ with
$(\widehat{\cW}_\bullet, \widehat{\cS}, 2\pi \cdot\widehat{\cN})$ on
$\oplus_a S^{can}_{(a,0)}(E)$.
The filtrations $\widehat{W}_\bullet$ and $\widehat{\cW}_\bullet$ are
defined as usual by the nilpotent morphisms $\widehat{N}$ and $\widehat{\cN}$,
so that it suffices to show the identification of $\widehat{S}$ with $\widehat{\cS}$
and of $\widehat{N}$ with $2\pi \cdot\widehat{\cN}$.

In order to compare the constructions of the parings and the nilpotent
morphisms, remember that Mochizuki starts from a general VPTS
(called $(\cE^\Delta, \bD^\Delta, h)$ on $\Delta^*$ in \cite[Corollary 11.8]{Mo2})
which is in our case constructed from the variation of TERP-structures $(G,\nnn,G'_\dR, P, w)$.
In particular, the pairing $\widehat{S}$ on $\widehat{G}$ is precisely the one from lemma \ref{t3.14}.
The pairing $\widehat{\cS}$ on $\oplus_a S^{can}_{(a,0)}(E)$ is defined in \cite[11.3.5]{Mo2} through
a series on lemmas and intermediate computations. If one follows these
definitions carefully, it becomes clear that the isomorphism
$\oplus_a S^{can}_{(a,0)}(E)\cong \cO(\widehat{H^{el}})$ identifies
$\widehat{\cS}$ with $\widehat{S}$.
The map $\widehat{\cN}$ is defined on $\oplus_a \cG_a$
resp. $\oplus_a \overline{\cG}_{-a}$ by the morphism $\frac{1}{iz}\cN_E$ resp.
by $iz\cN_{\overline{E}}$ \cite[11.3.6]{Mo2} with $\cN_E$ from \ref{deftheoKMSS}. This implies that
$\widehat{\cN}:\oplus_a S^{can}_{(a,0)}(E)\rightarrow S^{can}_{(a,0)}(E)\otimes\cO_{\dP^1}(2)$.
Under the isomorphism
$\oplus_a \cG_a \cong \cO(H^{el})$ resp.
$\oplus_a \overline{\cG}_{-a} \cong \cO(\widehat{H}^{el}_{|\dP^1\backslash\{0\}})$,
$\frac{1}{z}\cN_E$ corresponds to $-\frac{1}{2\pi i}\widehat{N}_{|\dC}$
and $z\cN_{\overline{E}}$ to $\frac{1}{2\pi i}\widehat{N}_{|\dP^1\backslash\{0\}}$.
Therefore, $\widehat{\cN}$ on $\oplus_a S^{can}_{(a,0)}(E)$ corresponds precisely
to $\frac1{2\pi}\widehat{N}$ on $\cO(\widehat{H}^{el})$.

\vspace*{1cm}

\textbf{Remark:} In \cite{Mo2} there is another construction of a limit mixed twistor structure
which uses the same bundles $\cG_a$ and $\overline{\cG}_{-a}$ as above but a different
gluing. This gluing procedure depends this time on a point $x$ ``near
the origin in $\Delta^*$'' and produces a $\dP^1$-bundle denoted by $S_{(a,0)}(E,x)$.
Although we do not need this construction in order to prove theorem \ref{theoMainResultRegSing},
we will give here the corresponding statement for the case of nilpotent orbits of TERP-structures.
\begin{lemma}
Consider the ``rescaled'' TERP-structure
$\pi^*_x(H,H'_\dR,P,w)$, then the corresponding twistor
$\widehat{(\pi^*_x H)^{el}}$ is isomorphic to
$\bigoplus_{0<a \leq 1} S_{(a,0)}(E,x)$
and under this isomorphism $\widehat{W}_\bullet,\widehat{S}$ get
identified with $\widehat{\cW}_\bullet,\widehat{\cS}$ and
$\widehat{N}$ is mapped to $2\pi \cdot \widehat{\cN}$.
\end{lemma}
\begin{proof}
The gluing of $\cG_a$ and $\overline{\cG}_{-a}$ to
$S_{(a,0)}(E,x)$ is done via an identification
$$
(\cG_a)_{|\dC^*} \cong \cG_a(\cE)_{|\dC^*\times \{0\}}
\xymatrix@!0{\ar[r]^\cong_{\Phi_x}&}
\overline{\cG}_{-a}(\overline{\cE})_{|\dC^*\times \{0\}}
\cong (\overline{\cG}_{-a})_{|\dC^*}
$$
where $\cG_a(\cE)$ resp. $\cG_{-a}(\overline{\cE})$ are locally free sheaves over $\cO_{\dC^*\times\Delta}$
resp. $\cO_{\dC^*\times\overline{\Delta}}$, defined
in \cite[sections 10.1.4 and 10.3.2]{Mo2}. In our situation, the fact that the filtration $\cF_\bullet$
on $\cH^\infty$ is trivial shows that
$$
\cG_a(\cE) \cong \left\{s\in i_*\cO(G'_{e^{2\pi i a}})\,|\, |s|_{p^*h}\in O(|r|^{-a-\epsilon})\,\forall\epsilon>0 \right\}
$$
where $i:\dC^*\times\Delta^*\hookrightarrow\dC^*\times\Delta$ and $G'_\lambda$ denotes
the flat generalized eigensubbundle of $G'$ with respect to either horizontal or
vertical monodromy.
From the proof of lemma \ref{lemIdentG-Hel} it is obvious that $\cG_a(\cE)$ resp.
$\overline{\cG}_{-a}(\overline{\cE})$ are
generated by sections of the form $r^{a\mathit{Id}-N/2\pi i}es(A,a+k)$
resp. $\overline{r}^{a\mathit{Id}+N/2\pi i}es(\overline{A},-a+l)$ with $k,l\in \dZ$.
The identification $\Phi_x$ is defined in \cite[section 11.3.3]{Mo2} as
a composition $\Phi_x={\Phi^\dag}^{-1}_{x,0}\circ \Phi_{x,0}$ where
$$
\begin{array}{c}
\Phi_{x,0}: \cG_a(\cE)_{|\dC^*\times\{0\}}\cong \cG_a(\cE)_{|\dC^*\times\{x\}} \\ \\
\Phi^\dag_{x,0}: \overline{\cG}_{-a}(\overline{\cE})_{|\dC^*\times\{0\}}\cong \overline{\cG}_{-a}(\overline{\cE})_{|\dC^*\times\{x\}}
\end{array}
$$
are isomorphisms of fibres which in our situation boil down to restricting
the sections of the above type to $\dC\times\{x\}$.
The isomorphism $\cG_a(\cE)_{|\dC^*\times\{x\}}\cong
\overline{\cG}_{-a}(\overline{\cE})_{|\dC^*\times\{x\}}$ in \cite{Mo2} is
then simply the identification $\tau:(\pi_x^*H)^{el}\cong\gamma^*(\overline{\pi_x^*H})^{el}$
which shows that $S_{(a,0)}(E,x) \cong \widehat{(\pi_x^*H)^{el}}$.
The statements on $\widehat{N}, \widehat{S}, \widehat{\cN}$ and $\widehat{\cS}$
are proved as in theorem \ref{theoMainResultRegSing}.
\end{proof}

\section{Sabbah's mixed Hodge structures}\label{c7}
\setcounter{equation}{0}

The results of the previous chapter strongly rely on the
assumption of regularity of the given TERP-structure.
In the general case, one can no longer define a filtration
on the space $H^\infty$ as in corollary \ref{corFiltRegSing}.
However, there exists another procedure due to Sabbah (\cite{Sa2})
which applies to arbitrary bundles with meromorphic connections.
It was used to construct MHS for tame functions, see chapter \ref{c11}.
It uses global sections with moderate growth at infinity.
We recall it briefly below. We will get
a statement similar to theorem \ref{theoMainResultRegSing}
but this time PHMS will correspond to Sabbah orbits instead
of nilpotent orbits (which explains their name).

Consider an arbitrary TERP-structure $(H,\nnn,H'_\dR,P,w)$.
In the spirit of the beginning of chapter \ref{c6}, we will consider
the Deligne-extensions of $H$ but this time over infinity. More precisely,
denote $\widetilde{i}:\dC\hookrightarrow \dP^1$ and put
$$
\widetilde{i}_*\cO(H)\supset \cO(H_{<\infty}) :=
\left\{s\in \cO(H)\,|\,s\mbox{ has moderate growth at }\infty\right\}
$$
Then $H_{<\infty}$ is a bundle on $\dP^1$, meromorphic at infinity
(i.e., a locally free $\cO_{\dP^1}(*\infty)$-module). The rational
bundle corresponding to it by GAGA is denoted by $G_0$. It is a free $\dC[z]$-module.
It is a lattice at zero inside of $G:=G_0\otimes\dC[z,z^{-1}]$, which is rational
with poles at zero and infinity. Moreover, for any $\alpha\in \dC$ we can consider the following
lattices at infinity
$$
\begin{array}{c}
V_\alpha^{Sab}:=\left\{\omega\in G\,\left|\, \omega=\sum_{\beta\leq \alpha} s(\omega,\beta) \in C^\beta\right.\right\}, \\ \\
V_{<\alpha}^{Sab}:=\left\{\omega\in G\,\left|\, \omega=\sum_{\beta <\alpha} s(\omega,\beta) \in C^\beta\right.\right\}.
\end{array}
$$
Due to the formula $z^{-1}\nabla_{\partial_{z^{-1}}}=-z\nabla_{\partial_z}$,
it does not matter whether we consider elementary sections at zero
or at infinity. Therefore the above formulas defines the Deligne lattices at infinity
of $H_{<\infty}$. $V^{Sab}_\alpha$ and $V^{Sab}_{<\alpha}$
are free $\dC[z^{-1}]$-modules inside of $G$ of maximal rank
corresponding to algebraic bundles on $\dP^1\backslash\{0\}$.
The advantage of working algebraically is that we can use
the V-filtration at infinity to define a TERP-structure generated by
elementary sections.
Namely, we write $Gr^{V^{Sab}}_\alpha:=V^{Sab}_\alpha/V^{Sab}_{<\alpha} \cong C^\alpha$
as before and use the filtration induced by the $V^{Sab}_\alpha$'s on $G_0$.
\begin{definition}
Let $(H,\nnn,H'_\dR,P,w)$ be any TERP-structure. Define
$$
G_0^{Sab}=\oplus_\alpha \dC[z]Gr_a^{V^{Sab}}G_0
$$
as the algebraic bundle
generated by the Sabbah-principal parts
of sections if $G_0$. Let $H^{Sab}$ be its analytic counterpart,
i.e., $\cO(H^{Sab})=\oplus \cO_\dC Gr_a^{V^{Sab}}G_0$.
The spectrum in this situation is defined as
$$
\mathit{Sp}^{Sab}(H,\nabla) = \sum_{\alpha \in\dC}\nu(\alpha)\alpha\in\dZ[\dC]
\;\;\;
\mbox{with}
\;\;\;
\nu(\alpha):=\dim_\dC\left(\frac{Gr_\alpha^{V^{Sab}} G_0}{Gr_\alpha^{V^{Sab}} zG_0}\right)
$$
Write as before the spectral numbers as an ordered $\rank(H)$-tuple of (possibly repeated)
numbers $\alpha_1\leq...\leq \alpha_{\rank(H)}$.
We call them spectral numbers at infinity.
\end{definition}
The following result is the analogue of lemma \ref{lemHelTERP}.
\begin{lemma}
Let $(H,H'_\dR, \nabla,P,w)$ be a TERP-structure.
\begin{enumerate}
\item
The spectral numbers derived from the V-filtration at infinity satisfy the
symmetry $\alpha_i+\alpha_{\rank(H)+1-i}=w$.
\item
$(H^{Sab},H'_\dR,\nabla,P,w)$ is again a TERP-structure. It is generated by
elementary sections and thus regular singular.
\end{enumerate}
\end{lemma}
\begin{proof}
$(\cO(H^{Sab}),\nabla)$ has a pole of order at most two at zero by
the same proof as in lemma \ref{lemHelTERP}.
Moreover, it is obvious that the pairing $P$ induces a pairing
$G\otimes j^*G\rightarrow \dC[z,z^{-1}]$ sending $G_0$
to $z^w\dC[z]$. Then we are exactly in the situation considered
in \cite{Sa2}. It is shown in section 3 of loc.cit. that
the spectral numbers are symmetric and the properties of $P$ as a pairing
on $G_0^{Sab}$ we are after are a direct consequence.
\end{proof}

As in the regular singular case, the last lemma allows us to define
a filtration on the space $H^\infty$: The TERP-structure $(H^{Sab},H'_\dR,\nabla,P,w)$
is generated by elementary sections and is therefore equivalent to the data considered in
lemma \ref{lemCorElTERP}. In particular, the formula
\begin{eqnarray}\label{7.1}
F^p_{Sab}H^\infty_{e^{-2\pi i \alpha}}
:=\psi^{-1}_\alpha z^{p+1-w}Gr_{\alpha+w-1-p}^{V^{Sab}}G_0
\cong \psi^{-1}_\alpha\left(z^{p+1-w}(C^{\alpha+w-1-p}
\cap\cO(H^{Sab})_0)\right)
\end{eqnarray}
for $\alpha\in(0,1]+i\dR$ defines a decreasing filtration on $H^\infty$ such that
its twisted version $\widetilde{F}_{Sab}^\bullet H^\infty := G^{-1} (F_{Sab}^\bullet H^\infty)$
satisfies the orthogonality conditions \eqref{eqOrthogonal} and therefore gives
an element of the classifying space $\check{D}$.

The next result is the ``Sabbah-orbit'' version of theorem \ref{theoMainResultRegSing}.
\begin{theorem}\label{theoMainResultSabbah}
For an arbitrary TERP-structure $(H,H'_\dR, \nabla,P,w)$, the following conditions are equivalent.
\begin{enumerate}
\item
$(H,H'_\dR, \nabla,P,w)$ induces a Sabbah orbit.
\item
$(H^\infty, H_\dR^\infty, N, S, \widetilde{F}_{Sab}^\bullet)$ defines
a PMHS of weight $w-1$ resp. $w$ on $H^\infty_{\arg\neq 0}$ resp. $H^\infty_{\arg = 0}$.
\end{enumerate}
\end{theorem}
\begin{proof}
The proof of 2) ---> 1) is similar to the proof of \cite[theorem 7.20]{He2}.
The parts (I)--(III) in that proof are unchanged, the parts (IV) and (V)
can be adapted easily.

The proof of 1)--->2) is virtually the same as the one in chapter \ref{c6}. We
give only comments on the necessary adjustments. The first point
is that by lemma \ref{lemTameness}, if we consider the
variation of pure polarized twistor structures associated to the
orbit $K:=(\pi')^* H$, then taking global sections along the
projection $p:\dP^1\times \Delta^*\rightarrow \Delta^*$ gives a
harmonic bundle $E$ on $\Delta^*$ with tame behavior at $0\in\Delta$.

The discussion of $\cU_{|r}$ in the proof of lemma \ref{lemTameness} shows that
on any extension $_a\cE^0$ of the harmonic bundle $E$ the residue at 0
of the Higgs field is nilpotent. Therefore the KMSS-spectrum
of the variation of twistor structures is precisely as described
in lemma \ref{lemPropKMSSNilpOrb}. All the properties of
the extension sheaves $_a\cE$ remain true.
The vertical and horizontal monodromy are not equal, but
inverse to each other. We continue to denote by $\cH^\infty$
the flat bundle $\psi_r((K')^\nabla)$, and by $\cH^\infty_\lambda$
the eigenbundle of the \emph{horizontal} monodromy. With this
convention, the definition of the extensions $\cV^\alpha$ remains unchanged.
However, we see that now the bundle
$\cO(K)$ is generated by sections of the form
$$
\sum_{i\geq 1} r^{-\alpha_{ij}\mathit{Id}+\frac{N}{2\pi i}} es(A_i, \alpha_{ij})
$$
where $\sum_{i\geq 1} es(A_{ij}, \alpha_{ij})$ for $j=1,\ldots,\rank(H)$
is a set of generating sections of $\cO(H^{Sab})$ and where this time
$\alpha_{1j}>\alpha_{2j}>\ldots$.
The remaining part of the proof is exactly the same. We identify the
limit polarized mixed twistor structure $(\oplus_a S_{(a,0)}^{can}(E), \widehat{\cW}_\bullet,
\widehat{\cN}, \widehat{\cS})$ of Mochizuki with the
twistor $(\widehat{H^{Sab}},\widehat{W}_\bullet, -\frac1{2\pi}
\widehat{N}, \widehat{S})$ which is generated by elementary sections
and corresponds by lemma \ref{lemCorPMHS-PMTS}
to a sum of two PMHS.
Note that the different signs $+N$ instead of $-N$ in the PMHS in 2. and
$-\frac{1}{2\pi}$ instead of $\frac{1}{2\pi}$ in the PMTS result from the fact
that now horizontal and vertical monodromy are inverse, not equal.
\end{proof}

\section{Formal structure and Stokes structure}\label{c8}
\setcounter{equation}{0}

From now on we will consider TERP-structures which are not necessarily
regular singular. This chapter is devoted to discuss
both formal and Stokes structures of these objects. Our
main references for the facts discussed below are
\cite{Mal1} and \cite{Sa4}.

Let $(H,\nnn,H'_\dR,P)$ be a TERP-structure.
We will forget about the real subbundle
$H'_\dR$ for a moment. The remaining object $(H,\nnn,P)$ is called a
TEP-structure.
We might as well consider only the germ $(\cH_0,\nabla,P)$ at zero. It turns
out that for TEP-structures, it is sufficient to work at the
level of formal power series,
i.e., to consider only the formalization
$(\cH_0,\nabla,P)\otimes \dC[[z]][z^{-1}]$. However,
to incorporate the real subbundle of a TERP-structure, we will need
to use the Stokes structure associated to the irregular connection $\nabla$.

To start with, let us simplify the situation even more and restrict our attention
to the germ $\cH_0[z^{-1}]=\cH_0\otimes\dC\{z\}[z^{-1}]$. Remember that
the connection has a pole of order two on $\cH_0[z^{-1}]$.
By the theorem of Turrittin (e.g. \cite[2.1]{Mal1})
there is a finite ramification $r_n:\dC\to\dC,\ z\mapsto z^n,$ such that
$r_n^*(\HH_0[z^{-1}],\nnn)$ is formally isomorphic to a sum
$\bigoplus_{i=1}^l e^{f_i}\otimes (\RR_i,\nnn_i)$, where $f_i\in\dC[z^{-1}]$,
$\RR_i$ is a $\Czm$-vector space and $(\RR_i,\nnn_i)$ is regular singular.
The symbol $e^{f_i}$ denotes the bundle of rank one equipped with the
connection given by the one form $df_i$ (equivalently, $e^{f_i}\otimes \RR_i$ is the
germ of the meromorphic extension at zero given by multiplying
sections in $\RR_i$ by $e^{f_i}$).

In the sequel, we will always make the simplifying assumption that the
ramification is unnecessary.
This is satisfied in all examples (singularity theory) and
potential examples (quantum cohomology) which we have in mind.
In this situation, the exponents that occur in the
formal decomposition are simply $f_i=-\frac{u_i}{z}$,
where $u_i$ is an eigenvalue of the pole part $\UU$.
This is encapsulated in the following definition.

\begin{definition}\label{def-ramification}
A TEP-structure (or its formalization at zero) $(H,\nnn,P)$
is said to require no ramification
iff the germ $(\HH_0[z^{-1}],\nnn)$ is formally isomorphic to a sum
$\bigoplus_{i=1}^l e^{-u_i/z}\otimes (\RR_i,\nnn_i)$  where
$\RR_i$ is a $\Czm$-vector space, $(\RR_i,\nnn_i)$ is regular
singular, and $u_1,...,u_l$ are the different eigenvalues of the
pole part $\UU=[z\nnn_\zdz]$. A TERP-structure is said to require
no ramification iff this is true for the corresponding TEP-structure.
\end{definition}

Given a TEP-structure $(\HH_0,\nnn,P)$, then the $\Cz$-lattice
$e^{-u/z}\cdot \HH_0$ with the induced connection (in the above sense) and the
induced pairing is also a TEP structure.

\begin{lemma}\label{lemma-TEP}
Let $(\HH_0,\nnn,P)$ be a TEP-structure which does not require
a ramification and $\Psi$
a formal isomorphism from $(\HH_0[z^{-1}],\nnn)$ to a sum
$\bigoplus_{i=1}^l e^{-u_i/z}\otimes (\RR_i,\nnn_i)$ as above.
Then there are unique $\Cz$-lattices $(\HH_i)_0\subset \RR_i$
and pairings $P_i$
such that $((\HH_i)_0,\nnn_i,P_i)$ is a  regular singular TEP-structure
and $\Psi$ is an isomorphism of the formal TEP-structures
$(\HH_0\otimes \dC[[z]],\nnn,P)$ and
$\bigoplus_{i=1}^l e^{-u_i/z}\otimes
((\HH_i)_0\otimes \dC[[z]],\nnn_i,P_i)$.
The TEP-structures $((\HH_i)_0,\nnn_i,P_i)$ are called the regular singular
pieces of $(\HH_0,\nnn,P)$.
\end{lemma}

\begin{proof}
Exercise 5.9 in \cite[II]{Sa4} shows that the $\dC[[z]]$-lattice
$\Psi(\HH_0\otimes \dC[[z]])\subset \Psi(\HH_0\otimes \dC[[z]][z^{-1}])$
is compatible with the splitting
$\Psi(\HH_0\otimes \dC[[z]][z^{-1}])=\bigoplus_{i=1}^l
e^{-u_i/z}\otimes (\RR_i\otimes \dC[[z]][z^{-1}])$,
that means, it splits into $l$ $\dC[[z]]$-lattices.
Because of \cite[III.2.1]{Sa4} these $l$ $\dC[[z]]$-lattices arise as
$e^{-u_i/z}\otimes ((\HH_i)_0\otimes \dC[[z]])$ from some
$\Cz$-lattices $(\HH_i)_0\subset\RR_i$.

The connection on $\RR_i$ has a pole of order $\leq 2$ with respect to
$(\HH_i)_0$, because the same holds for the connection on the factor
$e^{-u_i/z}$ and for the connection on
$e^{-u_i/z}\otimes (\RR_i\otimes \dC[[z]][z^{-1}])$
with respect to the $\dC[[z]$-lattice
$e^{-u_i/z}\otimes ((\HH_i)_0\otimes \dC[[z]])$.

The pairing $\Psi(P)$ on $\Psi(\HH_0\otimes \dC[[z]])$ induced by $P$
gives the first of the following isomorphisms
(compare the proof of lemma \ref{lemHelTERP}).
\begin{eqnarray}\nonumber
(z^{-w}\Psi(\HH_0\otimes \dC[[z]]),\nnn)
&\cong&
j^*(\Psi(\HH_0\otimes\dC[[z]])^*,\nnn^*)\\
&\cong&
j^*\left( \bigoplus_{i=1}^l e^{u_i/z}\otimes((\HH_i)_0^*,\nnn_i^*)\otimes
\dC[[z]]\right) \nonumber\\
&\cong&
\bigoplus_{i=1}^l e^{-u_i/z}\otimes j^*((\HH_i)_0,\nnn_i^*)\otimes \dC[[z]].
\label{8.1}
\end{eqnarray}
As the formal decomposition is unique (e.g. \cite[II.5.5]{Sa4}),
these isomorphisms respect the decomposition and give rise to isomorphisms
$(z^{-w}(\HH_i)_0,\nnn_i) \cong j^*((\HH_i)_0^*,\nnn_i^*)$ corresponding to pairings $P_i$ on $(\HH_i)_0$ which
make $((\HH_i)_0,\nnn_i,P_i)$ into regular singular TEP-structures.
Moreover, for $i\neq j$, the pairing $\Psi(P)$ vanishes on
$e^{-u_i/z}\otimes (\HH_i)_0\times e^{-u_j/z}\otimes (\HH_j)_0$.
\end{proof}

As remarked above, the formal structure of the connection $(\cH_0,\nabla)$
is not sufficient to treat the real subbundle $H'_\dR$. The necessary analytic information
is provided by what is called the Stokes structure associated to $(\cH_,\nabla)$.
Roughly speaking, it consists of data which keep track of the difference
of analytic liftings of the formal decompositions in sectors of $\dC$.
To discuss this Stokes structure, we need some notations.

For any $I\subset S^1$, let $\widehat I:=\{z\in \dC^*\ |\ \frac{z}{|z|}\in I\}$.
If $I\subset S^1$ happens to be open and connected, then $\widehat I$ is a sector.
Denote by $\cA[z^{-1}]$ the sheaf on $S^1$ of holomorphic functions
in neighborhoods of 0 in sectors which have an asymptotic development
in the sense of $\cite[3.]{Mal1}$.
For any $\xi\in S^1$ the Taylor development of functions in  $\cA[z^{-1}]_\xi$
yields a map $T:\cA[z^{-1}]_\xi\to \dC[[z]][z^{-1}]$
which is surjective by the lemma of Borel-Ritt.
We will need the subsheaves $\cA:=T^{-1}(\dC[[z]])$ and $\cA^{<0}:=\ker T$
of $\cA[z^{-1}]$. Remark that $\gamma^*(\oooo{\cA[z^{-1}]})$ is the sheaf on $S^1$ of
holomorphic functions in neighborhoods of $\infty$ in sectors
such that the functions have asymptotic developments at $\infty$.

Fix a TERP-structure $(H,\nnn,H'_\dR,P)$ which does not require a ramification and
a formal isomorphism
$$
\Psi: (\HH_0,\nnn,P)\otimes \dC[[z]] \
\stackrel{\cong}{\rightarrow} \
\bigoplus_{i=1}^l e^{-u_i/z}\otimes ((\HH_i)_0,\nnn_i,P_i)\otimes \dC[[z]]
$$
of TEP-structures. In this situation, we say that
$\xi\in S^1$ is a Stokes direction if there exist $i\neq j$ with
$\Re(\frac{u_i-u_j}{\xi})=0$. These Stokes directions form a finite set
$\Sigma\subset S^1$. For $\xi\in S^1 \backslash \Sigma$, the component of
$S^1-\Sigma$ which contains $\xi$ is denoted by $I(\xi)$ and we define
an order on $\{1,...,l\}$ depending on $\xi$ as follows
\begin{eqnarray}\label{8.2}
i\leq_\xi j &\stackrel{\textup{def}}{\iff}&  \Re\left(\frac{u_i}{\xi}\right)
< \Re\left(\frac{u_j}{\xi}\right) \quad \textup{ or }\quad i=j\\
&\iff& (z\mapsto e^{\frac{u_i-u_j}{z}})\in \cA_\xi
\iff (z\mapsto e^{(\oooo{u_i}-\oooo{u_j})z})\in
\gamma^*(\oooo{\cA_\xi}).\nonumber
\end{eqnarray}

In the sequel, we will work with a covering of $S^1$
by two open sets $I_\pm(\xi)$
defined in the following way: Choose once and for all a value
$\xi\in S^1\backslash \Sigma$ and put
$$
I_\pm(\xi):=I(\xi)\cup I(-\xi)\cup\{z\in S^1\ |\
\pm \Im(z /\xi)\leq 0\}.
$$
Each of the sets $I_+(\xi)$ and $I_-(\xi)$ contains exactly one of
two Stokes directions $\pm \xi'$ for any $\xi'\in \Sigma$.
Let $\cL$ be the local system of $(H,\nabla)_{|S^1}$, $\cL_\dR\subset\cL$
the local system of $(H'_\dR,\nabla)_{|S^1}$ and denote by $\cL_i$ the
local system of flat sections of the $i$-th regular singular piece $H_i$.
The theorem of Hukuhara (e.g. \cite[3.]{Mal1} or \cite[II.5.12]{Sa4}) and
the discussion in \cite[4.+5.]{Mal1} yield the following result.

\begin{lemma}\label{lemma-Stokes1}
Let $(H,\nnn,H'_\dR,P)$ be a TERP-structure requiring no ramification
and let $\Psi$ and $\bigoplus_{i=1}^l e^{-u_i/z}\otimes
((\HH_i)_0,\nnn_i,P_i)$ and $\xi\in S^1$
be as above.
There is a unique lift $\Psi_\pm$ of $\Psi$ to an isomorphism of sheaves
on $I_\pm(\xi)$,
\begin{eqnarray}\label{8.3}
\Psi_\pm:\cA[z^{-1}]_{|I_\pm(\xi)}\otimes (\cH_0[z^{-1}],\nnn)\longrightarrow
\cA[z^{-1}]_{|I_\pm(\xi)}\otimes
\left(\bigoplus_{i=1}^l e^{-u_i/z}\otimes ((\HH_i)_0[z^{-1}],\nnn_i)
\right).
\end{eqnarray}
The underlying isomorphism of local systems on $I_\pm( \xi)$ is denoted
by the same letter
$$
\Psi_\pm: \cL_{|I_\pm(\xi)}\to \bigoplus_{i=1}^l \cL_{i|I_\pm(\xi)}.
$$
Let $\cL_{\pm,i}:=\Psi_\pm^{-1}(\cL_{i|I_\pm(\xi)})$, so that
$\cL_{|I_\pm(\xi)} = \bigoplus_{i=1}^l \cL_{\pm,i}.$
At $\pm \xi$ the two splittings induce (restrictions of) projections
$t_{ij}^{(\pm)}: (\cL_{+,i})_{\pm\xi} \to (\cL_{-,j})_{\pm \xi}$ satisfying
$$
t_{ij}^{(+)}=0\;\;\;\Longleftarrow \;\;\;i<_\xi j \;\;\;\iff
\;\;\;j<_{-\xi}i \;\;\;\Longrightarrow \;\;\;t_{ji}^{(-)}=0.
$$
Therefore $t_{ii}^{(\pm)}$ are isomorphisms, and the local system obtained
by gluing $\cL_{+,i}$ and $\cL_{-,i}$ with $t_{ii}^{(\pm)}$
is mapped by $\Psi_\pm$ to $\cL_i$.
\end{lemma}

The following lemma extends the decomposition of the meromorphic bundle
to the lattice $\cH_0$ and describes the behavior of the pairing.
\begin{lemma}\label{lemma-Stokes2}
We consider the same situation as in the last lemma. The formal isomorphism
$\Psi$ respects the $\dC[[z]]$-lattice $\cH_0\otimes\dC[[z]]$
and thus induces isomorphisms
\begin{eqnarray}\label{8.4}
\Psi_\pm:\cA_{|I_\pm(\xi)}\otimes (\HH_0,\nnn) \to
\cA_{|I_\pm(\xi)}\otimes
\left(\bigoplus_{i=1}^l e^{-u_i/z}\otimes ((\HH_i)_0,\nnn_i)\right).
\end{eqnarray}

The two splittings $\bigoplus_{i=1}^l \cL_{\pm,i}$ of $\cL_{|I_\pm(\xi)}$
are dual with respect to the pairing $P$.
The maps $t_{ij}^{(+)}$ and $t_{ji}^{(-)}$ determine each other by
\begin{eqnarray}\label{8.5}
P(t_{ij}^{(+)}-,-)=P(-,-)=P(-,t_{ji}^{(-)}-):
(\cL_{+,i})_\xi \times (\cL_{+,j})_{-\xi}\longrightarrow \dC.
\end{eqnarray}
$P$ restricts to the local system obtained by gluing $\cL_{\pm,i}$ with
$t_{ii}^{(\pm)}$. It is identified with $P_i$ by the isomorphism
with $\cL_i$.
\end{lemma}

\begin{proof}
The first statement is clear. For the second point,
suppose $i\neq j$, $\sigma_i\in (\HH_i)_0$ and
$\sigma_j\in (\HH_j)_0$. The last statement in the proof of lemma
\ref{lemma-TEP} shows that
$P\left(\Psi_+^{-1}(e^{-u_i/z}\sigma_i),\Psi_-^{-1}(e^{-u_j/z}\sigma_j)\right)$
is an element in $\Gamma(I_+(\xi),\cA^{<0})$.
On the other hand, it takes the form
$e^{(u_j-u_i)/z}\cdot \sum_{\alpha\geq \alpha_0,k}
a_{\alpha,k}z^\alpha(\log z)^k$. Therefore it vanishes.
This shows the duality property of the two splittings
and also yields both equalities in \eqref{8.5}.
The remaining statements are obvious.
\end{proof}
The next definition introduces the condition we need in order to
equip the regular singular pieces with a real structure, i.e., make
them into TERP-structures.

\begin{definition}\label{def-compatibility}
Let $(H,\nnn,H'_\dR,P)$ be a TERP-structure requiring no ramification.
Consider the data from lemma \ref{lemma-Stokes1}. We say that
real structure and Stokes structure are compatible
if $\cL_{\pm,i}=\oooo{\cL_{\pm,i}}$.
\end{definition}

\begin{lemma}\label{lemma-Stokes3}
Let $(H,\nnn,H'_\dR,P)$ be a TERP-structure requiring no ramification
and which has compatible real structure and Stokes structure.
Then the regular singular pieces $(H_i,\nnn_i,P_i)$,
which are TEP-structures by lemma \ref{lemma-TEP}, are naturally equipped
with real structures and become TERP-structures.
\end{lemma}

\begin{proof}
The identity map on the stalks $\cL_\pm$ can be decomposed as
$\id =\sum_{j\leq_{\pm\xi} i} t_{ij}^{(\pm)}:\cL_{\pm\xi}\to \cL_{\pm\xi}$.
As the identity obviously respects the real structure, the same holds
true for the individual maps
$t_{ij}^{(\pm)}:(\cL_{+,i})_{\pm\xi}\to (\cL_{-,j})_{\pm\xi}$.
By lemma \ref{lemma-Stokes1}, $\cL_i$ is canonically isomorphic to the
local system obtained by gluing $\cL_{\pm,i}$ with $t_{ii}^{(\pm)}$.
Therefore it carries a canonical real structure.
$P_i$ maps the real local system in $\cL_i$ to $i^w\dR$
because of the identification in lemma \ref{lemma-Stokes2}
with the restriction of $P$.
This gives a real flat subbundle $H'_{i,\dR}\subset H'_i$ such that
$(H_i,\nnn_i,H'_{i,\dR},P_i)$ is a TERP-structure.
\end{proof}

The following lemma is now an immediate consequence of the
preceding results.
\begin{lemma}\label{lemma-Stokes4}
Let $(H,\nnn,H'_\dR,P)$ be a TERP-structure requiring no ramification.
The integrable twistor $(\widehat{\cH},\nabla)$ has a
pole of order two at infinity and
the map $\tau \circ \Psi_\pm\circ\tau$ is an isomorphism
of sheaves on $I_\pm(\xi)$:
\begin{eqnarray}\label{8.6}
\tau\circ \Psi_\pm\circ \tau :
\gamma^{*}\left(\oooo{\cA}_{|I_\pm(\xi)}\right)\otimes
({\widehat\HH}_\infty,\nnn) \longrightarrow
\gamma^{*}\left(\oooo{\cA}_{|I_\pm(\xi)}\right)\otimes
\left(\bigoplus_{i=1}^l e^{-\oooo{u_i}\cdot z}\otimes
((\widehat{\HH_i})_\infty,\nnn_i)\right).
\end{eqnarray}
It is the unique lift on $I_\pm(\xi)$ of the corresponding formal
isomorphism.
On the level of local systems it is given by the composition
$\oooo{-}\circ \Psi_\pm \circ\oooo{-} =:\oooo{\Psi_\pm}$.
Therefore it induces the splitting
$\cL_{|I_\pm(\xi)} = \bigoplus_{i=1}^l \oooo{\cL_{\pm,i}}$.
The morphisms $\Psi_\pm$ and $\tau\circ \Psi_\pm\circ \tau$ glue to an isomorphism
of sheaves on the real blow up $[0,\infty]\times I_\pm(\xi)$
at $0$ and $\infty$ of the sector $\widehat I_\pm(\xi)$
iff real structure and Stokes structure are compatible.
\end{lemma}

\section{Mixed TERP-structures}
\setcounter{equation}{0}\label{c9}

The correspondence between nilpotent orbits of Hodge structures
and PHMS due to Cattani, Kaplan and Schmid (theorem \ref{theoSchmidCorres})
was generalized to regular singular TERP-structures in chapter \ref{c6}
(theorem \ref{theoMainResultRegSing}).
In conjecture \ref{main-conjecture} we propose a further generalization
to arbitrary TERP-structures.
We will give a complete proof of one direction of this correspondence.
In order to state the general correspondence, we need
replace the notion of PHMS by what we call mixed TERP-structure. The definition
is straightforward after what has been said in the last chapter.
\begin{definition}\label{def-mixedTERP}
A TERP-structure is a mixed TERP-structure if it
does not require a ramification,
if real structure and Stokes structure are compatible and if the
regular singular pieces of lemma \ref{lemma-TEP}, which are TERP-structures
by lemma \ref{lemma-Stokes3}, induces PMHS as in theorem
\ref{theoMainResultRegSing}.
\end{definition}
With this definition at hand, the main conjecture is very simple to state,
and takes precisely the same form as in the regular singular case of
chapter \ref{c6}.
\begin{conjecture}\label{main-conjecture}
A TERP-structure which does not require a ramification
is a mixed TERP-structure if and only if it induces a nilpotent orbit.
\end{conjecture}

The following is the main result of this paper.

\begin{theorem}\label{main-theorem}
\begin{enumerate}
\item
The conjecture is true if the TERP-structure is regular singular.
\item
The implication $\Rightarrow$ is true.
\end{enumerate}
\end{theorem}
As already said, the first statement is precisely theorem
\ref{theoMainResultRegSing}.
The remaining part of this chapter is concerned with a proof
of the second claim.

Part 2. is already known in two cases which are opposite to each other
in a certain sense: If the TERP-structure is regular singular,
i.e., if $\cU$ is nilpotent, it is precisely \cite[theorem 7.20]{He2}.
If all eigenvalues of $\cU$ are different
(we call the TERP-structure semi-simple in that case),
then it follows from \cite[proposition 2.2]{Du}.
We will use \cite[theorem 7.20]{He2}
and combine it with ideas from \cite{Du} to get the general case.
This will contain a new proof of the semi-simple case.

Let $(H,\nnn,H'_\dR,P)$ be a mixed TERP-structure of weight $w$.
We use all notations and objects from chapter \ref{c8}.
In particular, $\xi\in S^1 \backslash \Sigma$ is chosen,
giving rise to the covering
$S^1=I_+(\xi)\cup I_-(\xi)$ with
$I_+(\xi)\cap I_-(\xi)=I(\xi)\cup (-I(\xi))$.
$\cL$ is the local system of $(H,\nabla)_{|S^1}$.
There are the canonical splittings
$\cL_{|I_\pm(\xi)} = \bigoplus_{i=1}^l \cL_{\pm,i}$ which are induced
from the unique isomorphisms
$\Psi_\pm:\cL_{|I_\pm(\xi)}\to \bigoplus_{i=1}^l \cL_{i|I_\pm(\xi)}$
in \eqref{8.4}.
We choose the numbering $1,\ldots,l$ of the block decomposition
in such a way
that for all $i,j\in\{1,\ldots,l\}$ we have $i\leq j\iff i\leq_\xi j$.
By assumption, the regular singular pieces
$(H_i,\nnn_i,H'_{i,\dR},P_i)$ are TERP-structures which induce
PMHS as in theorem \ref{theoMainResultRegSing}.
Also lemma \ref{lemma-Stokes4} applies.

Denote $n:=\rank(H)$ and $n_i:=\rank(H_i)$.
Let $\uuuu{e}^p=(e_1^\pm,...,e_n^\pm)$ be bases of $\cL_{|I_\pm(\xi)}$
satisfying $P(e_i^+,e_j^-)=(-i)^w\delta_{ij}$
and which are adapted to the splittings, i.e.,
$(e_1^\pm,...,e_{n_1}^\pm)$ is a basis of $\cL_{\pm,1}$ etc.
Such bases exist by lemma \ref{lemma-Stokes2}.
The images $\uuuu{f}^\pm=(f_1^\pm,...,f_n^\pm)
=(\Psi_\pm(e_1^\pm),...,\Psi_\pm(e_n^\pm))$ are bases of the local systems
$\bigoplus_{i=1}^l \cL_{i|I_\pm(\xi)}$, adapted to the splittings.
Because of \eqref{8.5} there is a unique matrix $T\in GL(n,\dC)$ with
\begin{eqnarray}\label{9.1}
(\uuuu{e}^-)_\xi = (\uuuu{e}^+)_\xi\cdot T\quad \textup{ and }\quad
(\uuuu{e}^-)_{-\xi} = (\uuuu{e}^+)_{-\xi}\cdot (-1)^wT^{tr}.
\end{eqnarray}
Write $T=(T_{jk})_{j,k=1,...,l}$ with blocks
$T_{jk}\in M(n_j\times n_k,\dC)$.
Lemma \ref{lemma-Stokes1} and the chosen ordering of the factors
of the decomposition implies that the matrix $T$ is block upper triangular.
Denote $T^{model}:=\textup{diag}(T_{11},...,T_{ll})$ and
$T^{Stokes}:=(T^{model})^{-1}\cdot T$. Then
\begin{eqnarray}\label{9.2}
(\uuuu{f}^-)_\xi = (\uuuu{f}^+)_\xi \cdot T^{model}
\quad \textup{ and }\quad
(\uuuu{f}^-)_{-\xi} = (\uuuu{f}^+)_{-\xi} \cdot (-1)^w (T^{model})^{tr}.
\end{eqnarray}

We want to prove that $\pi_r^*(H,\nnn,H'_\dR,P)$ is a polarized
pure TERP-structure for $|r|\ll 1$. As has been shown in lemma \ref{t4.4},
the corresponding twistors are isomorphic for $r$ having the same absolute value,
so it is sufficient to prove the case where $r\in\dR_{>0}$.

The family $\bigcup_{r>0}\pi_r^*(H,\nnn,H'_\dR,P)$ is isomonodromic.
This implies that data $\cL$, $\cL_{\pm,i}$, $\cL_i$, $\uuuu{e}^\pm$, $\uuuu{f}^\pm$
defined above for $r=1$ can be identified with the analogous data for any
$r>0$.
The family has constant Stokes structure (e.g. \cite[II.6.c]{Sa4}).
Therefore the isomorphisms
$\Psi_\pm:\cL_{|I_\pm(\xi)}\to \bigoplus_{i=1}^l \cL_{i|I_\pm(\xi)}$
of local systems extend for any $r>0$ to isomorphisms
\begin{eqnarray}\label{9.3}
\Psi_\pm(r):\cA_{|I_\pm(\xi)}\otimes \pi_r^*(\cH_0,\nnn)\to
\cA_{|I_\pm(\xi)}\otimes \left(\bigoplus_{i=1}^l e^{-\frac{u_i}{z\cdot r}}
\otimes \pi_r^*((\cH_i)_0,\nnn_i)\right).
\end{eqnarray}
As explained in the last chapter, we will use the notation $e^f \otimes$
to denote extensions of sheaves over zero or infinity (or both) twisted by multiplying
sections by the function $(z\mapsto e^f)$ on $\dC^*$.
Obviously $e^{-u_i/(z\cdot r)}\otimes \pi_r^*(H_i,\nnn_i,H'_{i,\dR},P_i)$
is a TERP-structure and is equal to
$\pi_r^*(e^{-u_i/z}\otimes (H_i,\nnn_i,H'_{i,\dR},P_i))$.
The extension to $\infty$ by $\tau$ satisfies
\begin{eqnarray}\label{9.4}
[\pi_r^*(e^{-\frac{u_i}{z}}\otimes \cH_i)]\widehat{\hspace*{0.3cm}} =
[e^{-\frac{u_i}{z\cdot r}}\otimes \pi_r^*\cH_i]\widehat{\hspace*{0.3cm}} =
e^{-\frac{u_i}{z\cdot r}-\frac{\oooo {u_i}\cdot z}{r}}\otimes
\widehat{\pi_r^*\cH_i}.
\end{eqnarray}
Lemma \ref{lemma-Stokes4} applies and shows that $\Psi_\pm(r)$ and
$\tau\circ\Psi_\pm\circ\tau$ glue to an isomorphism
\begin{eqnarray}\label{9.5}
\Psi_\pm(r):\widehat\cA_\pm \otimes \widehat{\pi_r^*\cH} \mapsto
\widehat\cA_\pm\otimes \left(\bigoplus_{i=1}^l
e^{-\frac{u_i}{z\cdot r}-\frac{\oooo {u_i}\cdot z}{r}} \otimes
\widehat{\pi_r^*\cH_i}\right).
\end{eqnarray}
Here $\widehat\cA_\pm$ is a sheaf on the real blow up
$[0,\infty]\times I_\pm(\xi)$ of $\widehat I_\pm(\xi)$ at $0$ and $\infty$.
It is the extension of $\cO_{\widehat I_\pm(\xi)}$
by $\cA_{|I_\pm(\xi)}$ at $0$ and $\gamma^*(\oooo{\cA_{|I_\pm(\xi)}})$
at $\infty$.

By \cite[theorem 7.20]{He2} resp. theorem \ref{theoMainResultRegSing}
2.$\Rightarrow$1., $\pi_r^*(H_i,\nnn_i,H'_{i,\dR},P_i)$ is a
polarized pure TERP-structure for small $r>0$.
One can choose a basis $\uuuu{\sigma}(r)=(\sigma_1(r),...,\sigma_n(r))$
of $\bigoplus_{i=1}^l \Gamma(\dP^1,\widehat{\pi_r^*\cH_i})$
compatible with the splitting.
The matrix
$P_{mat}(r):=z^{-w}\cdot P((\uuuu\sigma)^{tr},\tau(\uuuu{\sigma}))$
is independent of $z$, hermitian, positive definite and block diagonal.
The entries of the matrices $C^\pm$ defined by
\begin{eqnarray}\label{9.6}
\uuuu{\sigma}_{|\widehat I_\pm(\xi)} = \uuuu{f}^\pm\cdot C^\pm
\end{eqnarray}
take the form
$\sum_{\alpha_1\leq \alpha\leq \alpha_n} \sum_{k\in \dN\cup\{0\}}
a_{\alpha,k}(r)\cdot z^\alpha\cdot (\log z)^k$,
where $\alpha_1$ and $\alpha_n$ are the minimal and the maximal spectral number
of the regular singular TERP-structure
$\bigoplus_{i=1}^l \pi_r^*(H_i,\nnn_i,H'_{i,\dR},P_i)$.
Note that theses spectral numbers are real because
by assumption the TERP-structure induces PMHS (see the beginning of the
proof of lemma \ref{lemCorPMHS-PMTS}).

The proof of \cite[theorem 7.20]{He2}, and more precisely the formulas
(7.107), (7.117) and (7.123), provides a family of bases
$\uuuu{\sigma}(r)$ for small
$r>0$ such that the coefficients $a_{\alpha,k}(r)$ in all entries of
$C^\pm$ and $(C^\pm)^{-1}$ are real analytic and of order
$O(|\log r|^N)$ for some $N>0$.

Also $\bigoplus_{i=1}^l e^{-u_i/(z\cdot r)}\otimes
\pi_r^*(H_i,\nnn_i,H'_{i,\dR},P_i)$ is a polarized pure TERP-structure
for small $r>0$. A basis is $\uuuu{\sigma}(r)\cdot R(r)$ where
\begin{eqnarray}\label{9.7}
R(r):=\textup{diag}
(\exp(-\frac{u_1}{z\cdot r}-\frac{\oooo{u_1}\cdot z}{r})\cdot {\bf 1}_{n_1},
.., \exp(-\frac{u_l}{z\cdot r}-\frac{\oooo{u_l}\cdot z}{r})\cdot
{\bf 1}_{n_l}).
\end{eqnarray}
The identities
$R(r) = \gamma^*(\oooo{R(r)}) = R(r)^{tr} = R(r,-z)^{-1}$ show the
second equality in
\begin{eqnarray}\label{9.8}
z^{-w}\cdot P((\uuuu{\sigma}(r)\cdot R(r))^{tr},
\tau(\uuuu{\sigma}(r)\cdot R(r)))
= R(r)^{tr}\cdot P_{mat}(r) \cdot \gamma^*(\oooo{R(r,-z)})
= P_{mat}(r).
\end{eqnarray}

We want to show that $\pi_r^*(H,\nnn,H'_\dR,P)$ is a polarized
pure TERP-structure for small $r>0$.
Suppose for a moment that this holds and that
$\uuuu{\omega}=(\omega_1,...,\omega_n)$ is a basis of
$\Gamma(\dP^1,\widehat{\pi_r^*\cH})$. Then
\begin{eqnarray}\label{9.9}
\uuuu{\omega} = \Psi_\pm^{-1}(\uuuu{\sigma}\cdot R)\cdot A^\pm
= \uuuu{e}^\pm\cdot C^\pm\cdot R\cdot A^\pm
\end{eqnarray}
on $\widehat I_\pm(\xi)$, where $A^\pm$ are matrices with entries in
$\Gamma([0,\infty]\times I_\pm(\xi),\widehat\cA_\pm)$.
Formula \eqref{9.1}, that is,
$(\uuuu{e}^-)_{|I(\xi)} = (\uuuu{e}^+)_{|I(\xi)}\cdot T$,
shows on the small sector $\widehat I(\xi)$
\begin{eqnarray}\label{9.10}
\uuuu{\omega}_{|\widehat I(\xi)}
&=& \left( \uuuu{e}^+\cdot C^+\cdot R\cdot A^+\right)_{|\widehat I(\xi)}\\
&=& \left( \uuuu{e}^-\cdot C^-\cdot R\cdot A^-\right)_{|\widehat I(\xi)}
= \left( \uuuu{e}^+\cdot T\cdot C^-\cdot R\cdot A^-
\right)_{|\widehat I(\xi)}.\nonumber
\end{eqnarray}
\eqref{9.2} and  \eqref{9.6} give $C^+_{|\widehat I(\xi)}
=T^{model}\cdot C^-_{|\widehat I(\xi)}$.
Therefore
\begin{eqnarray}\label{9.11}
A^+_{|\widehat I(\xi)}\cdot (A^-_{|\widehat I(\xi)})^{-1}
= \left(R^{-1}\cdot (C^-)^{-1}\cdot T^{Stokes}\cdot C^-\cdot R
\right)_{|\widehat I(\xi)}.
\end{eqnarray}
Analogously on $-\widehat I(\xi)=\widehat I(-\xi)$
\begin{eqnarray}\label{9.12}
A^+_{|\widehat I(-\xi)}\cdot (A^-_{|\widehat I(-\xi)})^{-1}
= \left(R^{-1}\cdot (C^+)^{-1}\cdot (-1)^w(T^{Stokes})^{tr}\cdot C^+\cdot R
\right)_{|\widehat I(-\xi)}.
\end{eqnarray}

Obviously, the bundle $\widehat{\pi_r^*\cH}$ is trivial iff there are
invertible matrices $A^\pm$ with entries in
$\Gamma([0,\infty]\times I_\pm(\xi),\widehat\cA_\pm)$
and satisfying \eqref{9.11} and \eqref{9.12}.
The requirement $A^\pm(0)={\bf 1}_n$ makes them unique.
This is a Riemann boundary value problem. We will argue that the matrices
on the right hand side of \eqref{9.11} and \eqref{9.12} are close to
${\bf 1}_n$ for small $r>0$ and that therefore the problem has a solution.

The matrix $T^{Stokes}=(T^{model})^{-1}\cdot T$ is block upper triangular
with diagonal blocks $T_{jj}^{Stokes}={\bf 1}_{n_j}$
for $j=1,..,l$. The matrices $C^\pm$ are block diagonal.
Therefore the matrix on the right hand side of \eqref{9.11} is also
block upper triangular, and the $(j,k)$-block for $j\leq k$ is
\begin{eqnarray}\label{9.13}
(C^-_j)^{-1}\cdot T_{jk}^{Stokes}\cdot C^-_k\cdot
\exp\left(
\frac{u_j-u_k}{z\cdot r}+ \frac{(\oooo{u_j}-\oooo{u_k})\cdot z}{r}\right).
\end{eqnarray}
The diagonal block for $j=k$ is ${\bf 1}_{n_j}$.

The functions $\left(z\mapsto  \exp\left(
\frac{u_j-u_k}{z\cdot r}+ \frac{(\oooo{u_j}-\oooo{u_k})\cdot z}{r}\right)
\right)$ on $\widehat I(\xi)$ for $j<k$ have asymptotic developments equal to zero
at the origin and at infinity.
Remark that $\Re(\frac{u_j-u_k}{z})<0$ and
$\Re((\oooo{u_j}-\oooo{u_k})\cdot z)<0$ for $z\in \widehat I(\xi)$.
Therefore with $r\to 0$ the functions and all their $z$-derivatives
tend to 0 pointwise.
The moderate behavior in $z$ and $r$ of the entries of
$C^\pm$ and $(C^\pm)^{-1}$ had been discussed above.

We obtain: The restrictions to $\oooo{\xi\cdot \dR_{>0}}\subset \dP^1$
(here $\oooo -$ denotes closure)
of all the entries of a $(j,k)$-block with $j<k$ are real analytic on
$\xi\cdot \dR_{>0}$ and $C^\infty$ at $0$ and $\infty$;
if $r\to 0$, they and all their derivatives tend uniformly to 0.

The same analysis holds for the nondiagonal entries of the right hand
side of \eqref{9.12}.
Here one remarks that this matrix is lower triangular and that
$z\in \widehat I(-\xi)$. One considers the restrictions to
$\oooo{\xi\cdot \dR_{<0}}$.

By a M\"obius transformation one can map $\oooo{\xi\cdot \dR}\subset \dP^1$
to $S^1$.
The Birkhoff decomposition in \cite[(8.1.2)]{PS} gives a certain unique
decomposition for all $C^\infty$ loops $S^1\to GL(n,\dC)$
in an open dense set of the loop space $LGL(n,\dC)$,
which contains 0; and this decomposition depends smoothly on the loop.

In our case it shows the existence and uniqueness of matrices
$A^\pm$ with the following properties:
$A^\pm$ is continuous and invertible on the set
$\oooo{\{z\in \dC^*\ |\ \pm\arg\frac{z}{\xi}<0\}}$;
it is holomorphic on $\{z\in \dC^*\ |\ \pm\arg\frac{z}{\xi}<0\}$;
it satisfies $A^\pm(0)={\bf 1}_n$;
for $z\in\oooo{\xi\cdot \dR}$ \eqref{9.11} and \eqref{9.12} hold.
Furthermore, for $r\to 0$ the restrictions $A^\pm_{|\oooo{\xi\cdot \dR}}$
tend uniformly to ${\bf 1}_n$.

The right hand sides of \eqref{9.11} and \eqref{9.12}
are holomorphic in $\pm \widehat I(\xi)$.
With the theorem of Morera one sees that $A^\pm$ extend to holomorphic
matrices on $\widehat I_\pm(\xi)$ (this argument is taken from
\cite[page 296]{Do}).

We still have to show that $(A^\pm)_0\in \Gamma(I_\pm(\xi),GL(\cA))$
and $(A^\pm)_\infty\in \Gamma(I_\pm(\xi),GL(\gamma^*(\oooo\cA)))$.
But any $\cO_{\dC,0}$-basis of $\cH_0$
gives rise to matrices $\www A^\pm\in \Gamma(I_\pm(\xi),GL(\cA))$ which also
satisfy \eqref{9.11} and \eqref{9.12} for
$z\in \widehat I(\xi)$ close to $0$.
Therefore the matrices $(\www A^\pm)^{-1}\cdot A^\pm$ glue to a matrix
which is continuous and invertible close to $0$ and holomorphic
outside $0$. Thus it is holomorphic at 0.
Therefore $(A^\pm)_0\in \Gamma(I_\pm(\xi),GL(\cA))$.
The same applies at $\infty$.
The Riemann boundary value problem is solved for $r\ll 1$,
and $\pi_r^*(H,\nnn,H'_\dR,P)$ is a pure TERP-structure.

It remains to show that it is polarized.
Let now $\uuuu\omega$ be a global basis as above. The hermitian pairing
is given by the hermitian and $z$-independent matrix
\begin{eqnarray}\label{9.14}
z^{-w}\cdot P(\uuuu\omega^{tr},\uuuu{\omega})
= (A^+(z))^{tr}\cdot P_{mat}(r) \cdot \oooo{A^-(-\frac{1}{\oooo z})}
= P_{mat}(r)\cdot \oooo{A^-(\infty)}.
\end{eqnarray}
The matrix $P_{mat}(r)$ is hermitian and positive definite.

Now one way to conclude is to go through the construction simultaneously for
all Stokes matrices $T(t):=(1-t)\cdot T^{model}+t\cdot T$ with
$t\in[0,1]$. For sufficiently small $r$ it works, and the matrix
$\oooo{A^-(\infty)}(t) $ depends continuously on $t$, with
$\oooo{A^-(\infty)}(t=0)={\bf 1}_n$.
Therefore for all $t\in [0,1]$ the hermitian matrix
$P_{mat}(r)\cdot \oooo{A^-(\infty)}(t)$ is positive definite.

\section{Semi-simple case and ADE-singularities}\label{c10}
\setcounter{equation}{0}

A TERP-structure is called semi-simple if the eigenvalues of the
pole part are all different. Such TERP-structures automatically
do not require a ramification (e.g. \cite[II.5.7]{Sa4}).
Furthermore, a semi-simple TERP-structure is determined by very elementary
data, as is shown in the following lemma.
Let us call two matrices $T$ and $T'$ in $M(n\times n,\dC )$
sign equivalent if there is a matrix
$B=\textup{diag}(\pm 1,...,\pm 1)$ such that $BTB=T'$.
\begin{lemma}\label{t10.1}
\begin{enumerate}
\item
Fix a weight $w\in \dZ$,
$n$ different values $u_1,...,u_n\in\dC$ and $\xi\in S^1$
with $\Re(\frac{u_i-u_j}{\xi})<0$ for $i<j$.
There is a natural 1-1 correspondence between
the set of semi-simple TEP-structures (i.e. no real structure)
of weight $w$ with pole part having eigenvalues $u_1,...,u_n$,
and the set of sign equivalence classes of upper triangular matrices
$T\in M(n\times n,\dC)$ with diagonal entries equal to $1$.
The matrices $T$ are called Stokes matrices of the TEP-structure.
\item
The correspondence in 1. restricts to a correspondence between
semi-simple mixed TERP-structures and sign equivalence classes
of matrices with real entries.
\item
Given a semi-simple TERP-structure of weight $w$ with
eigenvalues $u_1,...,u_n\in\dC$ and $\xi\in S^1$ as above, the corresponding Stokes matrix
$T$ is constructed as follows.
We use the notations and results of chapter \ref{c8},
in particular, lemma \ref{lemma-Stokes1}.
The pieces $\cL_{\pm,j}$ of the canonical (after the choice of
$\xi\in S^1$) splittings $\cL_{|I_\pm(\xi)}=\bigoplus_{j=1}^n\cL_{\pm,j}$
have rank $1$. There are bases $\uuuu{e}^\pm=(e_1^\pm,...,e_n^\pm)$
of $\cL_{|I_\pm(\xi)}$ unique up to the common signs of the pairs
$e_j^+$ and $e_j^-$ with the following properties:
they are compatible with the splittings, they satisfy
$P((\uuuu{e}^+)^{tr},\uuuu{e}^-)=(-i)^w\cdot{\bf 1}_n$ and
$(\uuuu{e}^-)_\xi = (\uuuu{e}^+)_\xi\cdot T$ with $T_{jj}=1$ for all $j$.
Then $T$ is a Stokes matrix.
\end{enumerate}
\end{lemma}
\begin{proof}
First we prove the last part:
If $\uuuu{e}^\pm$ are any bases of $\cL_{|I_\pm(\xi)}$
which are compatible with the splittings, then $T$ is upper triangular
with $\det T\neq 0$ by lemma \ref{lemma-Stokes1},
and $P((\uuuu{e}^+)^{tr},\uuuu{e}^-)$ is diagonal and invertible
by lemma \ref{lemma-Stokes2}.
The constraints $T_{jj}=1$ and $P(e_j^+,e_j^-)=(-i)^w$ determine
$e_j^+$ and $e_j^-$ up to a common sign.
This shows part 3.

Now suppose that the semi-simple TEP-structure is a mixed TERP-structure.
Real structure and Stokes structure are compatible, thus there exist
$\lambda_j^\pm\in S^1$ such that
$(\lambda_1^\pm e_1^\pm,...,\lambda_n^\pm e_n^\pm)$ are real bases of
$\cL_{|I_\pm(\xi)}$.
Because of $T_{jj}=1$ one can choose $\lambda_j^-=\lambda_j^+$.
Then $(\lambda_j^+)^2\cdot (-i)^w =
P(\lambda_j^+e_j^+,\lambda_j^+e_j^-)\in i^w\cdot \dR$
shows $\lambda_j^+\in\{\pm 1\}$, so $\uuuu{e}^\pm$
are real bases and $T$ has real entries.
This proves one direction in 2.

For the other direction and for the first part, we start from
$u_1,...,u_n,\xi$ and $T$ as in 1.
Let us first construct the topological data $(H',\nnn,P)$,
and $H'_\dR$ in case $T$ is real.
Let $\cL_\pm$ be local systems of rank $n$ on $I_\pm(\xi)$
with bases $\uuuu{e}^\pm=(e_1^\pm,...,e_n^\pm)$.
They are glued to a local system $\cL$ on $S^1$ by
$(\uuuu{e}^-)_\xi = (\uuuu{e}^+)_\xi\cdot T$ and
$(\uuuu{e}^-)_{-\xi} = (\uuuu{e}^+)_{-\xi}\cdot (-1)^w T^{tr}$.
If $T$ is real then the bases $\uuuu{e}^\pm$ induce a real
structure on $\cL$.
In any case the formulas
$P((\uuuu{e}^+)^{tr},\uuuu{e}^-):=(-i)^w\cdot{\bf 1}_n$ and
$P((\uuuu{e}^-)^{tr},\uuuu{e}^+):=i^w\cdot{\bf 1}_n$
give a well-defined flat $(-1)^w$-symmetric nondegenerate pairing
on opposite fibers.
This yields the topological data
$(H',\nnn,H'_\dR\textup{ if }T\textup{ is real},P)$.

In the special case $T={\bf 1}_n$ they decompose as
$\bigoplus_{j=1}^n (H'_j,\nnn_j,H'_{j,\dR},P_j)$.
Then we write $\uuuu{f}^\pm$ instead of $\uuuu{e}^\pm$, and
$\cL_j$ denotes the local system on $S^1$ generated by $f_j^\pm$.
Put $\cO(H_j)=\cH_j:=\cO_{\dC}\cdot z^{w/2}\cdot f_j^\pm$,
then $(H_j,\nnn_j,H'_{j,\dR},P_j)$
is a regular singular TERP-structure.
Using \eqref{2.3}, \eqref{5.18} and \eqref{5.19} one checks that it induces
a PHS of weight $w$ if $w$ is even and of weight $w-1$ if $w$ is odd.
Therefore $\bigoplus_{j=1}^n e^{-u_j/z} \otimes (H_j,\nnn_j,H'_{j,\dR},P_j)$
is a mixed TERP-structure of weight $w$ with Stokes matrix $T={\bf 1}_n$.

It remains to construct, for arbitrary $T$, an extension $H\in\mathit{VB}_\dC$ of $H'$
such that $(H,\nnn,P)$ is a TEP-structure formally isomorphic to
$\bigoplus_{j=1}^n e^{-u_j/z} \otimes (H_j,\nnn_j,P_j)$.
The key step is the construction of two invertible matrices
$A^\pm\in \Gamma(I_\pm(\xi),\mathit{GL}(\cA))$ satisfying
\begin{eqnarray}\label{10.1}
(\uuuu{e}^+\cdot z^{w/2}\cdot R_0\cdot A^+)_{|\widehat I(\pm\xi)} =
(\uuuu{e}^-\cdot z^{w/2}\cdot R_0\cdot A^-)_{|\widehat I(\pm\xi)}.
\end{eqnarray}
Here $R_0:=\textup{diag} (e^{-u_1/z},...,e^{-u_n/z})$.
Then $\uuuu{\omega}:=\uuuu{e}^\pm\cdot z^{w/2}\cdot R_0\cdot A^\pm$
defines an extension $H\in\mathit{VB}_\dC$ of $H'$ by
$\cH=\bigoplus_{j=1}^n \cO_{\dC}\cdot\omega_j$
with all desired properties, namely:
It is a semi-simple TEP-structure.
The isomorphisms $\Psi_\pm:\cL_{|I_\pm(\xi)} \to
\bigoplus_{j=1}^n \cL_{j|I_\pm(\xi)}$, $e_j^\pm\mapsto f_j^\pm$,
extend to isomorphisms as in \eqref{8.3} and \eqref{8.4}.
The pairing $P$ satisfies
\begin{eqnarray}\label{10.2}
P(\uuuu{\omega}^{tr}(z),\uuuu{\omega}(-z)) &=&
z^w\cdot A^\pm(z)^{tr}\cdot A^\mp(-z) \hspace*{1cm} \textup{ if }
z\in\widehat I_\pm(\xi),
\end{eqnarray}
the entries of this matrix are elements in
$z^w\Gamma(S^1,\cA)=z^w\Cz$, and the matrix
$z^{-w}\cdot P(\uuuu{\omega}^{tr}(z),\uuuu{\omega}(-z))$
is invertible at $0$.

The two conditions on $A^\pm$ in \eqref{10.2} are equivalent to the
two conditions
\begin{eqnarray}\label{10.3}
(A^+\cdot (A^-)^{tr})_{|\widehat I(\xi)} &=&
(R_0^{-1}\cdot T\cdot R_0)_{|\widehat I(\xi)},\\
(A^+\cdot (A^-)^{tr})_{|\widehat I(-\xi)}&=&
(R_0^{-1}\cdot (-1)^wT^{tr}\cdot R_0)_{|\widehat I(-\xi)}.\label{10.4}
\end{eqnarray}
As in the proof of theorem \ref{main-theorem} 2. one checks that both
matrices on the right hand side have the asymptotic development ${\bf 1}_n$.
\cite[proposition A.1]{Mal1} applies and gives the existence
of $A^\pm$.
By construction, $T$ is a Stokes matrix of the TEP-structure
$(H,\nnn,P)$.
This gives the correspondences in 1. and 2.
\end{proof}
\textbf{Remarks:}
\begin{enumerate}
\item
In many interesting cases the Stokes matrix has actually entries in $\dZ$.
This holds in singularity theory for the TERP-structures defined by
function germs (e.g. \cite{Ph2}\cite[ch. 8]{He2})
and by tame functions \cite{Sa2}\cite{DS}\cite{Sa8},
in quantum cohomology at least for the mirror partners of tame functions,
and also in the massive supersymmetric field theories considered in
\cite{CV1}\cite{CV2}.
\item
It is a major point in \cite{CV1}\cite{CV2} that for
field theories considered therein, the quite elementary data
$(u_1,...,u_n,\xi,T)$ in lemma \ref{t10.1} 2. determine a mixed semi-simple
TERP-structure and thus allow to recover most of the geometry of the
field theory.
\item
Contrary to PMHS, a semi-simple mixed TERP-structure has at most
one compatible lattice in $H'_\dR$,
and that exists precisely iff the Stokes matrix has entries in $\dZ$.
\item
One might ask which of the data $(u_1,...,u_n,\xi,T)$ in lemma
\ref{t10.1} 2.
give rise to a polarized pure TERP-structure.
\cite[proposition 2.2]{Du} resp. theorem \ref{main-theorem} 2.
say that it is sufficient to have that all differences
$|u_i-u_j|$ ($i\neq j$)
are sufficiently large.
The following conjecture proposes another partial answer.
\end{enumerate}

\begin{conjecture}\label{t10.2}
Fix $u_1,...,u_n\in\dC$, $\xi\in S^1$ with $\Re(\frac{u_i-u_j}{\xi})<0$
for $i<j$ and $T\in M(n\times n,\dR)$ upper triangular with $T_{ii}=1$
and such that $T+T^{tr}$ is positive definite.
Then the corresponding mixed TERP-structure of weight $w$
is pure and polarized.
Its spectral numbers at infinity are in the interval
$(\frac{w-1}{2},\frac{w+1}{2})$.
\end{conjecture}

\textbf{Remarks:}
Actually, the matrix $(-1)^wT^{-1}T^{tr}$ gives the monodromy
with respect to the basis $\uuuu{e}^-$ in lemma \ref{t10.1}.
This matrix leaves invariant the pairing, which is given by the matrix
$T+T^{tr}$.
Therefore the monodromy is semi-simple with eigenvalues in
$S^1-\{(-1)^{w+1}\}$.

If the conjecture is true then it applies to all TERP-structures
in a family $\bigcup_{r>0}\pi^*_{r^{-1}}(H,\nnn,H'_\dR,P)$
because they have the same Stokes data and the eigenvalues
are given by $r\cdot u_1,...,r\cdot u_n$.
Then the TERP-structure $(H,\nnn,H'_\dR,P)$ induces a Sabbah orbit,
theorem \ref{theoMainResultSabbah} applies and yields
a sum of  pure PHS of weight $w$ and $w-1$.

\begin{theorem}\label{t10.3}
The conjecture is true if $T$ arises as follows.
Consider a root system of type ADE, its root lattice $L\cong \dZ^n$
with pairing $(-,-)$ and any basis of $L$ consisting of roots
$\beta_1,...,\beta_n$ such that $s_{\beta_1}\circ ...\circ s_{\beta_n}$
is a Coxeter element. Let $T$ be the upper-triangular matrix
defined by $T_{ii}:=1$ and
$T_{ij}:=(\beta_i,\beta_j)$ for $i<j$.
\end{theorem}

\begin{proof}
The proof consists in putting together different results
on the ADE-singularities.
We proceed in eight steps.

\medskip
{\bf Step 1.}
There is a standard universal unfolding
$F(x,t)=f_t(x)=f_0(x)+\sum_{i=1}^nt_im_i(x)$
on $\dC^w\times \dC^n$ of an ADE-singularity $f_0$
such that $F(x,t)$ is weighted homogeneous in $x$ and $t$
with weighted degrees $\deg F(x,t)=\deg(f_0)=1$ and $\deg m_i<1$
(e.g. \cite[8.4]{Ar}\cite[(2.2)]{Lo}).
The parameter space is $M=\dC^n$.
Any function $f_t:\dC^w\to \dC$ is tame
(see the next chapter for this notion).

\medskip
{\bf Step 2.}
Fix an arbitrary set of $n$ different numbers $u_1,...,u_n\in\dC$.
By \cite[(2.4)]{Lo} the number of parameters $t\in \dC^n$ such that
the critical values of $f_t$ are $u_1,...,u_n$ is
$(n+1)^{n-1}$ for $A_n$, $2(n-1)^n$ for $D_n$, $2^9\cdot3^4$ for $E_6$,
$2\cdot 3^{12}$ for $E_7$ and $2\cdot 3^5\cdot 5^7$ for $E_8$.

\medskip
{\bf Step 3.}
By \cite{Ph4}\cite[ch. 8]{He2}
the oscillating integrals of any function $f_t$ in the universal unfolding
induce a TERP-structure of weight $w$.
It does not require a ramification,
and the regular singular pieces are essentially
the Brieskorn lattices of the local singularities.
Because they induce PMHS (\cite{Va1}\cite{SchSt}\cite{SM}, for the
polarization see \cite{He1}),
it is a mixed TERP-structure.
The eigenvalues of the pole part are the critical values of $f_t$.

\medskip
{\bf Step 4.}
By \cite{DS} the TERP-structures for $f_t$, $t\in M$,
fit together to a variation of TERP-structures.
Because of $\deg_w m_i<1=\deg f_0$, the spectral numbers at infinity
(and also the filtration $\www F_{Sab}^\bullet$) are constant.
At $t=0$ they coincide with the usual spectral numbers.
They lie in the interval $(\frac{w-1}{2},\frac{w+1}{2})$.

\medskip
{\bf Step 5.}
By \cite[theorem 4.9]{Sa8} any tame function on an affine manifold
gives rise to a polarized pure TERP-structure via its
oscillating integrals.
This applies to the TERP-structure of $f_t$.

\medskip
{\bf Step 6.}
Fix $n$ different numbers $u_1,...,u_n\in\dC$ and $\xi\in S^1$ with
$\Re(\frac{u_i-u_j}{\xi})<0$ for $i<j$.
Choose a function $f_t$ with critical values $u_1,...,u_n$.

\smallskip
{\bf Claim:} {\it A Stokes matrix $T$ (of the sign equivalence class
in lemma \ref{t10.1}) of the TERP-structure of $f_t$ takes the
following form:
$(-1)^{(w-1)(w-2)/2}\cdot (T-{\bf 1}_n)$ is the strictly upper triangular
part of the intersection matrix of a certain distinguished basis
$\uuuu{\delta}$ of the Milnor lattice
$H_{w-1}(f_t^{-1}((-i)\xi\cdot r),\dZ)$ for some $r\gg 0$.
The distinguished basis $\uuuu{\delta}$ is constructed below.}

\smallskip
For the construction of $\uuuu{\delta}$, we choose $r\gg \max|u_i|$.
We connect $u_i$ with $u_i+(-i)\xi\cdot r$ by a straight line and
$u_i+(-i)\xi\cdot r$ with $(-i)\xi\cdot r$ by a straight line.
Along this path from $u_i$ to $(-i)\xi\cdot r$ the vanishing cycle
of the simple critical point of $f_t$ with value
$u_i$ is shifted to the Milnor fiber
$H_{w-1}(f_t^{-1}((-i)\xi\cdot r),\dZ)$.
This gives a distinguished basis
$\uuuu{\delta}=(\delta_1,...,\delta_n)$ of
$H_{w-1}(f_t^{-1}((-i)\xi\cdot r),\dZ)$ \cite[pages 14 and 31]{AGV}.
It is unique up to the signs $\pm \delta_i$.

For the proof of the claim one needs several facts and formulas:
\begin{list}{}{}
\item[(i)]
The bundle $H'$ of the TERP-structure is the bundle dual to the bundle
of Lefschetz thimbles, and
$P=(-1)^{(w-1)w/2}\cdot \frac{1}{(2\pi i)^w}\cdot P_{Lef}^*$
where $P_{Lef}$ is the intersection form for Lefschetz thimbles,
see \cite[ch. 8]{He2}.
\item[(ii)]
The Stokes structure is determined by the Lefschetz thimbles along
$\bigcup_{i=1}^n (u_i+(-i)\xi\cdot \dR_{>0})$ and the Lefschetz thimbles
along $\bigcup_{i=1}^n (u_i+i\xi\cdot \dR_{>0})$.
\item[(iii)]
One needs the precise relations between Lefschetz thimbles,
vanishing cycles, their intersection forms and the Seifert form,
see \cite[ch. 2]{AGV}.
\end{list}
We omit the details of the proof of the claim.
The claim is consistent with the conventions in \cite{AGV}.
The Picard-Lefschetz transformations $s_{\delta_i}=s_{-\delta_i}$
satisfy $s_{\delta_1}\circ ...\circ s_{\delta_n}=\textup{ monodromy}$.

\medskip
{\bf Step 7.}
Fix $u_1,...,u_n\in\dC$, $\xi\in S^1$ and $r\gg 0$ as in step 6
and denote
$$
\textup{Par}:=\{t\in\dC^n\ |\ f_t\textup{ has critical values }
u_1,...,u_n\}.
$$
The Milnor lattices $H_{w-1}(f_t^{-1}((-i)\xi\cdot r),\dZ)$
for $f_t\in\textup{Par}$ are canonically isomorphic.

Now suppose that $w$ is odd.
It is well known (e.g. \cite{Ar}) that in that case an isomorphism
$H_{w-1}(f_t^{-1}((-i)\xi\cdot r),\dZ)\stackrel{\cong}{\to} L$ exists
which maps the intersection form to $(-1)^{(w-1)/2}(-,-)$
and the monodromy to a Coxeter element $c$ in the Weyl group of $L$.

By \cite[\S 2]{De2} this isomorphism identifies the set of all
distinguished bases of the Milnor lattice
(definition in \cite[pages 14 and 31]{AGV})
with the set
$\cB:=\{(\beta_1,...,\beta_n)\ |\ \beta_i\in L,(\beta_i,\beta_i)=2,
s_{\beta_1}\circ ...\circ s_{\beta_n}=c\}$,
and the natural braid group action on this set is transitive.

Define an equivalence relation on $\cB$ by
$\uuuu{\beta}\sim \uuuu{\beta}' \iff
(\beta_1,...,\beta_n)=(\pm \beta_1',...,\pm\beta_n')$.
Denote by $\uuuu{\delta}(t)\in \cB$ the distinguished basis of $f_t$,
$t\in\textup{Par}$, which was constructed in step 6.
In \cite[\S 3]{De2} the formula $|\textup{Par}|=|\cB/\sim|$ is proved.
This and the discussion in \cite[\S 3]{Lo} show that the map
\begin{eqnarray}\label{10.5}
\textup{Par}\to \cB/\sim,\quad t\mapsto [\uuuu{\delta}(t)]
\end{eqnarray}
is a bijection.

\medskip
{\bf Step 8.}
The steps 7 and 6 show that any matrix $T$ in theorem \ref{t10.3}
is realized as Stokes matrix of the TERP-structure of a function $f_t$
with $t\in\textup{Par}$.
The steps 1 to 7 show theorem \ref{t10.3} for odd $w$.
The case of the even number $w+1$ can be treated by twisting the TERP-structure with
$z^{1/2}$ or by using the formulas in \cite[2.8]{AGV}
for the intersection forms of $f_t(x)$ and $f_t(x)+x_{w+1}^2$
with respect to distinguished bases.
\end{proof}

\textbf{Remarks:}
\begin{enumerate}
\item
The number of parameters in $\textup{Par}$
resp. of equivalence classes of distinguished
bases with fixed sign equivalence class of Stokes matrices is larger
than one. It is
$n+1$ for $A_n$, $2(n-1)$ for $D_n$,
$12$ for $E_6$, $9$ for $E_7$ and $15$ for $E_8$.
This follows from \cite{Vo}\cite{Klu}.
There the numbers of distinguished bases
with fixed Coxeter-Dynkin diagrams are listed.
It turns out to be twice the number above:
The Stokes matrices correspond to the Coxeter-Dynkin diagrams
of the distinguished bases. In one equivalence class
in $\cB/\sim$ exactly two distinguished bases have the same Coxeter-Dynkin
diagram, $(\beta_1,...,\beta_n)$ and $(-\beta_1,...,-\beta_n)$.
This follows from the connectedness of the Coxeter-Dynkin diagrams.

The number above is also the order of the automorphism group
$R_{f_0}$ studied in \cite[13.2]{He1}.

\item
The results on Landau-Ginzburg models in \cite{CV1}\cite{CV2}
are closely related to the crucial step 5, \cite[theorem 4.9]{Sa8}.
The ADE case follows also from \cite{CV1}\cite{CV2}.
But the proofs are completely different.
\end{enumerate}

\begin{proposition}\label{t10.4}
Conjecture \ref{t10.2} is true in the rank two case.
\end{proposition}

\begin{proof}
We only sketch the proof.
In \cite{Du} and \cite{CV1}\cite{CV2} it is proved that there is a
correspondence between semi-simple TERP-structures $(H,\nnn,H'_\dR,P)$
such that $\pi_{r^{-1}}^*(H,\nnn,H'_\dR,P)$ is pure polarized
for all $r>0$ and real smooth solutions on $(0,\infty)$ of the
sinh-Gordon equation
$(\partial_r^2+\frac{1}{r}\partial_r)u(r) = \sinh u(r)$.
It also follows implicitly from \cite{IN}.
The hermitian metric with respect to a certain basis is given by the
matrix
$$\begin{pmatrix} \cosh(\frac{u(r)}{2})& -i\sinh(\frac{u(r)}{2})\\
i\sinh(\frac{u(r)}{2})& \cosh(\frac{u(r)}{2})\end{pmatrix}.$$

In \cite{MTW} a family depending on one real parameter of such solutions has been
studied. In \cite[ch. 11]{IN} it is shown that the Stokes matrices
of these solutions are
$T=\begin{pmatrix} 1&t\\0&1\end{pmatrix}$
with $t\in (-2,2)$;
it is the case $A=B$ and $Q=\begin{pmatrix} 0 & {*} \\ {*} & 0\end{pmatrix}$
in \cite[ch. 11]{IN}.
This shows the first claim of the conjecture.

The behavior of the solutions for $r\to 0$
\cite{MTW}\cite[(11.2)]{IN} shows that the spectral numbers at infinity
are in the interval $(-\frac{1}{2},\frac{1}{2})$.
This shows the conjecture for $w=0$.
For other $w$ one twists the TERP-structure with $z^{w/2}$.
\end{proof}

\textbf{Remarks:}
Also singular solutions of the sinh-Gordon equation correspond to
families $\pi_{r^{-1}}^*(H,\nnn,H'_\dR,P)$ of semi-simple TERP-structures
of rank 2. At a singularity $r_k>0$ a solution has the asymptotic
form $u(r)= -2\log(r-r_k)-O(r-r_k)$ \cite[(11.6)]{IN}.
On one side of $r_k$ it is real, on the other side it takes values
in $2\pi i+\dR$.
Then $\pi^*_{r_k^{-1}}(H,\nnn,H'_\dR,P)$ is not pure,
and the hermitian form defined by $\pi_{r^{-1}}^*(H,\nnn,H'_\dR,P)$ is positive definite for $r$
on one side and negative definite for $r$ on the other side.

In \cite{MTW} also the solutions corresponding to Stokes matrices
$T=\begin{pmatrix} 1&t\\0&1\end{pmatrix}$ with $t\in \dR-[-2,2]$
are studied. By theorem \ref{main-theorem} they are smooth for $r\to\infty$.
But they have infinitely many singularities for $r\to 0$.
Their distribution is studied in \cite[page 1090]{MTW}.
Conjecture \ref{main-conjecture} predicts that all solutions
with Stokes matrix
$T=\begin{pmatrix} 1&t\\0&1\end{pmatrix}$ with $t\in \dC \backslash \dR$
have singularities.
In \cite[ch. 11]{IN} only solutions which are smooth for small $r$
with asymptotics \cite[(11.2)]{IN} are studied.
For those the rank two case of the conjecture \ref{main-conjecture} is proved:
Only the solutions in the one parameter family of \cite{MTW}
do not have singularities.
All other solutions have infinitely many singularities
$r_k = \pi(k-\frac{1}{2}) + O(\log k)$, $k\to\infty$ \cite[(11.10)]{IN}.

So in the semi-simple cases of rank two which are studied,
there are infinitely many singularities, i.e. not pure TERP-structures,
if there are any singularities at all.
This is in sharp contrast to the regular singular rank two cases
discussed in \cite[8.3]{He2}, where the TERP structures in the families
$\bigcup_{r>0}\pi_{r^{-1}}(H,\nnn,H'_\dR,P)$
are pure everywhere except for one parameter.

\section{Remarks on applications}\label{c11}
\setcounter{equation}{0}

As we already pointed out at several places,
one of the major sources of examples of TERP-structures is singularity theory.
Oscillating integrals of a holomorphic function with isolated
singularities give rise to a TERP-structure.
Oscillating integrals have been studied since long time
(\cite{Ph4}\cite{AGV} and references therein).
It is well known that they arise by a Fourier transformation from
the Gauss-Manin system and the Brieskorn lattice.

The map $\tau$ extending the bundle $H$ to a bundle
$\widehat H$ on $\dP^1$ was first considered for the TERP-structures
from oscillating integrals in the work of Cecotti and Vafa
\cite{CV1}\cite{CV2} on Landau-Ginzburg models.
These papers were a source of inspiration for \cite{He2} and also for the present article.
A similar construction, but using a hermitian metric from the very beginning instead of
a real structure and a pairing $P$, is present
in the work of Simpson on harmonic bundles
\cite{Si2}\cite{Si4}\cite{Si5}. Recently, this theory was generalized by
Sabbah \cite{Sa7} to what he calls polarizable twistor
$\cD$-modules.

Let us give a very brief summary on the situation in singularity theory
leading to TERP-structures. We consider simultaneously two cases, which are called local and global:

{\bf Local case.} $f:(\dC^w,0)\to (\dC,0)$ is a holomorphic function
germ with an isolated singularity at 0 with Milnor number $n$.

{\bf Global case.} $f:Y\to \dC$ is a function on an affine manifold $Y$
of dimension $w$ such that $f$ is M-tame (definition in \cite{NS}
and \cite[2.a]{DS})
and cohomologically tame (definition in \cite{Sa2}).
We simply call such $f$ tame. Then $f$ has only isolated singularities.
The sum $n$ of the Milnor numbers of the singularities is the
(global) Milnor number of $f$.

In both cases there exists a semi-universal unfolding $F$
(where in the global case ``semi-universal'' refers to
the Kodaira-Spencer map being an isomorphism, see \cite{DS})
with smooth base space $M$, isomorphic to the germ $(\dC^n,0)$.
The deformed functions $F_t$ for $t\in M$ are defined on a some small or large Stein manifold
(in the local case a small ball in $\dC^w$).
In the global case here M-tameness is used \cite{DS}.

\begin{theorem}\label{t11.1}
In both cases one obtains a variation of mixed TERP-structures
$\bigcup_{t\in M} \mathit{TERP}(F_t)$ of rank $n$ on $M$.
The notation $\mathit{TERP}(F_t)=(H(t),\nnn,H'_\dR,P)$ makes sense because
the topological data $(H'=H_{|\dC^*},\nnn,H'_\dR,P)$
are canonically isomorphic for all TERP-structures $\mathit{TERP}(F_t)$.
$H'$ is the bundle dual to the bundle of Lefschetz thimbles,
$\nnn$ is the natural flat connection from shifting Lefschetz thimbles.
The pairing $P$ is defined up to a constant by the intersection form
of Lefschetz thimbles.
The sections in $\cO(H)$ come from the Fourier transform of the
Gauss-Manin system.
\end{theorem}

The precise construction in the local case is described in \cite[8.1]{He2}.
In the global case it is given for $F_0=f$ in \cite{Sa2} and for any
$F_t$ in \cite{DS}. It builds on the work of many people
on the Gauss-Manin system and the Brieskorn lattice.

The TERP-structure $\mathit{TERP}(F_t)$ is a mixed TERP-structure:
By construction, it does not require a ramification and the regular singular
pieces are essentially the local Brieskorn lattices of the singularities
of $F_t$. Compatibility of real structure and Stokes structure is trivial,
because the splittings of the local system comes from topology,
from the Lefschetz thimbles and vanishing cycles
associated to the singularities.
The fact that the regular singular pieces give rise to PMHS is due to
\cite{Va1}\cite{SchSt} (for the polarization see \cite[ch. 10]{He1}).

The Landau-Ginzburg models of Cecotti and Vafa \cite{CV1}\cite{CV2}
involve a lot of additional structure from physics, but the central
objects of study are the TERP-structures of tame functions $f:Y\to\dC$
(or at least a substantial subfamily of them).
From physical considerations Cecotti and Vafa derive the following.
An independent completely different purely mathematical
proof was given recently by Sabbah \cite[theorem 4.9]{Sa8}\cite{Sa7}.

\begin{theorem}\label{t11.2}
The TERP-structure of a tame function $f:Y\to\dC$ on an affine manifold
$Y$ is pure and polarized.
\end{theorem}

In \cite{CV1}\cite{CV2} the resulting positive definite hermitian metric
$h$ is their ground state metric. This fundamental theorem will certainly
play an important role in the future study of tame functions.
It can be considered as the analogue for tame functions of the fact
that the (primitive part of the) cohomology of a compact K\"ahler manifold
carries (polarized) Hodge structures.

If we consider a semi-universal unfolding $F$ of a tame
function $f:Y\to\dC$, then most of the functions $F_t$, $t\in M$,
have to be restricted to a Stein subset of $Y$, as they would have
additional singularities ``from infinity'' on $Y$.
Therefore the theorem applies only to a certain subfamily of all
functions in the semi-universal unfolding, the subfamily along
which the global Milnor number $n$ and the tameness condition are preserved.

Let us discuss some applications of the main results of this paper
(theorems \ref{theoMainResultSabbah} and \ref{main-theorem}, the latter
containing theorem \ref{theoMainResultRegSing}) to TERP-structures coming
from local and global singularities. First we state a simple lemma.
\begin{lemma}\label{t11.3}
Consider a function $F_t$ in a semi-universal unfolding of a function $f$,
in the local or the global case. Then for any $r\in\dC^*$
\begin{eqnarray}\label{11.1}
\mathit{TERP}(r\cdot F_t) = \pi^*_{r^{-1}}(\mathit{TERP}(F_t)).
\end{eqnarray}
Furthermore, the one parameter unfolding
$r\cdot F_t$, $r\in\dC^*$, of $F_t$ is
(for $r$ close to $1$) isomorphic to the
Euler field orbit of $F_t$ in the universal unfolding,
and the pull back of the Euler field is the vector field
$r\partial_r$ on $\dC^*$.
\end{lemma}
\begin{proof}
The first part follows from the formulas for the Fourier transformation
of the Gauss-Manin system. The central point is
$\pi^*_{r^{-1}}(e^{-\theta/z})=e^{-r\cdot \theta/z}$;
here $\theta$ is the variable for the values of $F_t$.
The second part is proved in \cite[lemma 8.6]{He2}.
The central formula is $r\partial_r (r\cdot F_t) = r\cdot F_t$.
\end{proof}

The next result gives two major applications of our correspondence.
It shows nicely that both implications of the correspondence
are of interest.
\begin{corollary}\label{t11.4}
\begin{enumerate}
\item
$\mathit{TERP}(F_t)$ induces a nilpotent orbit.
For $|r|\gg 0$ the TERP-structure $\mathit{TERP}(r\cdot F_t)$ is pure and
polarized.
\item
In the case $F_0=f:Y\to \dC$ tame,
$\mathit{TERP}(f)$ induces a Sabbah orbit.
Sabbah's Hodge filtration makes the tuple
$(\hiii,\hiir,N,S,\www F^\bullet_{Sab})$ (see theorem \ref{theoMainResultSabbah})
into a PMHS of weight $w-1$ resp. $w$ on
$H^\infty_{\neq 1}$ resp. $H^\infty_1$.
\end{enumerate}
\end{corollary}

\begin{proof}
\begin{enumerate}
\item
$\mathit{TERP}(F_t)$ is a mixed TERP-structure by theorem \ref{t11.1}.
The second part of theorem \ref{main-theorem}  shows that it a nilpotent orbit.
Now formula \eqref{11.1} gives the second claim.
\item
$r\cdot f$ is tame for any $r\in\dC^*$. By theorem \ref{t11.2}
$\mathit{TERP}(r\cdot f)$ is pure and polarized for any $r\in\dC^*$.
In particular, $\mathit{TERP}(f)$ induces a nilpotent as well as a Sabbah orbit.
Now theorem \ref{theoMainResultSabbah} applies.
\end{enumerate}
\end{proof}
The first point of the corollary proves the main part of conjecture
8.3 in \cite{He2}. The second part strengthens a former result of Sabbah.
Namely, it was shown in \cite{Sa2} (building on \cite{Sa1}) that
$(\hiii,\hiir,W_\bullet, F^\bullet_{Sab})$ is a mixed Hodge structure,
where $W_\bullet$ is the weight filtration from the nilpotent part $N$
of the monodromy and $F^\bullet_{Sab}$ is the filtration
on $\hiii$ defined in \eqref{7.1}. Polarizations are not considered in these papers.
To obtain a PMHS, we need to work with the twisted filtration
$\www F^\bullet_{Sab} :=G^{-1}(F^\bullet_{Sab})$, but as already said,
it coincides with $F^\bullet_{Sab}$ on the quotient of the weight filtration.

\medskip
{\bf Remarks:}
It is interesting to ponder further on the logical interrelations.
\begin{itemize}
\item
We know already that the TERP-structure $\mathit{TERP}(f)$
for $f:Y\to\dC$ is mixed. But Theorem \ref{t11.2} and
the unproved direction of conjecture \ref{main-conjecture}
would give it again, because theorem \ref{t11.2} applied
to all $r\cdot f$ shows that $\mathit{TERP}(f)$
induces a nilpotent orbit. This would give a new proof
that Stokes structure and real structure are compatible and
that the Brieskorn lattices of the singularities of $f$ induce PMHS.
\item
The other way round, the first part of corollary \ref{t11.4}
applied to tame $f:Y\to\dC$ shows that theorem \ref{t11.2}
is true for $r\cdot f$ with $r\gg 0$.
\end{itemize}
\medskip

We conclude with some remarks on quantum cohomology.
The following ideas and speculations suggest that using mirror symmetry,
theorem \ref{t11.2} it is not so far away from the classical fact that
the cohomology of K\"ahler manifolds carries Hodge structures.

Mirror symmetry predicts that certain tame functions are related
to certain Fano manifolds, more precisely, that
the TEP-structures of the functions are isomorphic
to the TEP-structures in the structure connections
of the quantum cohomology of the Fano manifolds.

By these isomorphisms, the natural real structures on the singularity side
induce real structures on the quantum cohomology side.
However, there seems to be for the moment no intrinsic mathematical
description of this real structure on the quantum cohomology side
(possibly the structures in \cite{CV1}\cite{CV2}) have to be studied
carefully).

If there were such a description, one could hope that the TERP-structures from the
quantum cohomology of any manifold would be pure and polarized.
This would probably only hold for the TERP-structures for
parameters in or close to the small quantum cohomology,
because in some known examples of mirror symmetry small quantum
cohomology corresponds to the tame unfoldings $F_t$ of $F_0=f$.

The classical cohomology of a K\"ahler manifold
endowed with cup product and the sum of Hodge structures
is obtained as semiclassical limit from the quantum
cohomology ring. In some known examples the semiclassical limit
can be considered as the center of a
normal crossing divisor along which the family of TEP-structures
from quantum cohomology has logarithmic poles.
If the quantum multiplication could be used to define a TERP-structures,
one could hope to obtain a variation of
pure polarized TERP-structures on the complement of the divisor.
This may allow to get the classical Hodge structures as a limit
PMTS (or a quotient of it by some filtration, also an additional
twist is possible) using results from
\cite{Mo2} or \cite{Sa6}.

\bibliographystyle{amsalpha}

\begin{thebibliography}{AGZV88}

\bibitem[AGZV88]{AGV}
V.~I. Arnold, S.~M. Gusein-Zade, and A.~N. Varchenko, \emph{Singularities of
  differentiable maps. {V}ol. {II}}, Monographs in Mathematics, vol.~83,
  Birkh\"auser Boston Inc., Boston, MA, 1988, Monodromy and asymptotics of
  integrals, Translated from the Russian by Hugh Porteous, Translation revised
  by the authors and James Montaldi.

\bibitem[Arn72]{Ar}
V.I. Arnold, \emph{{Normal forms for functions near degenerate critical points,
  the Weyl groups of A$_k$, D$_k$, E$_k$ and lagrangian singularities.}},
  Funct. Anal. Appl. \textbf{6} (1972), 254--272 (Russian, English).

\bibitem[CFIV92]{CFIV}
Sergio Cecotti, Paul Fendley, Ken Intriligator, and Cumrun Vafa, \emph{A new
  supersymmetric index}, Nuclear Phys. B \textbf{386} (1992), no.~2, 405--452.

\bibitem[CK82]{CaK1}
Eduardo Cattani and Aroldo Kaplan, \emph{Polarized mixed {H}odge structures and
  the local monodromy of a variation of {H}odge structure}, Invent. Math.
  \textbf{67} (1982), no.~1, 101--115.

\bibitem[CK89]{CaK2}
\bysame, \emph{Degenerating variations of {H}odge structure}, Ast\'erisque
  (1989), no.~179-180, 9, 67--96, Actes du Colloque de Th\'eorie de Hodge
  (Luminy, 1987).

\bibitem[CKS86]{CaKS}
Eduardo Cattani, Aroldo Kaplan, and Wilfried Schmid, \emph{Degeneration of
  {H}odge structures}, Ann. of Math. (2) \textbf{123} (1986), no.~3, 457--535.
  \MR{MR840721 (88a:32029)}

\bibitem[CV91]{CV1}
Sergio Cecotti and Cumrun Vafa, \emph{Topological--anti-topological fusion},
  Nuclear Phys. B \textbf{367} (1991), no.~2, 359--461.

\bibitem[CV93]{CV2}
\bysame, \emph{On classification of {$N=2$} supersymmetric theories}, Comm.
  Math. Phys. \textbf{158} (1993), no.~3, 569--644.

\bibitem[Del]{De2}
Pierre Deligne, \emph{Letter to {L}ooijenga on {M}arch 9, 1974}, Reprinted in
  the diploma thesis of {P}. {K}luitmann: {G}eometrische {B}asen des
  {M}ilnorgitters einer einfach elliptischen {S}ingularität, {B}onn, 1983, pp.
  102--111.

\bibitem[Del71]{De}
\bysame, \emph{Th\'eorie de {H}odge. {II}}, Inst. Hautes \'Etudes Sci. Publ.
  Math. (1971), no.~40, 5--57.

\bibitem[Dim04]{Di}
Alexandru Dimca, \emph{Sheaves in topology}, Universitext, Springer-Verlag,
  Berlin, 2004.

\bibitem[Dou83]{Do}
Adrien Douady, \emph{Probl\`eme de {R}iemann-{H}ilbert. {II}. {S}olution pour
  des points singuliers r\'eels}, Mathematics and physics (Paris, 1979/1982),
  Progr. Math., vol.~37, Birkh\"auser Boston, Boston, MA, 1983, pp.~289--298.

\bibitem[DS03]{DS}
Antoine Douai and Claude Sabbah, \emph{Gauss-{M}anin systems, {B}rieskorn
  lattices and {F}robenius structures. {I}}, Ann. Inst. Fourier (Grenoble)
  \textbf{53} (2003), no.~4, 1055--1116.

\bibitem[Dub93]{Du}
Boris Dubrovin, \emph{Geometry and integrability of topological-antitopological
  fusion}, Comm. Math. Phys. \textbf{152} (1993), no.~3, 539--564.

\bibitem[Her02]{He1}
Claus Hertling, \emph{Frobenius manifolds and moduli spaces for singularities},
  Cambridge Tracts in Mathematics, vol. 151, Cambridge University Press,
  Cambridge, 2002.

\bibitem[Her03]{He2}
\bysame, \emph{{$tt\sp *$} geometry, {F}robenius manifolds, their connections,
  and the construction for singularities}, J. Reine Angew. Math. \textbf{555}
  (2003), 77--161.

\bibitem[IN86]{IN}
Alexander~R. Its and Victor~Yu. Novokshenov, \emph{The isomonodromic
  deformation method in the theory of {P}ainlev\'e equations}, Lecture Notes in
  Mathematics, vol. 1191, Springer-Verlag, Berlin, 1986.

\bibitem[Klu89]{Klu}
P.~Kluitmann, \emph{Addendum zu der {A}rbeit "{A}usgezeichnete {B}asen von
  {M}ilnorgittern einfacher {S}ingularit\"aten" von {E}. {V}oigt}, Abh. Math.
  Sem. Univ. Hamburg \textbf{59} (1989), 123--124.

\bibitem[Loo74]{Lo}
Eduard Looijenga, \emph{The complement of the bifurcation variety of a simple
  singularity}, Invent. Math. \textbf{23} (1974), 105--116.

\bibitem[Mal83]{Mal1}
B.~Malgrange, \emph{La classification des connexions irr\'eguli\`eres \`a une
  variable}, Mathematics and physics (Paris, 1979/1982), Progr. Math., vol.~37,
  Birkh\"auser Boston, Boston, MA, 1983, pp.~381--399.

\bibitem[Moc07]{Mo2}
Takuro Mochizuki, \emph{{Asymptotic behaviour of tame harmonic bundles and an
  application to pure twistor $\mathcal{{D}}$-modules, Part 1}}, Mem. Amer.
  Math. Soc. \textbf{185} (2007), no.~869, xi+324.

\bibitem[MTW77]{MTW}
Barry~M. McCoy, Craig~A. Tracy, and Tai~Tsun Wu, \emph{Painlev\'e functions of
  the third kind}, J. Mathematical Phys. \textbf{18} (1977), no.~5, 1058--1092.

\bibitem[NS99]{NS}
A.~N{\'e}methi and C.~Sabbah, \emph{Semicontinuity of the spectrum at
  infinity}, Abh. Math. Sem. Univ. Hamburg \textbf{69} (1999), 25--35.

\bibitem[Pha83]{Ph2}
Fr{\'e}d{\'e}ric Pham, \emph{Structures de {H}odge mixtes associ\'ees \`a un
  germe de fonction \`a point critique isol\'e}, Analysis and topology on
  singular spaces, II, III (Luminy, 1981), Ast\'erisque, vol. 101, Soc. Math.
  France, Paris, 1983, pp.~268--285.

\bibitem[Pha85]{Ph4}
\bysame, \emph{La descente des cols par les onglets de {L}efschetz, avec vues
  sur {G}auss-{M}anin}, Ast\'erisque (1985), no.~130, 11--47, Differential
  systems and singularities (Luminy, 1983).

\bibitem[PS86]{PS}
Andrew Pressley and Graeme Segal, \emph{Loop groups}, Oxford Mathematical
  Monographs, The Clarendon Press Oxford University Press, New York, 1986,
  Oxford Science Publications.

\bibitem[Sab]{Sa2}
Claude Sabbah, \emph{Hypergeometric period for a tame polynomial}, Preprint
  math.AG/9805077, short version published in: {C}omptes {R}endus de
  l'{A}cad\'emie des {S}ciences. {S}\'erie I. {M}ath\'ematique, vol. 328, no.
  7, 1999.

\bibitem[Sab97]{Sa1}
\bysame, \emph{Monodromy at infinity and {F}ourier transform}, Publ. Res. Inst.
  Math. Sci. \textbf{33} (1997), no.~4, 643--685.

\bibitem[Sab02]{Sa4}
\bysame, \emph{D\'eformations isomonodromiques et vari\'et\'es de {F}robenius},
  Savoirs Actuels, EDP Sciences, Les Ulis, 2002, Math\'ematiques.

\bibitem[Sab04]{Sa7}
\bysame, \emph{The {F}ourier-{L}aplace transform of irreducible regular
  differential systems on the {R}iemann sphere}, Uspekhi Mat. Nauk \textbf{59}
  (2004), no.~6(360), 161--176.

\bibitem[Sab05a]{Sa8}
\bysame, \emph{{F}ourier-{L}aplace transform of a variation of polarized
  complex {H}odge structure.}, Preprint math.AG/0508551, 2005.

\bibitem[Sab05b]{Sa6}
\bysame, \emph{Polarizable twistor {$\mathcal{D}$}-modules}, Ast\'erisque
  (2005), no.~300, vi+208.

\bibitem[Sai89]{SM}
Morihiko Saito, \emph{On the structure of {B}rieskorn lattice}, Ann. Inst.
  Fourier (Grenoble) \textbf{39} (1989), no.~1, 27--72.

\bibitem[Sch73]{Sch}
Wilfried Schmid, \emph{Variation of {H}odge structure: the singularities of the
  period mapping}, Invent. Math. \textbf{22} (1973), 211--319.

\bibitem[Ser66]{Se}
Jean-Pierre Serre, \emph{Prolongement de faisceaux analytiques coh\'erents},
  Ann. Inst. Fourier (Grenoble) \textbf{16} (1966), no.~fasc. 1, 363--374.

\bibitem[Sim88]{Si1}
Carlos~T. Simpson, \emph{Constructing variations of {H}odge structure using
  {Y}ang-{M}ills theory and applications to uniformization}, J. Amer. Math.
  Soc. \textbf{1} (1988), no.~4, 867--918. \MR{MR944577 (90e:58026)}

\bibitem[Sim90]{Si2}
\bysame, \emph{Harmonic bundles on noncompact curves}, J. Amer. Math. Soc.
  \textbf{3} (1990), no.~3, 713--770.

\bibitem[Sim92]{Si4}
\bysame, \emph{Higgs bundles and local systems}, Inst. Hautes \'Etudes Sci.
  Publ. Math. (1992), no.~75, 5--95.

\bibitem[Sim97]{Si5}
\bysame, \emph{Mixed twistor structures}, Preprint math.AG/9705006, 1997.

\bibitem[SS85]{SchSt}
J.~Scherk and J.~H.~M. Steenbrink, \emph{On the mixed {H}odge structure on the
  cohomology of the {M}ilnor fibre}, Math. Ann. \textbf{271} (1985), no.~4,
  641--665.

\bibitem[Var80]{Va1}
A.~N. Varchenko, \emph{Asymptotic behavior of holomorphic forms determines a
  mixed {H}odge structure}, Dokl. Akad. Nauk SSSR \textbf{255} (1980), no.~5,
  1035--1038.

\bibitem[Voi85]{Vo}
E.~Voigt, \emph{Ausgezeichnete {B}asen von {M}ilnorgittern einfacher
  {S}ingularit\"aten}, Abh. Math. Sem. Univ. Hamburg \textbf{55} (1985),
  183--190.

\end{thebibliography}
\providecommand{\bysame}{\leavevmode\hbox to3em{\hrulefill}\thinspace}
\providecommand{\MR}{\relax\ifhmode\unskip\space\fi MR }
\providecommand{\MRhref}[2]{%
  \href{http://www.ams.org/mathscinet-getitem?mr=#1}{#2}
}
\providecommand{\href}[2]{#2}

\vspace*{1cm}

\nd
Lehrstuhl f\"ur Mathematik VI \\
Institut f\"ur Mathematik\\
Universit\"at Mannheim,
A 5, 6 \\
68131 Mannheim\\
Germany

\vspace*{1cm}

\nd
hertling@math.uni-mannheim.de\\
sevenheck@math.uni-mannheim.de

\end{document}